\documentclass[english,3p]{elsarticle}
\usepackage[T1]{fontenc}
\usepackage[utf8]{inputenc}
\usepackage{units}
\usepackage{amsmath}
\usepackage{amsthm}
\usepackage{amssymb}
\usepackage{stmaryrd}
\usepackage{graphicx}

\makeatletter

\providecommand{\tabularnewline}{\\}

\numberwithin{equation}{section}
\numberwithin{figure}{section}

\usepackage{tikz}
\usetikzlibrary{shapes,arrows}

\@ifundefined{showcaptionsetup}{}{%
 \PassOptionsToPackage{caption=false}{subfig}}
\usepackage{subfig}
\makeatother

\usepackage{babel}
\begin{document}

\title{A Least-Squares Formulation of the Moving Discontinuous Galerkin
Finite Element Method with Interface Condition Enforcement}

\author[lcp]{Andrew D. Kercher}
\author[lcp]{Andrew Corrigan}
\address[lcp]{Laboratories for Computational Physics and Fluid Dynamics,  U.S. Naval Research Laboratory, 4555 Overlook Ave SW, Washington, DC 20375}

\begin{abstract}
A least-squares formulation of the Moving Discontinuous Galerkin Finite
Element Method with Interface Condition Enforcement (LS-MDG-ICE) is
presented. This method combines MDG-ICE, which uses a weak formulation
that separately enforces a conservation law and the corresponding
interface condition and treats the discrete geometry as a variable,
with the Discontinuous Petrov-Galerkin (DPG) methodology of Demkowicz
and Gopalakrishnan to systematically generate optimal test functions
from the trial spaces of both the discrete flow field and discrete
geometry. For inviscid flows, LS-MDG-ICE detects and fits a priori
unknown interfaces, including shocks. For convection-dominated diffusion,
LS-MDG-ICE resolves internal layers, e.g., viscous shocks, and boundary
layers using anisotropic curvilinear $r$-adaptivity in which high-order
shape representations are anisotropically adapted to accurately resolve
the flow field. As such, LS-MDG-ICE solutions are oscillation-free,
regardless of the grid resolution and polynomial degree. Finally,
for both linear and nonlinear problems in one dimension, LS-MDG-ICE
is shown to achieve optimal convergence of the $L^{2}$ solution error
with respect to the exact solution when the discrete geometry is fixed
and super-optimal convergence when the discrete geometry is treated
as a variable.
\end{abstract}
\maketitle\global\long\def\middlebar{\,\middle|\,}
\global\long\def\average#1{\left\{  \!\!\left\{  #1\right\}  \!\!\right\}  }
\global\long\def\cof{\operatorname{cof}}
\global\long\def\det{\operatorname{det}}
\global\long\def\adj{\operatorname{adj}}
\global\long\def\expnumber#1#2{{#1}\mathrm{e}{#2}}

\let\svthefootnote\thefootnote\let\thefootnote\relax\footnotetext{\\ \hspace*{120pt}Distribution A. Approved for public release: distribution unlimited.}\addtocounter{footnote}{-1}\let\thefootnote\svthefootnote

\section{Introduction\label{sec:Introduction} }

The Moving Discontinuous Galerkin Finite Element Method with Interface
Condition Enforcement (MDG-ICE) is a high-order, $r$-adaptive method
that is capable of computing flows with a priori unknown interfaces~\citep{Cor17,Cor18}.
The method simultaneously solves for the discrete flow field and discrete
geometry in order to satisfy the governing conservation law and corresponding
interface condition, which are enforced separately within a weak formulation.
As such, MDG-ICE is capable of resolving flows containing a priori
unknown interfaces, including shocks with non-trivial topology, with
high-order accuracy by moving the grid to fit the interfaces, avoiding
low-order errors introduced by artificial stabilization terms such
as shock capturing. The same mechanism that enables MDG-ICE to detect
and fit interfaces, i.e., satisfaction of the weak form, allows it
to automatically resolve sharp, but smooth, gradients using anisotropic
curvilinear $r$-adaptivity. Furthermore, since it is formulated in
terms of a general conservation law, incorporating additional physics
is straightforward and natural.

Recently, the method was applied to compressible viscous flows~\citep{Ker20},
extending the ability of MDG-ICE to move the grid in order to accurately
resolve complex flow features to the case of viscous flows with internal
and boundary layers. In the viscous setting, MDG-ICE separately enforces
a constitutive law and the corresponding interface condition, which
constrains the continuity of the state variable across the interface,
in addition to the conservation law and interface (Rankine-Hugoniot)
condition.

So far, MDG-ICE has been implemented in terms of a \emph{discrete
least squares} (DLS) formulation, cf.~\citep{Kei17}, where discretization
occurs prior to formulation of the least-squares system. In this work,
we derive a least-squares formulation of MDG-ICE by applying the DPG
methodology of Demkowicz and Gopalakrishnan~\citep{Dem10,Dem11,Dem12,Dem15_03,Dem15_20,Dem18,Gop13,Gop14,Car16}
to systematically generate optimal test functions from the trial spaces
of both the discrete flow field and discrete geometry. The proposed
approach, LS-MDG-ICE also extends least-squares finite element methods~\citep{Cai94,Boc98,Boc09}
by introducing the discrete geometry as a variable.

In this work, we apply LS-MDG-ICE to inviscid and viscous flow problems,
many of which were considered in our earlier work~\citep{Cor18,Ker20},
in order to study the ability of the least-squares formulation to
further improve the convergence properties of MDG-ICE. In particular,
we derive optimal test spaces in the context of the MDG-ICE formulation
for one-dimensional ordinary differential equation integration. We
assess the ability of LS-MDG-ICE to move the grid in order to fit
a priori unknown interfaces as well as resolve internal and boundary
layers via anisotropic curvilinear $r$-adaptivity. Finally, we study
convergence of $L^{2}$ solution error under both grid, $h$, and
polynomial, $p$, refinement for LS-MDG-ICE in both the case of a
static and variable discrete geometry.

\subsection{Background}

Whereas traditional finite element methods choose trial and test spaces
a priori, the DPG methodology~\citep{Dem10,Dem11,Dem12} pioneered
an approach for automatically generating optimal test functions from
trial functions for a given variational problem. The test spaces are
chosen to ensure stability and are optimal in the sense that they
provide the best approximation error in the energy norm~\citep[Theorem 2.2]{Dem11}.
The DPG method makes use of discontinuous or broken test spaces~\citep{Car16},
enabling so-called practical DPG methods~\citep{Dem11,Gop13,Gop14}.
In its ultra-weak formulation, global coupling is then achieved through
the introduction of single-valued flux variables defined on the element
interfaces. Unlike standard DG, optimal error estimates in the $L^{2}$
norm for both the mesh size, $h$, and polynomial degree, $p$, were
proven for DPG~\citep[Theorems 3.2 and 3.3]{Dem10}.

DPG has been successfully applied to a wide range of applications
governed by both linear and nonlinear equations. It was applied to
Burgers equation and the compressible Navier-Stokes equations in one
dimension~\citep{Cha10}, the two-dimensional Navier-Stokes equations~\citep{Cha14},
Stokes equations~\citep{Rob14}, the incompressible Navier-Stokes
equations~\citep{Rob15}, Maxwell's equations~\citep{Car16}, and
Friedrichs' systems~\citep{Bui13}. Here we apply the DPG methodology
to inviscid and viscous MDG-ICE formulations, which are weak formulations
in what are referred to as trivial, or strong, form, i.e., derivatives
remain on the trial functions, and not moved to the test functions
via integration by parts, cf.~\citep{Dem15_20}. The resulting LS-MDG-ICE
formulation with optimal test functions minimizes the residual in
the $L^{2}$ norm, reducing to a classical least-squares formulation~\citep{Cai94,Boc98,Boc09}.
Furthermore, since the $L^{2}$ space admits discontinuous solutions,
there is no need to introduce additional interface variables~\citep{Car16}.
Instead, LS-MDG-ICE accounts for discontinuous solutions and achieves
global coupling through the enforcement of interface conditions on
the element interfaces.

Previous work has also applied least-squares finite element methods
to resolve discontinuous solutions through enforcement of a jump condition.
Gerritsma and Proot~\citep{Ger02} presented a least squares spectral
element method capable of accurately handling prescribed discontinuities
at element interfaces by incorporating a term that penalizes jumps
in the state. The method was successfully applied to a one-dimensional
test problem with a prescribed discontinuity. Cao and Gunzburger~\citep{Cao98}
directly incorporated the elliptic interface conditions into the least-squares
functional that allowed them to consider non-conforming, with respect
to the interface, approximating subspaces. MDG-ICE builds on these
methods by introducing the grid as a variable in order to align the
grid with any interfaces present in the a priori unknown exact solution,
including shocks.

MDG-ICE falls within the class $r$-adaptive methods that reposition
the grid points to better resolve the flow, cf.~Budd et al~\citep{Bud09}
and the survey of Huang and Russell~\citep{Hua10}. MDG-ICE also
falls within a class of shock fitting or tracking methods that to
align the grid with shocks or other discontinuous interfaces present
in the exact solution. However, MDG-ICE is an \emph{implicit} shock
fitting method, a new form of shock fitting that moves the grid and
implicitly fits a priori unknown shocks with arbitrary topology via
satisfaction of the weak form, thus overcoming a key limitation of
earlier \emph{explicit} shock fitting methods, cf. Moretti~\citep{Mor02}
and Salas~\citep{Sal09,Sal11}.

The optimization-based approach of Zahr and Persson~\citep{Zah17,Zah18_AIAA,Zah18}
is another promising form of implicit shock fitting, or tracking.
Their approach retains a standard discontinuous Galerkin method as
the state equation, while employing an objective function to detect
and fit interfaces present in the flow, while also treating the discrete
geometry as a variable. It has been used to compute very accurate
solutions to discontinuous inviscid flows on coarse grids without
the use of artificial stabilization. Recently, Zahr et. al.~\citep{Zah19,Zah20_SCITECH}
have extended their implicit shock tracking framework with a new objective
function based on an enriched test space, an improved SQP solver,
and grid topology modification. Furthermore, they extend the Levenberg-Marquardt
elliptic or Laplacian regularization used by the present authors~\citep{Cor19_AVIATION},
to incorporate a factor proportional to the inverse volume of each
element in order to account for variations in the element size throughout
a given grid. In this work, we incorporate their improved scaling
into the regularization we employ for the solution of a Mach 5 compressible
Navier-Stokes flow over a blunt body at Reynolds numbers of $10^{3}$
and $10^{4}$ in Section~\ref{sec:viscous-bow-shock}.

Grid design, or grid point placement, is not only important at discontinuous
interfaces, it is critical for achieving optimal-order accuracy in
the case of problems that contain regions with sharp, but smooth,
gradients, e.g., internal and boundary layers, which are common features
in solutions to singularly perturbed boundary value problems of the
form

\begin{align}
\nabla\cdot\left(vy-\varepsilon\nabla y\right)-ay-f=0 & \textup{ in }\Omega,\label{eq:singularly-perturbed-strong-form}
\end{align}
where $\Omega\subset\mathbb{R}^{d}$ is a given domain, $y:\mathbb{R}^{d}\rightarrow\mathbb{R}$,
$v\in\mathbb{R}^{d}$ is uniform, $a\in\mathbb{R}$ is constant, and
$\varepsilon\in\left(0,1\right)$ is constant. Accurate approximations
to singularly perturbed problems are difficult to obtain since they
require either very small element sizes or a large number of collocation
points~\citep{Sch83,Bal99,Tan96}. Oscillations in the solution can
be reduced through the addition of stabilizing terms. However, these
terms can pollute high-order solutions with low-order errors and have
the effect of masking a poor quality approximation to a difficult
problem with a better approximation to a simpler problem~\citep{Gre81}.
We advocate an alternative approach, improving the approximation by
incorporating the discrete geometry as a solution variable and solving
for an optimal grid.

Additionally, one can seek the design of optimal grids, in the sense
of efficiency and accuracy, which, in the context of singularly perturbed
problems, has been an active area of research for the past 50 years~\citep{Lin03,Lin10}.
Boundary layer-adapted grids were first introduced by Bakhvalov~\citep{Bak69}.
The so-called Bakhvalov-meshes combine a nonuniform grid inside the
boundary layer with an otherwise uniform grid. Exponentially graded
grids that do not transition to a uniform grid have been proposed~\citep{Xen96,Xen99},
but were only derived in the context of the reaction-diffusion equation.
Simpler Shishkin-meshes~\citep{Shi90} combine two piecewise uniform
grids where the boundary layer is contained in the finer region. Shishkin-meshes
are usually less accurate than Bakhvalov meshes but easier to analyze.
Regardless, they share the difficulty of defining both the transition
point from fine-to-coarse as well as the grid resolution in the coarse
region. Furthermore, all three methods for generating layer-adapted
grids fail to account for the interior nodes of high-order curvilinear
element shape representations. In contrast, MDG-ICE accounts for these
issues as an intrinsic part of the solver.

Polynomial ($p$) refinement is an alternative to the quasi-uniform,
or exponential, grid ($h$) refinement used in the design of Bakhvalov,
Shishkin, and exponentially graded grids. Schwab and Suri~\citep{Sch96a}
showed that at least two elements are required for exponential convergence
in the energy norm under polynomial refinement using standard finite
element approximations for problems with boundary layer type solutions.
In addition, they showed that the size of the element adjacent to
the boundary layer must be on the order of the diffusive scale, i.e.
the size of the boundary layer~\citep{Sch96a}. For convection-reaction-diffusion
the size of the boundary layer is on the order of $\mathcal{O}\left(\varepsilon=\nicefrac{1}{\mathrm{Pe}}\right)$
while it is it is order $\mathcal{O}\left(\sqrt{\varepsilon}=\sqrt{\nicefrac{1}{\mathrm{Pe}}}\right)$
reaction-diffusion~\citep{Lin10,Syk19,Syk19hp}. The two element,
adaptive version is referred to as the $hp$ version, however, Schwab
and Suri note that $rp$ is a more appropriate term for the scheme
since the size of elements, not the number of elements, change as
the polynomial degree is refined. Their result provide strong evidence
of the intrinsic coupling between the solution to singularly perturbed
boundary value problems and the corresponding optimal mesh design.
To achieve spectral, or exponential convergence, the grid must adapt
to the approximate solution. Furthermore, to achieve both algebraic
and exponential convergence for problems with evolving or a priori
unknown layers, internal or boundary, the solution process must incorporate
grid design, thus further motivating the application of MDG-ICE to
viscous problems.

\section{Optimal Test Spaces\label{subsec:Optimal-Test-Spaces}}

In this section, we derive optimal test spaces for MDG-ICE using the
DPG methodology~\citep{Dem10,Dem11,Dem12,Dem15_03,Dem15_20,Dem18,Gop13,Gop14,Car16}
in a specialized setting. In earlier work~\citep{Cor18,Ker20}, we
proposed an MDG-ICE formulation using equal-order test spaces chosen
a priori, which was solved using a discrete least-squares solver.
In order to demonstrate optimal convergence under grid refinement
for flows with shocks and other discontinuous interfaces, we introduced
a generalized discretization strategy, cf.~\citep[Section 6]{Cor18}.
Using this generalized discretization strategy, we considered two
distinct discretizations. One discretization \emph{only} enforced
the interface condition on interfaces present in the coarse grid,
while employing a DG-type numerical flux on interfaces introduced
on refined grids, achieving optimal-order convergence for a flow with
a curved bow shock using higher-order elements, cf.~\citep[Section 6.3]{Cor18}.
The other discretization enforced the interface condition on \emph{all}
interfaces present in the discrete grid at \emph{all} levels of refinement,
cf.~\citep[Section 6.2]{Cor18}. Using this strategy, it was observed
for the case of two-dimensional, steady, linear advection, that an
ad hoc scaling of the test functions is required in order to achieve
optimal-order convergence, a problem which will be revisited in Section~\ref{subsec:sinusoidal-waves-convergence}.
In the present work, we seek to weight test functions in a systematic
fashion using the DPG methodology in order to ensure optimal-order
convergence. 

We follow the work of Demkowicz and Gopalakrishnan~\citep{Dem11,Gop13}
in order to derive optimal test spaces for MDG-ICE in the special
case of a linear, ordinary differential equation (ODE) defined on
a static grid. The generalization for multi-dimensional flow problems
on moving grids will be presented in Section~\ref{sec:Formulation}.
The ODE is defined in strong form for piecewise smooth $y$,
\begin{align}
y'=f & \textup{ in }\kappa\qquad\forall\kappa\in\mathcal{T},\label{eq:conservation-strong}\\
\left\llbracket y\right\rrbracket =0 & \textup{ on }\epsilon\qquad\forall\epsilon\in\mathcal{E}.\label{eq:interface-condition-strong}
\end{align}
We assume that $\Omega=\left(0,1\right)$ is the domain, which is
partitioned by $\mathcal{T}$, consisting of disjoint sub-domains
or cells $\kappa$, so that $\overline{\Omega}=\cup_{\kappa\in\mathcal{T}}\overline{\kappa}$,
with interfaces $\epsilon$, composing a set $\mathcal{E}$ so that
$\cup_{\epsilon\in\mathcal{E}}\epsilon=\cup_{\kappa\in\mathcal{T}}\partial\kappa$.
Here, $f$ is state independent source term data, i.e., it does not
depend on $y$. An inflow boundary condition $y=y_{\text{in}}$ is
imposed via~(\ref{eq:interface-condition-strong}) at the inflow
interface $\epsilon_{\text{in}}\in\mathcal{E}$ located at $x=0$.
Over interior interfaces, $\mathcal{E}_{0}=\left\{ \epsilon\in\mathcal{E}\bigl|\epsilon\cap\partial\Omega=\emptyset\right\} $,
the jump is defined as $\left\llbracket y\right\rrbracket =y^{+}-y^{-}$
between the right and left limits. We recall~\citep[Definition 3]{Gop13},
which is restated here as follows: for a given trial space $Y_{h}\subset Y$,
the optimal test space for the continuous bilinear form $b\left(\cdot,\cdot\right):Y\times V\rightarrow\mathbb{R}$
is

\begin{equation}
V_{h}=T\left(Y_{h}\right),\label{eq:optimal-test-space}
\end{equation}
where the trial-to-test operator $T:Y\rightarrow V$ is defined by

\begin{equation}
\left(Tz,v\right)_{V}=b\left(z,y\right)\qquad\forall v\in V,\:\forall z\in Y.\label{eq:trial-to-test}
\end{equation}
Using the above definition, we compute the optimal test space corresponding
to the MDG-ICE formulation for a single-element, $\Omega=\kappa=\left(0,1\right)$,
spectral discretization. Given $y\in Y$ and $\left(v,w\right)\in V\times W$,
the bilinear form $b:Y\times\left(V\times W\right)\rightarrow\mathbb{R}$
is defined by
\begin{equation}
b\left(y,\left(v,w\right)\right)=\left(y',v\right)_{L^{2}\left(\kappa\right)}-y\left(0\right)w,\label{eq:advection-1d-spectral-bilinear-form}
\end{equation}
while the data functional $\ell\in\left(V\times W\right)^{*}$ is
\begin{equation}
\ell\left(v,w\right)=\left(f,v\right)_{L^{2}\left(\kappa\right)}+y_{\text{in}}w,\label{eq:advection-1d-spectral-data}
\end{equation}
whose Riesz representation is 
\begin{equation}
\left(f,y_{in}\right)\in V\times W.\label{eq:advection-1d-spectral-data-riesz}
\end{equation}
The trial space is

\begin{equation}
Y=H^{1}\left(\kappa\right)=\left\{ y\in L^{2}\left(\kappa\right)\bigl|y'\in L^{2}\left(\kappa\right)\right\} ,\label{eq:advection-1d-spectral-trial-space}
\end{equation}
where $y'$ denotes the weak derivative, and the test space is

\begin{equation}
V\times W=L^{2}\left(\kappa\right)\times\mathbb{R}.\label{eq:advection-1d-spectral-test-space}
\end{equation}
For higher-dimensional spatial domains, the trace space $W$ would
have to be defined in terms of a fractional Sobolev space on the inflow
interface, denoted $\epsilon_{\text{in}}$, however in the present
setting we follow~\citep{Gop13} and treat the trace space $W$ as
$\mathbb{R}$. The norm is defined as
\begin{equation}
\left\Vert \left(v,w\right)\right\Vert _{V\times W}^{2}=\left\Vert v\right\Vert _{L^{2}\left(\Omega\right)}^{2}+\left|w\right|^{2}.\label{eq:advection-1d-spectral-norm}
\end{equation}
Following~\citep[Examples 18 and 19]{Gop13}, it can be shown that
the trial-to-test operator $T:Y\rightarrow V\times W$ is defined
in this case by

\begin{equation}
Ty=\left(y',y\left(0\right)\right).\label{eq:trial-to-test-advection-1d-spectral}
\end{equation}
In particular, if the finite-dimensional trial space $Y_{h}\subset Y,$
is given by $Y_{h}=\mathcal{P}_{p}\left(\Omega\right)$, where $\mathcal{P}_{p}$
denotes the space of polynomials of degree at most $p$, then 
\begin{equation}
TY_{h}=\mathcal{P}_{p-1}\left(\Omega\right)\times\mathbb{R},\label{eq:trial-to-test-advection-1d-spectral-polynomial}
\end{equation}
since if $y_{h}\in Y_{h}=\mathcal{P}_{p}\left(\Omega\right)$ then
$y'\in\mathcal{P}_{p-1}\left(\Omega\right)$ and $y\left(0\right)\in\mathbb{R}$. 

Following~\citep[Example 25]{Gop13}, we generalize the one-element
formulation to a discretization over a partitioned domain. We consider
a partition of $\Omega=\left(0,1\right)$ into a number of elements
or intervals $\kappa_{i}=\left(x_{i-1},x_{i}\right)$ with $x_{0}=0$
and $x_{m}=1$ for $i=1,\ldots,m$. We construct the MDG-ICE formulation
over the broken Sobolev space
\begin{equation}
Y=H^{1}\left(\mathcal{T}\right)=\left\{ y\in L^{2}\left(\Omega\right)\bigl|\forall\kappa\in\mathcal{T},\left.y\right|_{\kappa}\in H^{1}\left(\kappa\right)\right\} \label{eq:advection-1d-trial-space}
\end{equation}
with norm defined by 
\begin{equation}
\left\Vert \cdot\right\Vert _{Y}^{2}=\left\Vert \cdot\right\Vert _{H^{1}\left(\mathcal{T}\right)}^{2}=\sum_{\kappa\in\mathcal{T}}\left\Vert \cdot\right\Vert _{H^{1}\left(\kappa\right)}^{2}.\label{eq:advection-1d-norm}
\end{equation}
The test space is 
\begin{equation}
V\times W=L^{2}\left(\Omega\right)\times\mathbb{R}^{\left|\mathcal{E}_{0}\right|+1}\label{eq:advection-1d-test-space}
\end{equation}
where $\left|\mathcal{E}_{0}\right|+1$ is the number of interior
and inflow interfaces, which is equal to the number of cells $m$.
Given $y\in Y$ and $\left(v,w\right)\in V\times W$, the bilinear
form $b:Y\times\left(V\times W\right)\rightarrow\mathbb{R}$ is 
\begin{equation}
b\left(y,\left(v,w\right)\right)=\sum_{\kappa\in\mathcal{T}}\left(y',v\right)_{L^{2}\left(\kappa\right)}-y\left(0\right)w-\sum_{\epsilon\in\mathcal{E}_{0}}\left\llbracket y\right\rrbracket w\label{eq:advection-1d-bilinear-form}
\end{equation}
with data functional $\ell\left(v,w\right)\in\left(V\times W\right)^{*}$
\begin{equation}
\ell\left(v,w\right)=\sum_{\kappa\in\mathcal{T}}\left(f,v\right)_{L^{2}\left(\kappa\right)}+y_{\text{in}}w.\label{eq:advection-1d-data}
\end{equation}
that has Riesz representation 
\begin{equation}
\left(f,\left(y_{\text{in}},\ldots,\left.0\right|_{\epsilon\in\mathcal{E}_{0}},\ldots\right)\right)\in V\times W.\label{eq:advection-1d-data-riesz}
\end{equation}
 The trial-to-test operator $T:Y\rightarrow V\times W$ is

\begin{equation}
Ty=\left(y',\left(y\left(0\right),\ldots,\left.\left\llbracket y\right\rrbracket \right|_{\epsilon\in\mathcal{E}_{0}},\ldots\right)\right).\label{eq:trial-to-test-advection-1d}
\end{equation}
If the finite-dimensional trial space $Y_{h}\subset Y,$ is given
by

\begin{equation}
Y_{h}=\left\{ y\in Y\middlebar\forall\kappa\in\mathcal{T},\left.y\right|_{\kappa}\in\mathcal{P}_{p}\right\} ,\label{eq:advection-1d-discrete-solution-space-simplex}
\end{equation}
which we denote $\mathcal{P}_{p}\left(\mathcal{T}\right)$, then 
\begin{equation}
T\left(Y_{h}\right)=\mathcal{P}_{p-1}\left(\mathcal{T}\right)\times\mathbb{R}^{\left|\mathcal{E}_{0}\right|+1},\label{eq:trial-to-test-advection-1d-polynomial}
\end{equation}
since if $y_{h}\in\mathcal{P}_{p}\left(\mathcal{T}\right)$ then $y'\in\mathcal{P}_{p-1}\left(\mathcal{T}\right)$
while $\left(y\left(0\right),\ldots,\left.\left\llbracket y\right\rrbracket \right|_{\epsilon\in\mathcal{E}_{0}},\ldots\right)\in\mathbb{R}^{\left|\mathcal{E}_{0}\right|+1}$.

Having identified the optimal test space, in this specialized setting,
the LS-MDG-ICE method with optimal test functions can be implemented
in a number of equivalent ways. One way would be to discretize the
optimal space~(\ref{eq:trial-to-test-advection-1d-polynomial}) directly
using a standard finite element basis for $\mathcal{P}_{p-1}\left(\mathcal{T}\right)$:
find $y\in Y_{h}$ such that
\begin{equation}
b\left(y,\left(v,w\right)\right)=\sum_{\kappa\in\mathcal{T}}\left(y',v\right)_{L^{2}\left(\kappa\right)}-y\left(0\right)w-\sum_{\epsilon\in\mathcal{E}_{0}}\left\llbracket y\right\rrbracket w=0\qquad\forall\left(v,w\right)\in\mathcal{P}_{p-1}\left(\mathcal{T}\right)\times\mathbb{R}^{\left|\mathcal{E}_{0}\right|+1}.\label{eq:advection-1d-weak-form-pm1-test-space}
\end{equation}
We observe that in this case the dimensionality of the trial space
and optimal test space match,
\begin{align}
\dim\left(Y_{h}\right) & =\dim\left(\mathcal{P}_{p}\left(\mathcal{T}\right)\right)\nonumber \\
 & =m\left(p+1\right)\nonumber \\
 & =m\left(\left(p-1\right)+1\right)+m\nonumber \\
 & =\dim\left(\mathcal{P}_{p-1}\left(\mathcal{T}\right)\right)+\dim\left(\mathbb{R}^{\left|\mathcal{E}_{0}\right|+1}\right)\nonumber \\
 & =\dim\left(T\left(Y_{h}\right)\right).\label{eq:discrete-space-dimensionality}
\end{align}
 In comparison, in previous work, we used an MDG-ICE formulation that
assumed equal-order polynomial test spaces, leading to a weak formulation:
find $y\in Y_{h}$ such that

\begin{equation}
b\left(y,\left(v,w\right)\right)=\sum_{\kappa\in\mathcal{T}}\left(y',v\right)_{L^{2}\left(\kappa\right)}-y\left(0\right)w-\sum_{\epsilon\in\mathcal{E}_{0}}\left\llbracket y\right\rrbracket w=0\qquad\forall\left(v,w\right)\in\mathcal{P}_{p}\left(\mathcal{T}\right)\times\mathbb{R}^{\left|\mathcal{E}_{0}\right|+1},\label{eq:advection-1d-weak-form-p-test-space}
\end{equation}
for which the dimensionality of the test space is higher than the
trial space, necessitating a discrete least-squares solver strategy~\citep{Cor18}.

In a more general setting, it will not be possible to identify the
optimal test space as a standard polynomial space. The DPG methodology,
however, provides a systematic and generalizable approach to generating
optimal test spaces automatically from the chosen discrete trial space,
without explicit identification of the optimal test space. Since we
formulate MDG-ICE using a so-called trivial or strong formulation
as opposed to an ultra-weak formulation, cf.~\citep{Dem15_20}, we
compose the weak formulation with the trial-to-test operator to systematically
generate optimal test functions from a standard discontinuous finite
element basis of the trial space $\mathcal{P}_{p}\left(\mathcal{T}\right)$:
find $y\in Y_{h}$ such that

\begin{equation}
b\left(y,T\left(v\right)\right)=h\sum_{\kappa\in\mathcal{T}}\left(y',v'\right)_{L^{2}\left(\kappa\right)}-y\left(0\right)v\left(0\right)-\sum_{\epsilon\in\mathcal{E}_{0}}\left\llbracket y\right\rrbracket \left\llbracket v\right\rrbracket \qquad\forall v\in Y_{h},\label{eq:advection-1d-weak-form-trial-to-test}
\end{equation}
where $h$ is a scaling discussed in Section~\ref{subsec:dpg-weak-form}.

\begin{figure}
\subfloat[\label{fig:mdg-ice-convergence-linear-advection-p-minus-one}Convergence
of MDG-ICE($\mathcal{P}_{2}$) discretizations using the sub-optimal,
$\mathcal{P}_{2}\left(\mathcal{T}_{h}\right)\times\mathbb{R}^{\left|\mathcal{E}_{h,0}\right|+1}$,
test space given by weak formulation~(\ref{eq:advection-1d-weak-form-p-test-space-MESH-DEPENDENT})
and the optimal, $\mathcal{P}_{1}\left(\mathcal{T}_{h}\right)\times\mathbb{R}^{\left|\mathcal{E}_{h,0}\right|+1}$,
test space given by weak formulation~(\ref{eq:advection-1d-weak-form-pm1-test-space-MESH-DEPENDENT}).]{\begin{centering}
\includegraphics[width=0.45\linewidth]{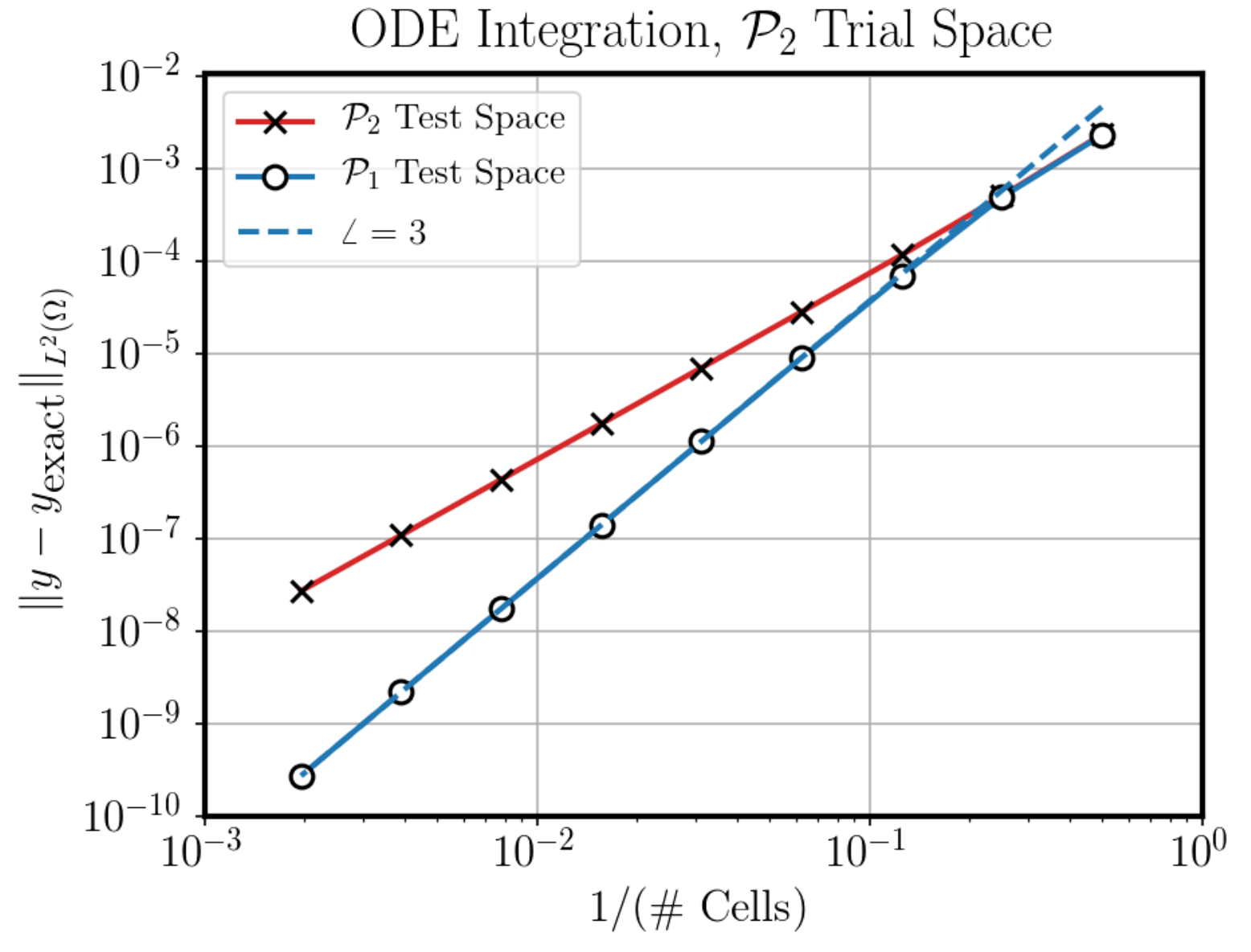}
\par\end{centering}
}\hfill{}
\begin{centering}
\subfloat[\label{fig:mdg-ice-convergence-linear-advection-dpg}Convergence of
MDG-ICE($\mathcal{P}_{2}$) discretizations using the sub-optimal,
$\mathcal{P}_{2}\left(\mathcal{T}_{h}\right)\times\mathbb{R}^{\left|\mathcal{E}_{h,0}\right|+1}$,
test space given by weak formulation~(\ref{eq:advection-1d-weak-form-p-test-space-MESH-DEPENDENT})
and the optimal, $T\left(Y_{h}\right)$, test space given by weak
formulation~(\ref{eq:advection-1d-weak-form-trial-to-test-MESH-DEPENDENT}).]{\begin{centering}
\includegraphics[width=0.45\linewidth]{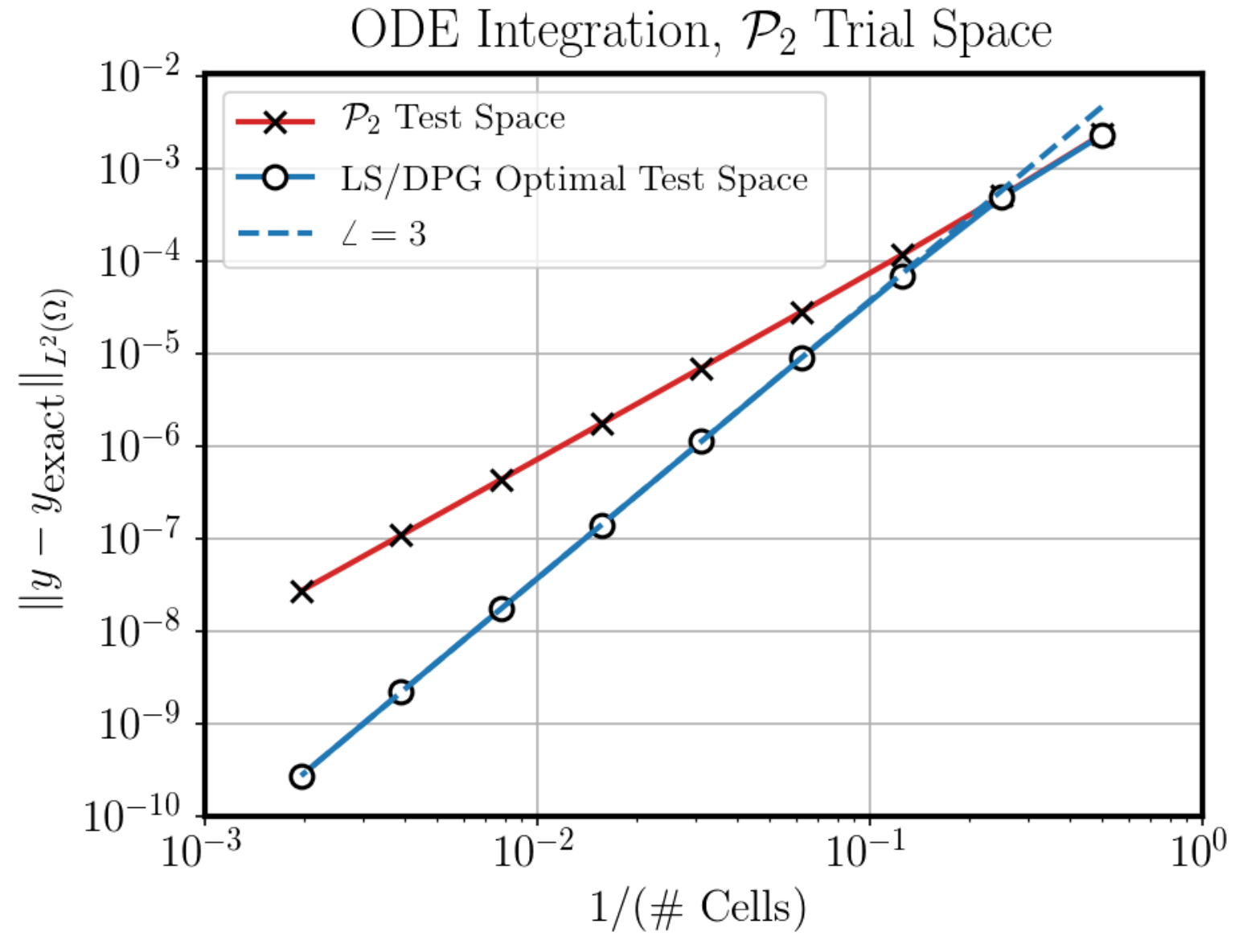}
\par\end{centering}
}
\par\end{centering}
\caption{\label{fig:mdg-ice-convergence-linear-advection}Convergence plots
for solving an ODE integration problem with exact solution~(\ref{eq:exact-solution})
using MDG-ICE($\mathcal{P}_{2}$) on a static grid with the sub-optimal,
$\mathcal{P}_{2}\left(\mathcal{T}_{h}\right)\times\mathbb{R}^{\left|\mathcal{E}_{h,0}\right|+1}$,
test space~(\ref{eq:advection-1d-weak-form-p-test-space-MESH-DEPENDENT}),
the optimal, $\mathcal{P}_{1}\left(\mathcal{T}_{h}\right)\times\mathbb{R}^{\left|\mathcal{E}_{h,0}\right|+1}$,
test space~(\ref{eq:advection-1d-weak-form-pm1-test-space-MESH-DEPENDENT}),
and the automatically generated optimal test space~(\ref{eq:advection-1d-weak-form-trial-to-test-MESH-DEPENDENT}).
The coarsest grid consisted of $2$ line cells, while the finest grid
consisted of $512$ line cells. Weak formulations~(\ref{eq:advection-1d-weak-form-pm1-test-space-MESH-DEPENDENT})
and~(\ref{eq:advection-1d-weak-form-trial-to-test-MESH-DEPENDENT})
achieve optimal, third-order convergence, while the weak formulation~(\ref{eq:advection-1d-weak-form-p-test-space-MESH-DEPENDENT})
exhibits sub-optimal convergence behavior.}
\end{figure}

In order to verify both the optimality of weak formulations~(\ref{eq:advection-1d-weak-form-pm1-test-space})
and~(\ref{eq:advection-1d-weak-form-trial-to-test}) and the sub-optimality
of the weak formulation~(\ref{eq:advection-1d-weak-form-p-test-space}),
we study their convergence under mesh refinement. Mesh refinement
introduces a family of meshes $\mathcal{T}_{h}$ for which the set
of interfaces $\mathcal{E}_{h}$ is increasingly dense as $h\rightarrow0$.
Ultra-weak DPG methods, which ``break'' the test space and introduce
a mesh-dependent variational form~cf.~\citep{Dem11,Car16}, introduce
interface unknowns on a mesh-dependent interface trial space over
$\mathcal{E}_{h}$. In an analogous fashion, the DPG formulation of
MDG-ICE, which ``breaks'' the trial space by defining $Y_{h}=\mathcal{P}_{p}\left(\mathcal{T}_{h}\right)$
and introduces a mesh-dependent variational form, enforces the interface
condition on a mesh-dependent interface test space over $\mathcal{E}_{h}$.
The mesh-dependent variational form is given as: find $y\in\mathcal{P}_{p}\left(\mathcal{T}_{h}\right)$
such that

\begin{equation}
b_{h}\left(y,\left(v,w\right)\right)=\sum_{\kappa\in\mathcal{T}_{h}}\left(y',v\right)_{L^{2}\left(\kappa\right)}-y\left(0\right)w-\sum_{\epsilon\in\mathcal{E}_{h,0}}\left\llbracket y\right\rrbracket w=0\qquad\forall\left(v,w\right)\in\mathcal{P}_{p-1}\left(\mathcal{T}_{h}\right)\times\mathbb{R}^{\left|\mathcal{E}_{h,0}\right|+1}.\label{eq:advection-1d-weak-form-pm1-test-space-MESH-DEPENDENT}
\end{equation}
or, equivalently,

\begin{equation}
b_{h}\left(y,T_{h}\left(v\right)\right)=h\sum_{\kappa\in\mathcal{T}_{h}}\left(y',v'\right)_{L^{2}\left(\kappa\right)}-y\left(0\right)v\left(0\right)-\sum_{\epsilon\in\mathcal{E}_{h,0}}\left\llbracket y\right\rrbracket \left\llbracket v\right\rrbracket \qquad\forall v\in\mathcal{P}_{p}\left(\mathcal{T}_{h}\right).\label{eq:advection-1d-weak-form-trial-to-test-MESH-DEPENDENT}
\end{equation}
As a result, the trial-to-test operator inherits a dependence on $\mathcal{T}_{h}$.
Such a dependence also arises in the context of ``practical'' DPG
methods where the inverse Riesz map defining the trial-to-test operator
is approximated in a mesh-dependent fashion, the effect of which has
been analyzed by Gopalakrishnan and Qiu~\citep{Gop14}. We defer
such analysis to future work and instead in the present study investigate
the mesh or grid convergence behavior of solutions to ~(\ref{eq:advection-1d-weak-form-pm1-test-space-MESH-DEPENDENT})
and (\ref{eq:advection-1d-weak-form-trial-to-test-MESH-DEPENDENT})
as well as discrete least-squares solutions to: find $y\in\mathcal{P}_{p}\left(\mathcal{T}_{h}\right)$
such that
\begin{equation}
b_{h}\left(y,\left(v,w\right)\right)=\sum_{\kappa\in\mathcal{T}_{h}}\left(y',v\right)_{L^{2}\left(\kappa\right)}-y\left(0\right)w-\sum_{\epsilon\in\mathcal{E}_{h,0}}\left\llbracket y\right\rrbracket w=0\qquad\forall\left(v,w\right)\in\mathcal{P}_{p}\left(\mathcal{T}_{h}\right)\times\mathbb{R}^{\left|\mathcal{E}_{h,0}\right|+1}.\label{eq:advection-1d-weak-form-p-test-space-MESH-DEPENDENT}
\end{equation}

We consider discrete trial spaces defined by $Y_{h}=\mathcal{P}_{2}\left(\mathcal{T}_{h}\right)$,
and measure the $L^{2}$ error of these formulations with respect
to an exact solution given by 
\begin{equation}
y\left(x\right)=\left(x-0.1\right)\left(x-0.2\right)\left(x-0.3\right)\left(x-0.4\right)\left(x-0.5\right)\left(x-0.9\right)\label{eq:exact-solution}
\end{equation}
generated by data $f=y'$ and $y_{\text{in}}=y\left(0\right)$. Figure~\ref{fig:mdg-ice-convergence-linear-advection-p-minus-one}
shows that the formulation~(\ref{eq:advection-1d-weak-form-pm1-test-space-MESH-DEPENDENT})
defined using a $\mathcal{P}_{1}\left(\mathcal{T}_{h}\right)\times\mathbb{R}^{\left|\mathcal{E}_{h,0}\right|+1}$
test space achieves optimal, third-order convergence. Likewise, Figure~\ref{fig:mdg-ice-convergence-linear-advection-dpg}
shows that the equivalent formulation~(\ref{eq:advection-1d-weak-form-trial-to-test-MESH-DEPENDENT})
defined using the trial-to-test operator also achieves optimal, third-order
convergence. We also observe the sub-optimal convergence behavior
of the formulation~(\ref{eq:advection-1d-weak-form-p-test-space-MESH-DEPENDENT})
using an equal-order test space $\mathcal{P}_{2}\left(\mathcal{T}_{h}\right)\times\mathbb{R}^{\left|\mathcal{E}_{h,0}\right|+1}$.
Having verified optimal-order convergence, in the following section,
we proceed to extend LS-MDG-ICE from formulation~(\ref{eq:advection-1d-weak-form-trial-to-test-MESH-DEPENDENT})
to a general setting.

\section{Least-Squares Moving Discontinuous Galerkin Method with Interface
Condition Enforcement\label{sec:Formulation}}

In this section, we derive the least-squares formulation of MDG-ICE
for conservation laws governed by convection-diffusion with state
independent source term data in arbitrary spatial dimensions. Let
$\Omega\subset\mathbb{R}^{d}$ be a given domain. Our earlier work
also considered unsteady problems where the domain is defined over
space-time $\Omega\subset\mathbb{R}^{d=d_{x}+1}$, cf.~\citep{Cor18,Ker20}.
Although LS-MDG-ICE readily applies to unsteady problems, in the present
work, we focus on steady problems defined over a spatial domain $\Omega\subset\mathbb{R}^{d=d_{x}}$.
Throughout this work, we assume that $\Omega$ is partitioned by $\mathcal{T}$,
consisting of disjoint sub-domains or cells $\kappa$, so that $\overline{\Omega}=\cup_{\kappa\in\mathcal{T}}\overline{\kappa}$,
with interfaces $\epsilon$, composing a set $\mathcal{E}$ so that
$\cup_{\epsilon\in\mathcal{E}}\epsilon=\cup_{\kappa\in\mathcal{T}}\partial\kappa$.
Furthermore, we assume that each interface $\epsilon$ is oriented
so that a unit normal $n:\epsilon\rightarrow\mathbb{R}^{d}$ is defined.
For consistency of notation with earlier work~\citep{Ker20}, we
also define the spatial normal $n_{x}:\epsilon\rightarrow\mathbb{R}^{d_{x}}$,
which coincides with the normal $n$ in the spatial setting.

In Section~\ref{subsec:Governing-equations-Primal-Formulation} we
present the governing equations. Then, in Section~\ref{subsec:formulation-physical-space}
we derive a least-squares formulation of MDG-ICE for both inviscid
and viscous flows in physical space with static, discrete geometry.
The most powerful aspect of MDG-ICE is its ability to fit a priori
unknown interfaces and resolve otherwise under-resolved flow features.
Therefore, in Section~\ref{subsec:variable-geometry} we transform
to reference space and introduce the discrete geometry as a variable,
so that using DPG we systematically generate optimal test functions
for the case of moving grid geometry.

\subsection{Governing equations\label{subsec:Governing-equations-Primal-Formulation}~\citep[Section 2.1]{Ker20} }

In this work, we will apply LS-MDG-ICE to governing equations also
considered in our earlier work~\citep{Ker20}, cf. Section~2.1.
Consider a nonlinear conservation law governing the behavior of smooth,
$\mathbb{R}^{m}$-valued, functions $y$,

\begin{align}
\nabla\cdot\mathcal{F}\left(y,\nabla_{x}y\right)-\mathcal{S}\left(y\right)=0 & \textup{ in }\Omega,\label{eq:conservation-strong-primal}
\end{align}
in terms of a given flux function, $\mathcal{F}:\mathbb{R}^{m}\times\mathbb{R}^{m\times d_{x}}\rightarrow\mathbb{R}^{m\times d}$
that depends on the flow state variable $y$ and its $d_{x}$-dimensional
spatial gradient,
\begin{equation}
\nabla_{x}y=\left(\frac{\partial y}{\partial x_{1}},\ldots,\frac{\partial y}{\partial x_{d_{x}}}\right).\label{eq:spatial-gradient}
\end{equation}
For unsteady flows, $\mathcal{F}$ is defined to be the space-time
flux with $d=d_{x}+1$, cf.~\citep{Ker20}. In the present work,
we focus on the case of steady flows, $d=d_{x}$, so that $\mathcal{F}$
coincides with a given spatial flux, $\mathcal{F}^{x}:\mathbb{R}^{m}\times\mathbb{R}^{m\times d_{x}}\rightarrow\mathbb{R}^{m\times d_{x}}$.
The spatial flux is defined in terms of a convective flux, which depends
on the state variable only, and a viscous, or diffusive flux, which
also depends on the spatial gradient of the state variable,
\begin{equation}
\mathcal{F}^{x}\left(y,\nabla_{x}y\right)=\mathcal{F}^{c}\left(y\right)-\mathcal{F}^{v}\left(y,\nabla_{x}y\right).\label{eq:navier-stokes-spatial-flux}
\end{equation}
The governing equation~(\ref{eq:conservation-strong-primal}) also
incorporates a source term $\mathcal{S}:\mathbb{R}^{m}\rightarrow\mathbb{R}^{m}$,
which in previous work we assumed to be zero~\citep{Cor18,Ker20}.
Here, we incorporate a source term of the form $\mathcal{S}\left(y\right)=f$,
defined in terms of data $f\in\mathbb{R}^{m}$ that is independent
of both the state $y$ and position $x\in\Omega$.

\subsubsection{Linear advection-diffusion\label{subsec:Linear-advection-diffusion}}

Linear advection-diffusion involves a single-component flow state
variable $y:\Omega\rightarrow\mathbb{R}^{1}$ with a linear diffusive
flux,
\begin{equation}
\mathcal{F}^{v}\left(y,\nabla_{x}y\right)=\epsilon\nabla_{x}y,\label{eq:linear-viscous-flux}
\end{equation}
that is independent of the state $y$, where the coefficient $\epsilon$
represents mass diffusivity. The convective flux is given as 
\begin{equation}
\mathcal{F}^{c}\left(y\right)=\left(v_{1}y,\ldots,v_{d_{x}}y\right),\label{eq:linear-convective-flux}
\end{equation}
where $\left(v_{1},\ldots,v_{d_{x}}\right)\in\mathbb{R}^{d_{x}}$
is a prescribed spatial velocity that in the present setting is assumed
to be uniform.

\subsubsection{One-dimensional Burgers flow\label{subsec:One-dimensional-Burgers-flow}}

As in the case of linear advection-diffusion, one-dimensional Burgers
flow involves a single-component flow state variable $y:\Omega\rightarrow\mathbb{R}^{1}$
with a linear viscous flux,
\begin{equation}
\mathcal{F}^{v}\left(y,\nabla y\right)=\epsilon\nabla y,\label{eq:linear-viscous-flux-Burgers}
\end{equation}
is independent of the state $y$, where the coefficient $\epsilon$
represents viscosity. The convective flux is given as 
\begin{equation}
\mathcal{F}^{c}\left(y\right)=\left(\frac{1}{2}y^{2}\right).\label{eq:convective-burgers-flux}
\end{equation}

\subsubsection{Compressible Navier-Stokes flow\label{subsec:Compressible-Navier-Stokes-flow}}

For compressible Navier-Stokes flow, the state variable $y:\Omega\rightarrow\mathbb{R}^{m}$,
where $m=d_{x}+2$, is given by

\begin{equation}
y=\left(\rho,\rho v_{1},\ldots,\rho v_{d_{x}},\rho E\right).\label{eq:navier-stokes-state}
\end{equation}
The $i$-th spatial component of the convective flux, $\mathcal{F}^{c}:\mathbb{R}^{m}\rightarrow\mathbb{R}^{m\times d_{x}}$,
is
\begin{equation}
\mathcal{F}_{i}^{c}\left(y\right)=\left(\rho v_{i},\rho v_{i}v_{1}+p\delta_{i1},\ldots,\rho v_{i}v_{d_{x}}+p\delta_{id_{x}},\rho Hv_{i}\right),\label{eq:ns-convective-flux-spatial-component}
\end{equation}
where $\delta_{ij}$ is the Kronecker delta, $\rho:\Omega\rightarrow\mathbb{R}_{+}$
is density, $\left(v_{1},\ldots,v_{d_{x}}\right):\mathbb{R}^{m}\rightarrow\mathbb{R}^{d_{x}}$
is velocity, $\rho E:\Omega\rightarrow\mathbb{R}_{+}$ is stagnation
energy per unit volume, and 
\begin{equation}
H=\left(\rho E+p\right)/\rho\label{eq:stagnation-enthalpy}
\end{equation}
 is stagnation enthalpy, where $H:\mathbb{R}^{m}\rightarrow\mathbb{R}_{+}$.
Assuming that the fluid is a perfect gas, the pressure $p:\mathbb{R}^{m}\rightarrow\mathbb{R}_{+}$
is defined as 
\begin{equation}
p=\left(\gamma-1\right)\left(\rho E-\frac{1}{2}\sum_{i=1}^{d_{x}}\rho v_{i}v_{i}\right),\label{eq:pressure}
\end{equation}
where the ratio of specific heats for air is given as $\gamma=1.4$.
The $i$-th spatial component of the viscous flux is given by
\begin{equation}
\mathcal{F}_{i}^{\nu}\left(y,\nabla y\right)=\left(0,\tau_{1i},\ldots,\tau_{d_{x}i},\sum_{j=1}^{d_{x}}\tau_{ij}v_{j}-q_{i}\right),\label{eq:navier-stokes-viscous-flux-spatial-component}
\end{equation}

\noindent where $q:\mathbb{R}^{m}\times\mathbb{R}^{m\times d_{x}}\rightarrow\mathbb{R}^{d_{x}}$
is the thermal heat flux, $\tau:\mathbb{R}^{m}\times\mathbb{R}^{m\times d_{x}}\rightarrow\mathbb{R}^{d_{x}\times d_{x}}$
is the viscous stress tensor. The $i$-th spatial component of the
thermal heat flux is given by

\noindent 
\begin{equation}
q_{i}=-k\frac{\partial T}{\partial x_{i}},\label{eq:navier-stokes-thermal-heat-flux-component}
\end{equation}
where $T:\mathbb{R}^{m}\rightarrow\mathbb{R}_{+}$ is the temperature
and $k$ is thermal conductivity. The temperature $T$ is defined
as

\begin{equation}
T=\frac{p}{R\rho},\label{eq:temperature}
\end{equation}

\noindent where $R=287$ is the mixed specific gas constant for air.
The $i$-th spatial component of the viscous stress tensor is given
by

\noindent 
\begin{equation}
\tau_{i}=\mu\left(\frac{\partial v_{1}}{\partial x_{i}}+\frac{\partial v_{i}}{\partial x_{1}}-\delta_{i1}\frac{2}{3}\sum_{j=1}^{d_{x}}\frac{\partial v_{j}}{\partial x_{j}},\ldots,\frac{\partial v_{d_{x}}}{\partial x_{i}}+\frac{\partial v_{i}}{\partial x_{d_{x}}}-\delta_{id_{x}}\frac{2}{3}\sum_{j=1}^{d_{x}}\frac{\partial v_{j}}{\partial x_{j}}\right),\label{eq:navier-stokes-viscous-stress-tensor-component}
\end{equation}

\noindent where $\mu$ is the dynamic viscosity coefficient.

\subsection{Formulation in physical space with fixed geometry\label{subsec:formulation-physical-space}}

We consider the first-order system corresponding to convection-diffusion
that includes a conservation law, constitutive law, and the corresponding
interface conditions. The viscous interface conditions, more specifically
the viscous contribution to~(\ref{eq:interface-condition-strong-viscous})
and~(\ref{eq:interface-condition-state-strong-viscous}), were derived
in Section~2.2 of our earlier work~\citep{Ker20}.

\subsubsection{Strong formulation\label{subsec:strong-form}~\citep[Section 2.3.1]{Ker20}}

We consider the strong formulation presented in Section~2.3.1 of
our earlier work~\citep{Ker20} involving a nonlinear conservation
law, generalized constitutive law, and their corresponding interface
conditions,

\begin{align}
\nabla\cdot\mathcal{F}\left(y,\sigma\right)-f=0 & \textup{ in }\kappa\qquad\forall\kappa\in\mathcal{T},\label{eq:conservation-strong-viscous}\\
\sigma-G(y)\nabla_{x}y=0 & \textup{ in }\kappa\qquad\forall\kappa\in\mathcal{T},\label{eq:constitutive-strong-viscous}\\
\left\llbracket n\cdot\mathcal{F}\left(y,\sigma\right)\right\rrbracket =0 & \textup{ on }\epsilon\qquad\forall\epsilon\in\mathcal{E},\label{eq:interface-condition-strong-viscous}\\
\average{G\left(y\right)}\left\llbracket y\otimes n_{x}\right\rrbracket =0 & \textup{ on }\epsilon\qquad\forall\epsilon\in\mathcal{E},\label{eq:interface-condition-state-strong-viscous}
\end{align}
governing the flow state variable $y$ and auxiliary variable $\sigma$.
The interface condition~(\ref{eq:interface-condition-strong-viscous})
corresponding to the conservation law~(\ref{eq:conservation-strong-viscous})
is the jump or Rankine-Hugoniot condition~\citep{Maj12}, which accounts
for both the convective and viscous flux, ensuring continuity of the
normal flux at the interface. The interface condition~(\ref{eq:interface-condition-state-strong-viscous})
corresponds to the constitutive law~(\ref{eq:constitutive-strong-viscous})
and enforces a constraint on the continuity of the state variable
at the interface. 

The flux $\mathcal{F}\left(y,\sigma\right)$ is now defined
\begin{equation}
\mathcal{F}\left(y,\sigma\right)=\mathcal{F}^{c}\left(y\right)-\mathcal{\tilde{F}}^{v}\left(y,\sigma\right),\label{eq:flux-definition-in-terms-of-auxiliary-variable-flux-formulation}
\end{equation}
in terms of the convective flux and a modified viscous flux $\tilde{F}^{v}:\mathbb{R}^{m}\times\mathbb{R}^{m\times d_{x}}\rightarrow\mathbb{R}^{m\times d_{x}}$
. The modified viscous flux is defined consistently with the primal
formulation of Section~\ref{subsec:Governing-equations-Primal-Formulation},
as
\begin{equation}
\mathcal{F}^{v}\left(y,\nabla_{x}y\right)=\mathcal{\tilde{F}}^{v}\left(y,G\left(y\right)\nabla_{x}y\right),\label{eq:modified-viscous-flux-consistency}
\end{equation}
where $G\left(y\right)\in\mathbb{R}^{m\times d_{x}\times m\times d_{x}}$
is a generalized \emph{constitutive tensor} that depends on the specific
choice of constitutive law, a detailed discussion of which is provided
in our earlier work~\citep{Ker20}. In this work, we continue to
use the \emph{flux formulation}~\citep{Ker20} in the case of linear
advection-diffusion and viscous Burgers flow, where the constitutive
tensor is defined so that
\begin{equation}
G\left(y\right)\nabla_{x}y=\mathcal{F}^{v}\left(y,\nabla_{x}y\right)=\mathcal{F}_{\nabla_{x}y}^{v}\left(y,\nabla_{x}y\right)\nabla_{x}y,\label{eq:flux-formulation-tensor}
\end{equation}
and the modified viscous flux is defined to be the auxiliary variable,
\begin{equation}
\mathcal{\tilde{F}}^{v}\left(y,\sigma\right)=\sigma.\label{eq:flux-formulation-modified-viscous-flux}
\end{equation}

In the case of compressible Navier-Stokes flow, we define $G\left(y\right)\nabla_{x}y$,
up to a factor $\mu_{\infty}^{-\nicefrac{1}{2}}$, as the viscous
stress tensor, $\tau$, given by~(\ref{eq:navier-stokes-viscous-stress-tensor-component})
and thermal heat flux, $q$, given by~(\ref{eq:navier-stokes-thermal-heat-flux-component}).,
which is an approach similar to that of Chan et al.~\citep{Cha14}
while also incorporating the scaling advocated by Broersen and Stevenson~\citep{Bro14,Bro15}
and later Demkowicz and Gopalakrishnan~\citep{Dem15_20}. Thus, the
constitutive tensor $G\left(y\right)$ is defined such that
\begin{equation}
\left(G\left(y\right)\nabla y\right)_{i}=\mu_{\infty}^{-\nicefrac{1}{2}}\left(0,\tau_{1i},\ldots,\tau_{d_{x}i},-q_{i}\right),\label{eq:navier-stokes-auxiliary-component}
\end{equation}
where $\mu_{\infty}$ is the freestream dynamic viscosity and the
viscous flux is defined in terms of the auxiliary variable as
\begin{equation}
\mathcal{F}_{i}^{v}\left(y,\sigma\right)=\mu_{\infty}^{\nicefrac{1}{2}}\left(\sigma_{1i},\sigma_{2i},\ldots,\sigma_{d_{x}+1i},\sigma_{i+1j}v_{j}+\sigma_{mi}\right).\label{eq:navier-stokes-viscous-flux-auxiliary-component}
\end{equation}

\subsubsection{Interior and boundary interfaces\label{subsec:Interior-and-Boundary}~\citep[Section 2.3.2]{Ker20}}

We  assume that $\mathcal{E}$ consists of two disjoint subsets: the
interior interfaces $\mathcal{E}_{0}=\left\{ \epsilon_{0}\in\mathcal{E}\middlebar\epsilon_{0}\cap\partial\Omega=\emptyset\right\} $
and exterior interfaces $\mathcal{E}_{\partial}=\left\{ \epsilon_{\partial}\in\mathcal{E}\middlebar\epsilon_{\partial}\subset\partial\Omega\right\} ,$
so that $\mathcal{E}=\mathcal{E}_{0}\cup\mathcal{E}_{\partial}$.
For interior interfaces $\epsilon_{0}\in\mathcal{E}_{0}$ there exists
$\kappa^{+},\kappa^{-}\in\mathcal{T}$ such that $\epsilon_{0}=\partial\kappa^{+}\cap\partial\kappa^{-}$.
On interior interfaces Equations~(\ref{eq:interface-condition-strong-viscous}),~(\ref{eq:interface-condition-state-strong-viscous})
are defined as

\begin{align}
\left\llbracket n\cdot\mathcal{F}\left(y,\sigma\right)\right\rrbracket =n^{+}\cdot\mathcal{F}\left(y^{+},\sigma^{+}\right)+n^{-}\cdot\mathcal{F}\left(y^{-},\sigma^{-}\right)=0, & \textup{ on }\epsilon\qquad\forall\epsilon\in\mathcal{E}_{0},\label{eq:interface-condition-strong-interior}\\
\average{G\left(y\right)}\left\llbracket y\otimes n_{x}\right\rrbracket =\frac{1}{2}\left(G\left(y^{+}\right)+G\left(y^{-}\right)\right)\left(y^{+}\otimes n_{x}^{+}+y^{-}\otimes n_{x}^{-}\right)=0, & \textup{ on }\epsilon\qquad\forall\epsilon\in\mathcal{E}_{0}.\label{eq:interface-condition-state-strong-interior}
\end{align}
where $n^{+},n^{-}$ denote the outward facing normal of $\kappa^{+},\kappa^{-}$
respectively, so that $n^{+}=-n^{-}$. For exterior interfaces

\begin{align}
\left\llbracket n\cdot\mathcal{F}\left(y,\sigma\right)\right\rrbracket =n^{+}\cdot\mathcal{F}\left(y^{+},\sigma^{+}\right)-n^{+}\cdot\mathcal{\mathcal{F}}_{\partial}\left(y^{+},\sigma^{+}\right)=0, & \textup{ on }\epsilon\qquad\forall\epsilon\in\mathcal{E}_{\partial},\label{eq:interface-condition-strong-exterior}\\
\average{G\left(y\right)}\left\llbracket y\otimes n_{x}\right\rrbracket =G_{\partial}\left(y^{+}\right)\left(y^{+}\otimes n_{x}^{+}-y_{\partial}\left(y^{+}\right)\otimes n_{x}^{+}\right)=0, & \textup{ on }\epsilon\qquad\forall\epsilon\in\mathcal{E}_{\partial}.\label{eq:interface-condition-state-strong-exterior}
\end{align}

Here $n^{+}\cdot\mathcal{F}_{\partial}\left(y^{+},\sigma^{+}\right)$
is the imposed normal boundary flux, $G_{\partial}\left(y^{+}\right)$
is the boundary modified homogeneity tensor, and $y_{\partial}\left(y^{+}\right)$
is the boundary state, which are functions chosen depending on the
type of boundary condition. Therefore, we further decompose $\mathcal{E}_{\partial}$
into disjoint subsets of inflow and outflow interfaces $\mathcal{E}_{\partial}=\mathcal{E}_{\text{in}}\cup\mathcal{E}_{\text{out}},$
so that at an outflow interface $\epsilon_{\textup{out}}$ the boundary
flux is defined as the interior convective flux, and the boundary
state is defined as the interior state,
\begin{align}
n^{+}\cdot\mathcal{F}_{\partial}\left(y^{+},\sigma^{+}\right)=n^{+}\cdot\mathcal{F}\left(y^{+},\sigma_{\text{out}}=0\right), & \textup{ on }\epsilon\qquad\forall\epsilon\in\mathcal{E}_{\text{out}},\label{eq:normal-boundary-flux-outflow}\\
G_{\partial}\left(y^{+}\right)=G\left(y^{+}\right), & \textup{ on }\epsilon\qquad\forall\epsilon\in\mathcal{E}_{\text{out}},\label{eq:boundary-tensor-outflow}\\
y_{\partial}\left(y^{+}\right)=y^{+}, & \textup{ on }\epsilon\qquad\forall\epsilon\in\mathcal{E}_{\text{out}},\label{eq:boundary-state-outflow}
\end{align}
and therefore~(\ref{eq:interface-condition-state-strong-exterior})
is satisfied trivially. At an inflow boundary $\epsilon_{\textup{in}}\in\mathcal{E}_{\text{in}}$,
the normal convective boundary flux and boundary state are prescribed
values independent of the interior state $y^{+}$ , while the normal
viscous boundary flux is defined as the interior normal viscous flux,
\begin{align}
n^{+}\cdot\mathcal{\mathcal{F}}_{\partial}\left(y^{+},\sigma^{+}\right)=n^{+}\cdot\mathcal{\mathcal{F}}_{\text{in}}^{c}-n^{+}\cdot\mathcal{\tilde{F}}^{v}\left(y_{\partial}\left(y^{+}\right),\sigma^{+}\right), & \textup{ on }\epsilon\qquad\forall\epsilon\in\mathcal{E}_{\text{in}},\label{eq:normal-boundary-flux-inflow}\\
G_{\partial}\left(y^{+}\right)=G\left(y_{\partial}\left(y^{+}\right)\right), & \textup{ on }\epsilon\qquad\forall\epsilon\in\mathcal{E}_{\text{in}},\label{eq:boundary-tensor-inflow}\\
y_{\partial}\left(y^{+}\right)=y_{\text{in}}, & \textup{ on }\epsilon\qquad\forall\epsilon\in\mathcal{E}_{\text{in}}.\label{eq:boundary-state-inflow}
\end{align}

\subsubsection{Least-squares weak formulation\label{subsec:dpg-weak-form}}

We now describe a least-squares weak formulation of MDG-ICE in physical
space. Before proceeding, we note that special care is required to
weight the interface condition terms in the weak formulation in physical
space. If $h$ represents an element length scale, then a scaling
$h^{-1}$ arises in an inverse inequality that applies to functions
$v$ in a finite-dimensional subspace $V_{h}$, 
\begin{equation}
\left(v,v\right)_{H^{1/2}\left(\epsilon\right)}=\left\Vert v\right\Vert _{H^{1/2}\left(\epsilon\right)}^{2}\leq Ch^{-1}\left\Vert v\right\Vert _{L^{2}\left(\epsilon\right)}^{2},\label{eq:fractional-inverse-inequality}
\end{equation}
where $C$ is a positive constant independent of $h$, cf.~the works
of Bochev and Gunzberger~\citep[(4.22)--(4.23)]{Boc98}, Cao and
Gunzburger~\citep[(3.6)]{Cao98}, and Guermond~\citep[(3.3)]{Gue09}.
Practically speaking, the element integrals that arise in LS-MDG-ICE
over $\kappa$ scale proportionally to $h^{d-2}$ while the interface
integrals over $\epsilon$ scale proportionally to $h^{d-1}$. This
scaling has been accounted for in other least-squares finite element
methods involving interface or jump condition enforcement, cf. the
works of Cao and Gunzburger~\citep{Cao98} and Gerritsma and Proot~\citep{Ger02}.
We present an LS-MDG-ICE formulation in physical space in order to
illustrate the methodology that will be used to derive the reference
space formulation, and to this end, we neglect the scaling. However,
the scaling \emph{will} be naturally accounted for in the process
of transforming to a reference space least-squares formulation in
Section~\ref{subsec:variable-geometry}. Therefore, the results presented
in Section~\ref{sec:Examples} will use the reference space formulation
of Section~\ref{subsec:variable-geometry} even in the case of a
static grid.

The solution spaces $Y$ and $\Sigma$ are the broken Sobolev spaces,
\begin{eqnarray}
Y & = & \left\{ y\in\left[L^{2}\left(\Omega\right)\right]^{m\hphantom{\times d_{x}}}\bigl|\forall\kappa\in\mathcal{T},\:\:\hphantom{\nabla_{x}\cdot}\left.y\right|_{\kappa}\in\left[H^{1}\left(\kappa\right)\right]^{m}\right\} ,\label{eq:solution-space-vector}\\
\Sigma & = & \left\{ \sigma\in\left[L^{2}\left(\Omega\right)\right]^{m\times d_{x}}\bigl|\forall\kappa\in\mathcal{T},\left.\nabla_{x}\cdot\sigma\right|_{\kappa}\in\left[\:L^{2}\left(\kappa\right)\right]^{m}\right\} ,\label{eq:solution-space-tensor}
\end{eqnarray}
defined over a mesh $\mathcal{T}$. The test spaces are defined as
$V_{y}=\left[L^{2}\left(\Omega\right)\right]^{m}$ and $V_{\sigma}=\left[L^{2}\left(\Omega\right)\right]^{m\times d_{x}}$,
with $W_{y}$ and $W_{\sigma}$ defined to be the corresponding single-valued
trace spaces, cf. Carstensen et al.~\citep{Car16}. We integrate~(\ref{eq:conservation-strong-viscous})-(\ref{eq:interface-condition-state-strong-viscous})
on each element and interface against separate test functions to define
a nonlinear state operator,
\begin{equation}
e:Y\times\Sigma\rightarrow\left(V_{y}\times V_{\sigma}\times W_{y}\times W_{\sigma}\right)^{*}\label{eq:state-equation-physical-space}
\end{equation}
where,
\begin{align}
\left\langle e\left(y,\sigma\right),\left(v,\tau,w_{y},w_{\sigma}\right)\right\rangle = & \;\;\;\;\sum_{\kappa\in\mathcal{T}}\left(\nabla\cdot\mathcal{F}\left(y,\sigma\right)-f,v\right)_{\kappa}\nonumber \\
 & +\sum_{\kappa\in\mathcal{T}}\left(\sigma-G(y)\nabla_{x}y,\tau\right)_{\kappa}\nonumber \\
 & -\sum_{\epsilon\in\mathcal{E}}\left(\left\llbracket n\cdot\mathcal{F}\left(y,\sigma\right)\right\rrbracket ,w_{y}\right)_{\epsilon}\nonumber \\
 & -\sum_{\epsilon\in\mathcal{E}}\left(\average{G\left(y\right)}\left\llbracket y\otimes n_{x}\right\rrbracket ,w_{\sigma}\right)_{\epsilon},\label{eq:state-equation-physical-definition}
\end{align}
for $\left(y,\sigma\right)\in Y\times\Sigma$ and $\left(v,\tau,w_{y},w_{\sigma}\right)\in V_{y}\times V_{\sigma}\times W_{y}\times W_{\sigma}$.
Due to the nonlinear state operator, the trial-to-test operator is
defined in terms of the linearization of the state operator, cf.~\citep{Cha14,Bui14}.
Given $\left(y,\sigma\right)\in Y\times\Sigma$, we have
\begin{equation}
T\left(y,\sigma\right):Y\times\Sigma\rightarrow V_{y}\times V_{\sigma}\times W_{y}\times W_{\sigma},\label{eq:trial-to-test-operator-spaces-viscous}
\end{equation}
defined for $\left(v_{y},v_{\sigma}\right)\in Y\times\Sigma$ by
\begin{equation}
T\left(y,\sigma\right)=\left(v_{y},v_{\sigma}\right)\mapsto R_{V}^{-1}e'\left(y,\sigma\right)\left(v_{y},v_{\sigma}\right),\label{eq:trial-to-test-operator}
\end{equation}
where make use of Riesz map notation, cf.~\citep{Dem11},
\begin{equation}
R_{V}=R_{V_{y}\times V_{\sigma}\times W_{y}\times W_{\sigma}}:\left(V_{y}\times V_{\sigma}\times W_{y}\times W_{\sigma}\right)\rightarrow\left(V_{y}\times V_{\sigma}\times W_{y}\times W_{\sigma}\right)^{*}.\label{eq:riesz-map}
\end{equation}
Since each test space is equipped with an $L^{2}$ inner product,
we can write the trial-to-test-operator as
\begin{equation}
T\left(y,\sigma\right):\left(v_{y},v_{\sigma}\right)\mapsto\left(v^{\text{opt}},\tau^{\text{opt}},w_{y}^{\text{opt}},w_{\sigma}^{\text{opt}}\right),\label{eq:trial-to-test-physical}
\end{equation}
with
\begin{align}
v^{\text{opt}}= & \nabla\cdot\left(\mathcal{F}'\left(y,\sigma\right)\left(v_{y},v_{\sigma}\right)\right),\nonumber \\
\tau^{\text{opt}}= & v_{\sigma}-\left(\left(G'\left(y\right)v_{y}\right)\nabla_{x}y+G\left(y\right)\nabla_{x}v_{y}\right),\nonumber \\
w_{y}^{\text{opt}}= & \left\llbracket n\cdot\left(\mathcal{F}'\left(y,\sigma\right)\left(v_{y},v_{\sigma}\right)\right)\right\rrbracket ,\nonumber \\
w_{\sigma}^{\text{opt}}= & \average{G'\left(y\right)v_{y}}\left\llbracket y\otimes n_{x}\right\rrbracket +\average{G\left(y\right)}\left\llbracket v_{y}\otimes n_{x}\right\rrbracket .\label{eq:trial-to-test-operator-viscous}
\end{align}
Therefore, the DPG/least-squares formulation is given as: find $\left(y,\sigma\right)\in Y\times\Sigma$
such that

\begin{align}
0= & \;\;\;\;\sum_{\kappa\in\mathcal{T}}\left(\nabla\cdot\mathcal{F}\left(y,\sigma\right)-f,v^{\text{opt}}\right)_{\kappa}\nonumber \\
 & +\sum_{\kappa\in\mathcal{T}}\left(\sigma-\left(G\left(y\right)\nabla_{x}y\right),\tau^{\text{opt}}\right)_{\kappa}\nonumber \\
 & -\sum_{\epsilon\in\mathcal{E}}\left(\left\llbracket n\cdot\mathcal{F}\left(y,\sigma\right)\right\rrbracket ,w_{y}^{\text{opt}}\right)_{\epsilon}\nonumber \\
 & -\sum_{\epsilon\in\mathcal{E}}\left(\average{G\left(y\right)}\left\llbracket y\otimes n_{x}\right\rrbracket ,w_{\sigma}^{\text{opt}}\right)_{\epsilon}\qquad\forall\left(v_{y},v_{\sigma}\right)\in Y\times\Sigma,\label{eq:dpg-weak-form-viscous}
\end{align}
or abstractly as: find $\left(y,\sigma\right)\in Y\times\Sigma$ such
that
\begin{equation}
\left\langle e\left(y,\sigma\right),R_{V}^{-1}\left(e'\left(y,\sigma\right)\left(v_{y},v_{\sigma}\right)\right)\right\rangle =0\qquad\forall\left(v_{y},v_{\sigma}\right)\in Y\times\Sigma.\label{eq:dpg-weak-form-concise}
\end{equation}
where we note that the scaling~(\ref{eq:fractional-inverse-inequality})
is not yet accounted for.

\subsection{Formulation in reference space with variable geometry\label{subsec:variable-geometry}}

In this section, we transform the strong and weak formulation to reference
space, leading to a DPG/least-squares formulation of MDG-ICE, while
introducing the grid as a variable in order to fit interfaces and
resolve initially under-resolved solution features. We assume that
there is a continuous, invertible mapping
\begin{equation}
u:\hat{\Omega}\rightarrow\Omega,\label{eq:shape-mapping}
\end{equation}
from a reference domain $\hat{\Omega}\subset\mathbb{R}^{d}$ to the
physical domain $\Omega\subset\mathbb{R}^{d}$. We assume that $\hat{\Omega}$
is partitioned by $\hat{\mathcal{T}}$, so that $\overline{\hat{\Omega}}=\cup_{\hat{\kappa}\in\hat{\mathcal{T}}}\overline{\hat{\kappa}}$.
Also, we consider the set of interfaces $\hat{\mathcal{E}}$ consisting
of disjoint interfaces $\hat{\epsilon}$, such that $\cup_{\hat{\epsilon}\in\hat{\mathcal{E}}}\hat{\epsilon}=\cup_{\hat{\kappa}\in\hat{\mathcal{T}}}\partial\hat{\kappa}$.
The strong form in reference space~\citep{Ker20} is

\begin{align}
\left(\cof\left(\nabla u\right)\nabla\right)\cdot\mathcal{F}\left(y,\sigma\right)-\det\left(\nabla u\right)f=0 & \textup{ in }\hat{\kappa}\qquad\forall\hat{\kappa}\in\hat{\mathcal{T}},\label{eq:conservation-strong-reference-viscous}\\
\det\left(\nabla u\right)\sigma-G(y)\left(\cof\left(\nabla u\right)\nabla\right)_{x}y=0 & \textup{ in }\hat{\kappa}\qquad\forall\hat{\kappa}\in\hat{\mathcal{T}},\label{eq:constitutive-strong-reference-viscous}\\
\left\llbracket s\left(\nabla u\right)\cdot\mathcal{F}\left(y,\sigma\right)\right\rrbracket =0 & \textup{ on }\hat{\epsilon}\qquad\forall\hat{\epsilon}\in\hat{\mathcal{E}},\label{eq:interface-condition-strong-reference-viscous}\\
\average{G\left(y\right)}\left\llbracket y\otimes s\left(\nabla u\right)_{x}\right\rrbracket =0 & \textup{ on }\hat{\epsilon}\qquad\forall\hat{\epsilon}\in\hat{\mathcal{E}},\label{eq:interface-condition-state-strong-reference-viscous}\\
b\left(u\right)-u=0 & \textup{ on }\hat{\epsilon}\qquad\forall\hat{\epsilon}\in\hat{\mathcal{E}}.\label{eq:geometric-boundary-condition}
\end{align}
where $\nabla u$ is the Jacobian of the mapping from reference to
physical space, $\det\left(\nabla u\right)$ is the determinant of
the Jacobian, and $\cof\left(\nabla u\right)=\det\left(\nabla u\right)\left(\nabla u\right)^{-\top}$
is the cofactor matrix. The scaled normal $s\left(\nabla u\right)$
is the physical space normal scaled by the magnitude of generalized
cross product of the tangent plane basis vectors of a parameterization
of the physical space interface, as detailed in previous work~\citep{Cor18,Ker20}.
On elements $\kappa$ and interfaces $\epsilon$ with length scale
$h$, $\left(\cof\left(\nabla u\right)\nabla\right)$ and $s\left(\nabla u\right)$
both scale proportionally to $h^{d-1}$ accounting for the scaling
discussed in Section~\ref{subsec:dpg-weak-form}. In a spatial setting,
$\left(\cof\left(\nabla u\right)\nabla\right)_{x}=\left(\cof\left(\nabla u\right)\nabla\right)$
and $s\left(\nabla u\right)_{x}=s\left(\nabla u\right)$. Equation~(\ref{eq:geometric-boundary-condition})
constrains points to the boundary of the physical domain via a projection
operator $b:U\rightarrow U$, where the space of mappings from reference
space to physical space, $U=\left[H^{1}\left(\hat{\Omega}\right)\right]^{d},$
is the $\mathbb{R}^{d}$-valued Sobolev space over $\hat{\Omega}$,
cf. the previous work~\citep{Cor18,Ker20}. We assume that $Y$ and
$\Sigma$ now consist of functions defined in $\mathbb{R}^{m}$-valued
and $\mathbb{R}^{m\times d_{x}}$-valued broken Sobolev spaces over
$\hat{\mathcal{T}}$ respectively. We further assume that the test
spaces $V_{y}=\left[L^{2}\left(\hat{\Omega}\right)\right]^{m}$ and
$V_{\sigma}=\left[L^{2}\left(\hat{\Omega}\right)\right]^{m\times d_{x}}$,
along with the single-valued trace test spaces, $W_{y}$ and $W_{\sigma}$,
now consist of functions defined over reference space.

As in our previous work~\citep{Ker20}, we define a provisional state
operator $\tilde{e}:Y\times\Sigma\times U\rightarrow\left(V_{y}\times V_{\sigma}\times W_{y}\times W_{\sigma}\right)^{*}$
for $\left(y,\sigma,u\right)\in Y\times\Sigma\times U$, by
\begin{align}
\tilde{e}\left(y,\sigma,u\right)=\left(v,\tau,w_{y},w_{\sigma}\right)\mapsto & \;\;\;\;\sum_{\hat{\kappa}\in\hat{\mathcal{T}}}\left(\left(\cof\left(\nabla u\right)\nabla\right)\cdot\mathcal{F}\left(y,\sigma\right)-\det\left(\nabla u\right)f,v\right)_{\hat{\kappa}}\nonumber \\
 & +\sum_{\hat{\kappa}\in\hat{\mathcal{T}}}\left(\det\left(\nabla u\right)\sigma-G(y)\left(\cof\left(\nabla u\right)\nabla\right)_{x}y,\tau\right)_{\hat{\kappa}}\nonumber \\
 & -\sum_{\hat{\epsilon}\in\hat{\mathcal{E}}}\left(\left\llbracket s\left(\nabla u\right)\cdot\mathcal{F}\left(y,\sigma\right)\right\rrbracket ,w_{y}\right)_{\hat{\epsilon}}\nonumber \\
 & -\sum_{\hat{\epsilon}\in\hat{\mathcal{E}}}\left(\average{G\left(y\right)}\left\llbracket y\otimes s\left(\nabla u\right)_{x}\right\rrbracket ,w_{\sigma}\right)_{\hat{\epsilon}}\label{eq:state-operator-weak-formulation-reference-viscous}
\end{align}
which has a Fréchet derivative defined for perturbation $\left(\delta y,\delta\sigma,\delta u\right)\in Y\times\Sigma\times U$,
and test functions $\left(v,\tau,w_{y},w_{\sigma}\right)\in V_{y}\times V_{\sigma}\times W_{y}\times W_{\sigma}$,
by its partial derivative with respect to the state variable $y$,
\begin{align}
\tilde{e}_{y}\left(y,\sigma,u\right)\delta y=\left(v,\tau,w_{y},w_{\sigma}\right)\mapsto & \;\;\;\;\sum_{\hat{\kappa}\in\hat{\mathcal{T}}}\left(\left(\cof\left(\nabla u\right)\nabla\right)\cdot\left(\mathcal{F}_{y}\left(y,\sigma\right)\delta y\right),v\right)_{\hat{\kappa}}\nonumber \\
 & -\sum_{\hat{\kappa}\in\hat{\mathcal{T}}}\left(\left(G'(y)\delta y\right)\left(\cof\left(\nabla u\right)\nabla\right)_{x}y+G(y)\left(\cof\left(\nabla u\right)\nabla\right)_{x}\delta y,\tau\right)_{\hat{\kappa}}\nonumber \\
 & -\sum_{\hat{\epsilon}\in\hat{\mathcal{E}}}\left(\left\llbracket s\left(\nabla u\right)\cdot\left(\mathcal{F}_{y}\left(y,\sigma\right)\delta y\right)\right\rrbracket ,w_{y}\right)_{\hat{\epsilon}}\nonumber \\
 & -\sum_{\hat{\epsilon}\in\hat{\mathcal{E}}}\left(\average{G'\left(y\right)\delta y}\left\llbracket y\otimes s\left(\nabla u\right)_{x}\right\rrbracket +\average{G\left(y\right)}\left\llbracket \delta y\otimes s\left(\nabla u\right)_{x}\right\rrbracket ,w_{\sigma}\right)_{\hat{\epsilon}},\label{eq:state-operator-state-derivative-weak-formulation-reference-viscous}
\end{align}
its partial derivative with respect to the auxiliary variable $\sigma$,
\begin{align}
\tilde{e}_{\sigma}\left(y,\sigma,u\right)\delta\sigma=\left(v,\tau,w_{y},w_{\sigma}\right)\mapsto & \;\;\;\;\sum_{\hat{\kappa}\in\hat{\mathcal{T}}}\left(\left(\cof\left(\nabla u\right)\nabla\right)\cdot\left(\mathcal{F}_{\sigma}\left(y,\sigma\right)\delta\sigma\right),v\right)_{\hat{\kappa}}\nonumber \\
 & +\sum_{\hat{\kappa}\in\hat{\mathcal{T}}}\left(\det\left(\nabla u\right)\delta\sigma,\tau\right)_{\hat{\kappa}}\nonumber \\
 & -\sum_{\hat{\epsilon}\in\hat{\mathcal{E}}}\left(\left\llbracket s\left(\nabla u\right)\cdot\left(\mathcal{F}_{\sigma}\left(y,\sigma\right)\delta\sigma\right)\right\rrbracket ,w_{y}\right)_{\hat{\epsilon}},\label{eq:state-operator-auxiliary-derivative-weak-formulation-reference-viscous}
\end{align}
and its partial derivative with respect to the geometry variable $u$,
\begin{align}
\tilde{e}_{u}\left(y,\sigma,u\right)\delta u=\left(v,\tau,w_{y},w_{\sigma}\right)\mapsto & \;\;\;\;\sum_{\hat{\kappa}\in\hat{\mathcal{T}}}\left(\left(\left(\cof'\left(\nabla u\right)\nabla\delta u\right)\nabla\right)\cdot\mathcal{F}\left(y,\sigma\right)-\left(\det'\left(\nabla u\right)\nabla\delta u\right)f,v\right)_{\hat{\kappa}}\nonumber \\
 & +\sum_{\hat{\kappa}\in\hat{\mathcal{T}}}\left(\left(\det'\left(\nabla u\right)\nabla\delta u\right)\sigma-G(y)\left(\left(\cof'\left(\nabla u\right)\nabla\delta u\right)\nabla\right)_{x}y,\tau\right)_{\hat{\kappa}}\nonumber \\
 & -\sum_{\hat{\epsilon}\in\mathcal{\hat{E}}}\left(\left\llbracket \left(s'\left(\nabla u\right)\nabla\delta u\right)\cdot\mathcal{F}\left(y,\sigma\right)\right\rrbracket ,w_{y}\right)_{\hat{\epsilon}}\nonumber \\
 & -\sum_{\hat{\epsilon}\in\hat{\mathcal{E}}}\left(\average{G\left(y\right)}\left\llbracket y\otimes\left(s'\left(\nabla u\right)\nabla\delta u\right)_{x}\right\rrbracket ,w_{\sigma}\right)_{\hat{\epsilon}}.\label{eq:state-operator-geometry-derivative-weak-formulation-reference-viscous}
\end{align}
In order to derive the LS-MDG-ICE formulation, we define the trial-to-test
operator
\begin{equation}
T\left(y,\sigma,u\right):Y\times\Sigma\times U\rightarrow V_{y}\times V_{\sigma}\times W_{y}\times W_{\sigma},\label{eq:trial-to-test-operator-spaces-viscous-reference}
\end{equation}
for $\left(v_{y},v_{\sigma},v_{u}\right)\in Y\times\Sigma\times U$
by
\begin{equation}
T\left(y,\sigma,u\right)=\left(v_{y},v_{\sigma},v_{u}\right)\mapsto R_{V}^{-1}\left(\tilde{e}'\left(y,\sigma,u\right)\left(v_{y},v_{\sigma},v_{u}\right)\right),\label{eq:trial-to-test-operator-1}
\end{equation}
where we continue to employ the inverse Riesz map~(\ref{eq:riesz-map}).
Analogous to the physical space formulation, since each test space
is equipped with an $L^{2}$ inner product, the trial-to-test operator
can be expressed as
\begin{equation}
T:\left(v_{y},v_{\sigma},v_{u}\right)\mapsto\left(v^{\text{opt}},\tau^{\text{opt}},w_{y}^{\text{opt}},w_{\sigma}^{\text{opt}}\right),\label{eq:trial-to-test-operator-viscous-reference-space}
\end{equation}
where
\begin{align}
v^{\text{opt}}= & \left(\cof\left(\nabla u\right)\nabla\right)\cdot\left(\mathcal{F}'\left(y,\sigma\right)\left(v_{y},v_{\sigma}\right)\right)\nonumber \\
 & +\left(\left(\cof'\left(\nabla u\right)\nabla v_{u}\right)\nabla\right)\cdot\mathcal{F}\left(y,\sigma\right)\nonumber \\
 & -\left(\det'\left(\nabla u\right)\nabla v_{u}\right)f,\label{eq:optimal-v}
\end{align}
\begin{align}
\tau^{\text{opt}}= & \det\left(\nabla u\right)v_{\sigma}+\left(\det'\left(\nabla u\right)\nabla v_{u}\right)\sigma\nonumber \\
 & -\left(G'\left(y\right)v_{y}\right)\left(\cof\left(\nabla u\right)\nabla\right)_{x}y\nonumber \\
 & -G\left(y\right)\left(\cof\left(\nabla u\right)\nabla\right)_{x}v_{y}\nonumber \\
 & -G\left(y\right)\left(\left(\cof'\left(\nabla u\right)\nabla v_{u}\right)\nabla\right)_{x}y,\label{eq:optimal-tau}
\end{align}
\begin{align}
w_{y}^{\text{opt}}= & \left\llbracket s\left(\nabla u\right)\cdot\left(\mathcal{F}'\left(y,\sigma\right)\left(v_{y},v_{\sigma}\right)\right)\right\rrbracket \nonumber \\
 & +\left\llbracket \left(s'\left(\nabla u\right)\nabla v_{u}\right)\cdot\mathcal{F}\left(y,\sigma\right)\right\rrbracket ,\label{optimal-w-y}
\end{align}
and
\begin{align}
w_{\sigma}^{\text{opt}}= & \average{G'\left(y\right)v_{y}}\left\llbracket y\otimes s\left(\nabla u\right)_{x}\right\rrbracket \nonumber \\
 & +\average{G\left(y\right)}\left\llbracket v_{y}\otimes s\left(\nabla u\right)_{x}\right\rrbracket \nonumber \\
 & +\average{G\left(y\right)}\left\llbracket y\otimes\left(s'\left(\nabla u\right)\nabla v_{u}\right)_{x}\right\rrbracket .\label{eq:optimal-w-sigma}
\end{align}
Using the trial-to-test operator (\ref{eq:trial-to-test-operator-viscous-reference-space}),
we can state a \emph{provisional, }in the sense that it does not yet
incorporate the geometric boundary condition~(\ref{eq:geometric-boundary-condition})\emph{,
}LS-MDG-ICE weak formulation as: find $\left(y,\sigma,u\right)\in Y\times\Sigma\times U$
such that

\begin{align}
 & \;\;\;\;\sum_{\hat{\kappa}\in\hat{\mathcal{T}}}\left(\left(\cof\left(\nabla u\right)\nabla\right)\cdot\mathcal{F}\left(y,\sigma\right)-\det\left(\nabla u\right)f,v^{\text{opt}}\right)_{\hat{\kappa}}\nonumber \\
 & +\sum_{\hat{\kappa}\in\hat{\mathcal{T}}}\left(\left(\det\left(\nabla u\right)\sigma\right)-\left(G\left(y\right)\left(\cof\left(\nabla u\right)\nabla\right)_{x}y\right),\tau^{\text{opt}}\right)_{\hat{\kappa}}\nonumber \\
 & -\sum_{\hat{\epsilon}\in\hat{\mathcal{E}}}\left(\left\llbracket s\left(\nabla u\right)\cdot\mathcal{F}\left(y,\sigma\right)\right\rrbracket ,w_{y}^{\text{opt}}\right)_{\hat{\epsilon}}\nonumber \\
 & -\sum_{\hat{\epsilon}\in\hat{\mathcal{E}}}\left(\average{G\left(y\right)}\left\llbracket y\otimes s\left(\nabla u\right)_{x}\right\rrbracket ,w_{\sigma}^{\text{opt}}\right)_{\hat{\epsilon}}\qquad\forall\left(v_{y},v_{\sigma},v_{u}\right)\in Y\times\Sigma\times U,\label{eq:state-operator-weak-formulation-reference-least-squares-1}
\end{align}
or abstractly as: find $\left(y,\sigma,u\right)\in Y\times\Sigma\times U$
such that
\begin{equation}
\left\langle \tilde{e}\left(y,\sigma,u\right),R_{V}^{-1}\left(\tilde{e}'\left(y,\sigma,u\right)\left(v_{y},v_{\sigma},v_{u}\right)\right)\right\rangle =0\qquad\forall\left(v_{y},v_{\sigma},v_{u}\right)\in Y\times\Sigma\times U.\label{eq:dpg-weak-form-concise-1}
\end{equation}
In order to incorporate the geometric boundary condition~(\ref{eq:geometric-boundary-condition}),
we compose the provisional weak formulation with the boundary projection
to obtain the weak formulation, find $\left(y,\sigma,u\right)\in Y\times\Sigma\times U$
such that
\begin{equation}
\left\langle \tilde{e}\left(y,\sigma,b\left(u\right)\right),R_{V}^{-1}\tilde{e}'\left(y,\sigma,b\left(u\right)\right)\left(v_{y},v_{\sigma},b'\left(u\right)v_{u}\right)\right\rangle =0\qquad\forall\left(v_{y},v_{\sigma},v_{u}\right)\in Y\times\Sigma\times U,\label{eq:dpg-weak-form-concise-1-1}
\end{equation}
so that the solution satisfying (\ref{eq:conservation-strong-reference-viscous})
and (\ref{eq:interface-condition-strong-reference-viscous}) weakly
and (\ref{eq:geometric-boundary-condition}) strongly is therefore
given as $\left(y,\sigma,b\left(u\right)\right)\in Y\times\Sigma\times U$.

\subsection{Discretization\label{subsec:Discretization}}

We consider a family of meshes, with discretization parameter $h$,
where each mesh is denoted $\hat{\mathcal{T}}_{h}$. We choose discrete
spaces $Y_{h}$, $\Sigma_{h}$, and $U_{h}$, and define a mesh-dependent
weak formulation of~(\ref{eq:dpg-weak-form-concise-1}) via the discrete
state operator,
\begin{equation}
\tilde{e}_{h}:Y_{h}\times\Sigma_{h}\times U_{h}\rightarrow T\left(Y_{h}\times\Sigma_{h}\times U_{h}\right).\label{eq:discrete-state-operator}
\end{equation}
The mesh-dependent state operator $\tilde{e}_{h}\left(y,\sigma,u\right)$
is defined according to (\ref{eq:state-operator-weak-formulation-reference-viscous}),
but with $\hat{\mathcal{T}}_{h}$ and $\hat{\mathcal{E}}_{h}$ substituted
for $\mathcal{\hat{T}}$ and $\hat{\mathcal{E}}$. As discussed in
Section~\ref{subsec:Optimal-Test-Spaces}, as the mesh is refined,
the interface condition is enforced over an increasingly dense set
of interfaces. We also discretize the boundary projection $b_{h}:U_{h}\rightarrow U_{h}$.
The discrete weak formulation is: find $\left(y,\sigma,u\right)\in Y_{h}\times\Sigma_{h}\times U_{h}$
such that
\begin{equation}
\left\langle \tilde{e}_{h}\left(y,\sigma,b_{h}\left(u\right)\right),R_{V_{h}}^{-1}\tilde{e}_{h}'\left(y,\sigma,b_{h}\left(u\right)\right)\left(v_{y},v_{\sigma},b_{h}'\left(u\right)v_{u}\right)\right\rangle =0\qquad\forall\left(v_{y},v_{\sigma},v_{u}\right)\in Y_{h}\times\Sigma_{h}\times U_{h},\label{eq:dpg-weak-form-concise-discretized}
\end{equation}
so that the discrete solution is given as $\left(y,\sigma,b_{h}\left(u\right)\right)\in Y_{h}\times\Sigma_{h}\times U_{h}$.

The discrete spaces are defined over the mesh $\mathcal{\hat{T}}_{h}$,
where $Y_{h}$ and $\Sigma_{h}$ are discontinuous finite element
spaces, while $U_{h}$ is a continuous finite element space. Let $\mathcal{P}_{p}$
denote the space of polynomials spanned by the monomials $\boldsymbol{x}^{\alpha}$
with multi-index $\alpha\in\mathbb{N}_{0}^{n}$ , satisfying $\sum_{i=1}^{n}\alpha_{i}\leq p$.
In the case of a simplicial grid,
\begin{eqnarray}
Y_{h} & = & \left\{ y\in\left[L^{2}\left(\hat{\Omega}\right)\right]^{m\hphantom{\times d_{x}}}\,\bigl|\forall\hat{\kappa}\in\hat{\mathcal{T}}_{h},\left.y\right|_{\hat{\kappa}}\in\left[\mathcal{P}_{p}\right]^{m\hphantom{\times d_{x}}}\right\} ,\label{eq:discrete-solution-state-space-simplex}\\
\Sigma_{h} & = & \left\{ \sigma\in\left[L^{2}\left(\hat{\Omega}\right)\right]^{m\times d_{x}}\bigl|\forall\hat{\kappa}\in\hat{\mathcal{T}}_{h},\left.\sigma\right|_{\hat{\kappa}}\in\left[\mathcal{P}_{p}\right]^{m\times d_{x}}\right\} .\label{eq:discrete-solution-flux-space-simplex}
\end{eqnarray}
The polynomial degree of the state space and flux space are in general
distinct.

The discrete space $U_{h}$ of mappings from reference space to physical
space are discretized into continuous, $\mathbb{R}^{d}$-valued piecewise
polynomials, in the case of a simplicial grid
\begin{equation}
U_{h}=\left\{ u\in\left[H^{1}\left(\hat{\Omega}\right)\right]^{d}\middlebar\forall\hat{\kappa}\in\hat{\mathcal{T}_{h}},\left.u\right|_{\hat{\kappa}}\in\left[\mathcal{P}_{p}\right]^{d}\right\} .\label{eq:discrete-shape-space-simplex}
\end{equation}
The case that the chosen polynomial degree of $U_{h}$ is equal to
that of $Y_{h}$ is referred to as isoparametric. It is also possible
to choose the polynomial degree of $U_{h}$ to be less (sub-parametric)
or greater (super-parametric) than that of $Y_{h}$. 

\subsection{Solver\label{subsec:minimum-residual-method}}

We employ an iterative, regularized Gauss-Newton solver to solve the
discretized nonlinear weak formulation~(\ref{eq:dpg-weak-form-concise-discretized}),
cf.~\citep{Cha14,Bui14}. Given an initialization $\left(y,\sigma,u\right)_{0}$
the solution is repeatedly updated
\begin{equation}
\left(y,\sigma,u\right)_{i+1}=\left(y,\sigma,u\right)_{i}+\left(\delta y,\delta\sigma,\delta u\right)_{i}\qquad i=0,1,2,\ldots,\label{eq:discrete-iteration}
\end{equation}
using an increment $\left(\delta y,\delta\sigma,\delta u\right)_{i}$
that is defined as the solution to a linear system of equations. We
employ a Gauss-Newton method that approximates the full linearization
of the weak formulation~(\ref{eq:dpg-weak-form-concise-discretized}),
neglecting second derivatives by only differentiating the term $e_{h}\left(y,\sigma,b_{h}\left(u\right)\right)$
in the dual pairing that appears in the weak formulation~(\ref{eq:dpg-weak-form-concise-discretized}).
The resulting bilinear form is symmetric and positive semi-definite.
To ensure rank-sufficiency, we incorporate a Levenberg-Marquardt regularization~\citep{Cor18}
via a positive-definite bilinear form. This leads to the linear problem:
find $\left(\delta y,\delta\sigma,\delta u\right)\in Y_{h}\times\Sigma_{h}\times U_{h}$
such that,
\begin{gather}
\left\langle e'_{h}\left(y,\sigma,b_{h}\left(u\right)\right)\left(\delta y,\delta\sigma,b'_{h}\left(u\right)\delta u\right),R_{V_{h}}^{-1}e'_{h}\left(y,\sigma,b_{h}\left(u\right)\right)\left(v_{y},v_{\sigma},b'_{h}\left(u\right)v_{u}\right)\right\rangle \nonumber \\
+I_{h,\lambda}\left(y,\sigma,u\right)\left(\left(\delta y,\delta\sigma,\delta u\right),\left(v_{y},v_{\sigma},v_{u}\right)\right)\\
=-\left\langle e_{h}\left(y,\sigma,b_{h}\left(u\right)\right),R_{V_{h}}^{-1}\left(e'_{h}\left(y,\sigma,b_{h}\left(u\right)\right)\left(v_{y},v_{\sigma},b'_{h}\left(u\right)v_{u}\right)\right)\right\rangle \qquad\forall\left(v_{y},v_{\sigma},v_{u}\right)\in Y_{h}\times\Sigma_{h}\times U_{h},\label{eq:regularized-gauss-newton-levenberg-marquardt}
\end{gather}
where $I_{h,\lambda}\left(y,\sigma,u\right):\left(Y_{h}\times\Sigma_{h}\times U_{h}\right)\times\left(Y_{h}\times\Sigma_{h}\times U_{h}\right)\rightarrow\mathbb{R}$
is a symmetric, positive-definite bilinear form that defines the choice
of regularization. An identity regularization defined for $\lambda_{y},\lambda_{\sigma},\lambda_{u}\ge0$
is given by 
\begin{equation}
I_{h,\lambda}\left(y,\sigma,u\right)\left(\left(\delta y,\delta\sigma,\delta u\right),\left(v_{y},v_{\sigma},v_{u}\right)\right)=\left(\delta y,\lambda_{y}v_{y}\right)+\left(\delta\sigma,\lambda_{\sigma}v_{\sigma}\right)+\left(\delta u,\lambda_{u}v_{u}\right).\label{eq:regularization-operator}
\end{equation}
Separate regularization coefficients $\lambda_{y},\lambda_{\sigma},\lambda_{u}\ge0$
are defined for each solution variable. In practice, the state and
auxiliary regularization coefficients can be set to zero, $\lambda_{y}=\lambda_{\sigma}=0$,
while the grid regularization coefficient $\lambda_{u}>0$ must be
positive in order to ensure rank sufficiency and to limit excessive
grid motion. Additional symmetric, positive (semi-)definite bilinear
forms can be incorporated into the regularization~\citep{Cor19_AVIATION},
for example
\begin{equation}
I_{h,\lambda}^{\Delta}\left(y,\sigma,u\right)\left(\left(\delta y,\delta\sigma,\delta u\right),\left(v_{y},v_{\sigma},v_{u}\right)\right)=-\left(\nabla\left(b_{h}'\left(u\right)\delta u\right),\lambda_{\Delta u}\nabla\left(b_{h}'\left(u\right)v_{u}\right)\right).\label{eq:regularization-operator-laplacian}
\end{equation}
The resulting linear system of equations is positive definite and
symmetric, which we solve using a sparse direct solver provided by
Eigen~\citep{Eig10}. 

Recently, Zahr et al.~\citep{Zah19} modified the weight $\lambda_{u}$
to include a factor proportional to the inverse of the element. For
the cases of supersonic and transonic flow over an airfoil, this weight
was found to be necessary for maintaining a valid increment during
the nonlinear solution procedure of their optimization based approach.
Without it, the line search would stall due to successive near zero
increments or produce excessively large increments and cause the grid
to become invalid. In Section~\ref{sec:viscous-bow-shock} we follow
this approach and confirm the observation of Zahr et al.~\citep{Zah19},
that regularization of this type can enhance the robustness of the
nonlinear solver for problems in which there is significant variation
in the size of the element.

Throughout the solution process, the discrete geometry may need to
be modified by the solver in order to ensure cell validity. As detailed
in our previous work~\citep{Cor18,Ker20}, we employ standard edge
refinement and edge collapse algorithms~\citep{Loh08}, which ensures
the interface topology of the partition $\hat{\mathcal{T}}$ remains
compatible.

\section{Examples\label{sec:Examples}}

The LS-MDG-ICE discretization is now applied to compute solutions
to both one-dimensional and two-dimensional problems corresponding
to linear advection-diffusion, viscous Burgers, inviscid Euler, and
compressible Navier-Stokes. We consider both the case of a static,
or fixed, discrete geometry and variable discrete geometry. Unless
otherwise noted, MDG-ICE solutions are assumed to include both the
discrete flow field and the discrete geometry. When clarification
is necessary, we will refer to the case of a fixed discrete geometry
as \emph{static} LS-MDG-ICE.

\subsection{Linear advection-diffusion in one dimension\label{sec:linear-advection-diffusion-1d}}

We consider the solution to the system corresponding to one-dimensional
linear advection-diffusion, subject to the following boundary conditions

\begin{align}
y\left(x=0\right) & =0\nonumber \\
y\left(x=1\right) & =1.\label{eq:steady-boundary-layer-boundary-conditions}
\end{align}

The exact solution is given by

\begin{equation}
y\left(x\right)=\frac{1-\exp\left(x\cdot\mathrm{Pe}\right)}{1-\exp\left(\mathrm{Pe}\right)}.\label{eq:steady-boundary-layer-exact-solution}
\end{equation}
where $\mathrm{Pe}=\frac{1}{\varepsilon}=\frac{v\ell}{\mu}$ is the
Péclet number, $v$ is the characteristic velocity, $\ell$ is the
characteristic length, and $\mu$ is the mass diffusivity.
\begin{figure}
\centering{}%
\begin{tabular}{cc}
\includegraphics[width=0.4\linewidth]{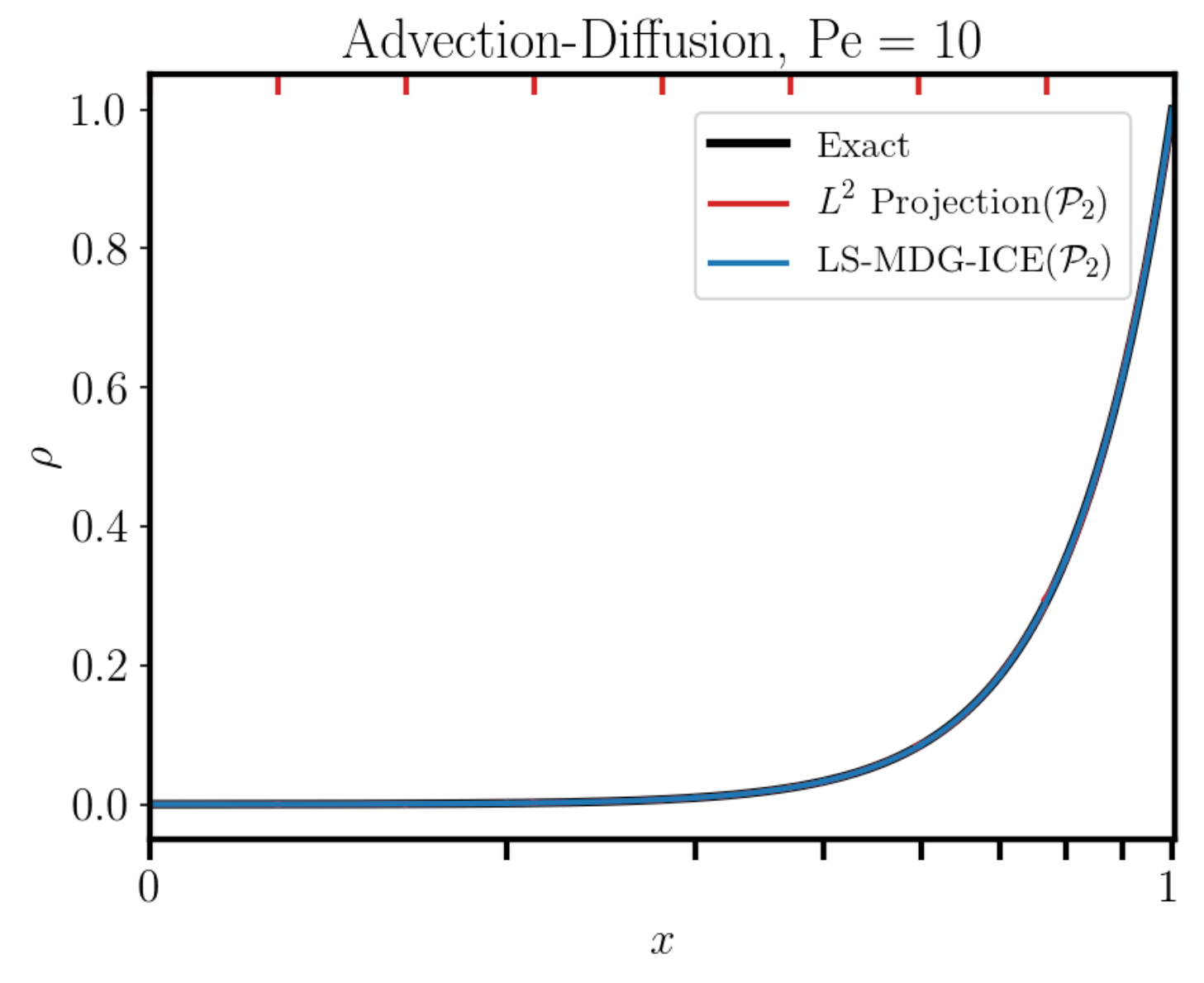} & \includegraphics[width=0.4\linewidth]{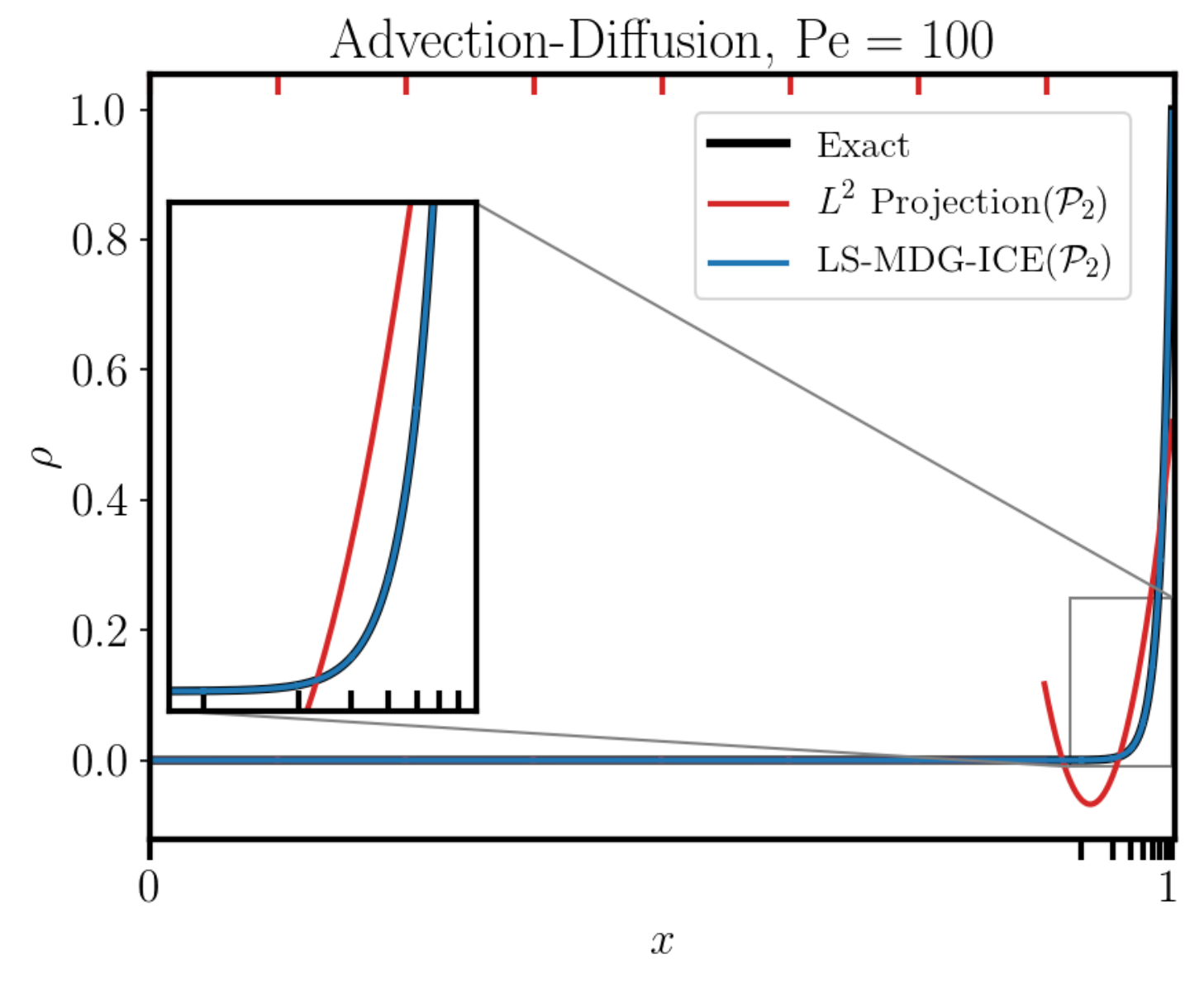}\tabularnewline
\includegraphics[width=0.4\linewidth]{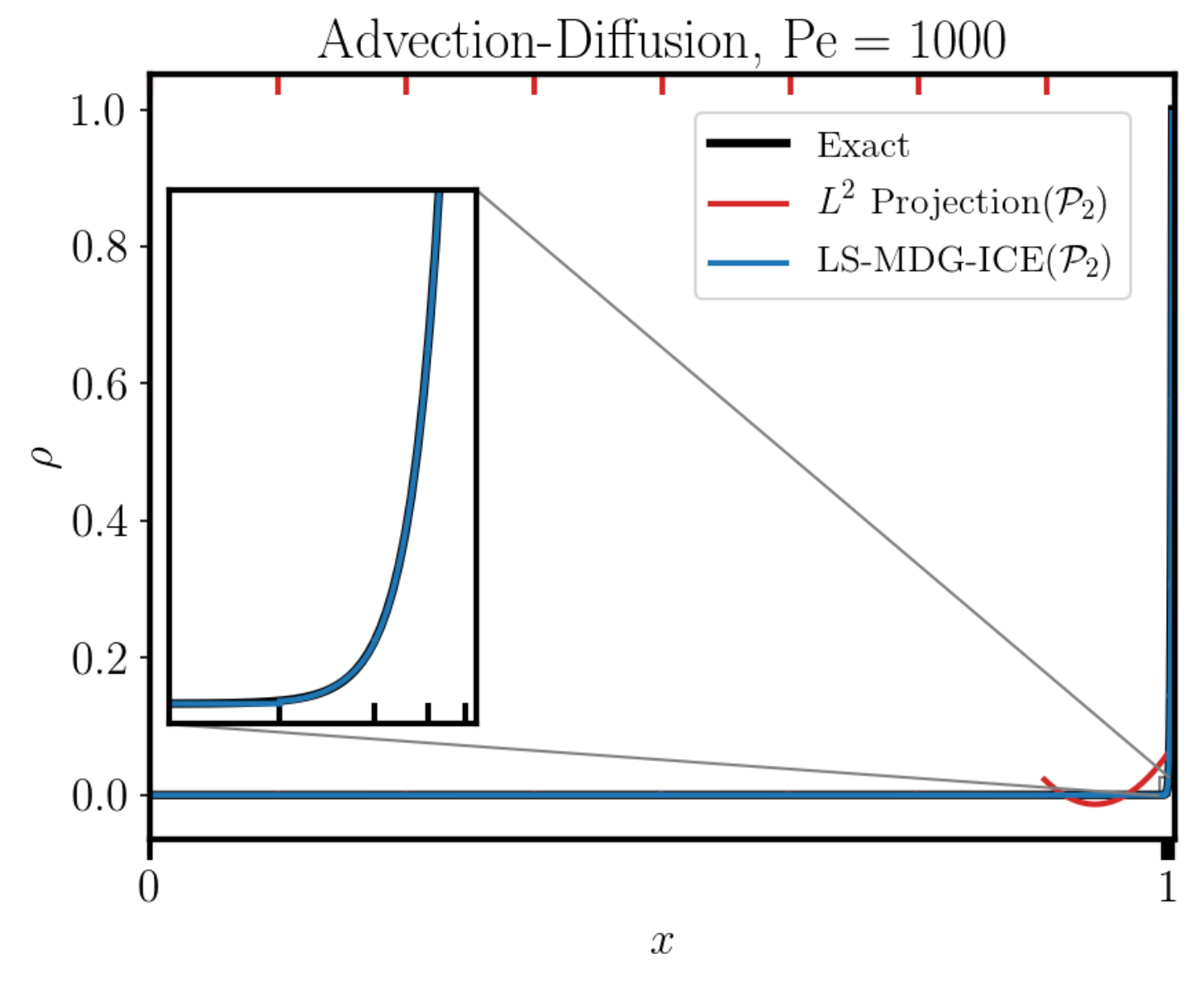} & \includegraphics[width=0.4\linewidth]{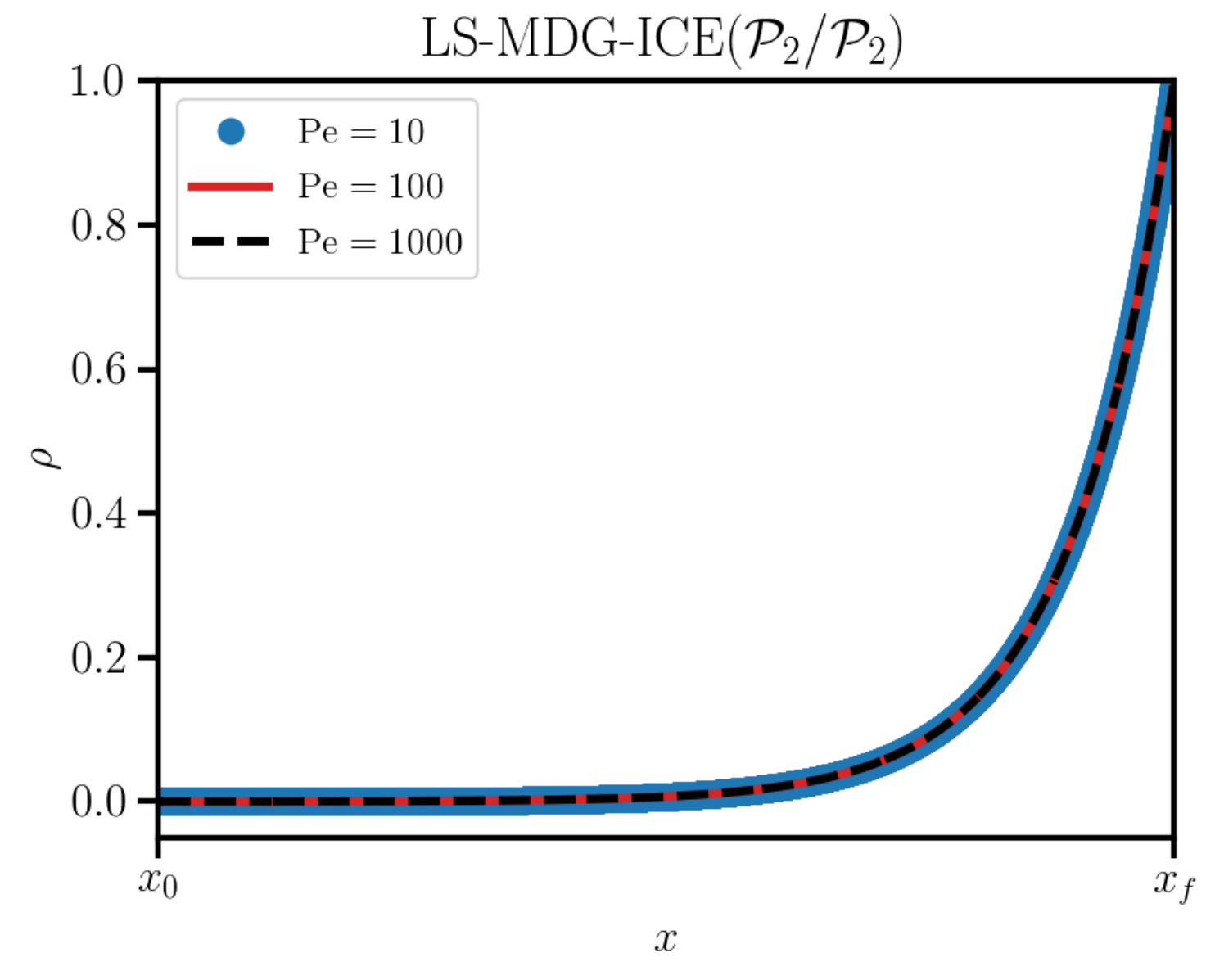}\tabularnewline
\end{tabular}\caption{\label{fig:steady-boundary-layer-p2}Grid adapted, self-similar solutions
computed with LS-MDG-ICE$\left(\mathcal{P}_{2}\right)$ for a one-dimensional
boundary layer type profile for $\mathrm{Pe}=10$, $100$, and $1000$
on a grid consisting of eight isoparametric line cells are compared
to the exact solution and the $L^{2}$ projection of the exact solution
onto a uniform grid. The solution was initialized with a linear profile
on a uniform grid. The grid markers correspond to the adapted grid
given by the LS-MDG-ICE solutions. The initial uniform grid is indicated
with red tick marks. The LS-MDG-ICE$\left(\mathcal{P}_{2}/\mathcal{P}_{2}\right)$
solutions for $\mathrm{Pe}=10$, $100$, and $1000$ plotted together
on rescaled domains (bottom-right).}
\end{figure}
\begin{figure}
\centering{}%
\begin{tabular}{cc}
\includegraphics[width=0.4\linewidth]{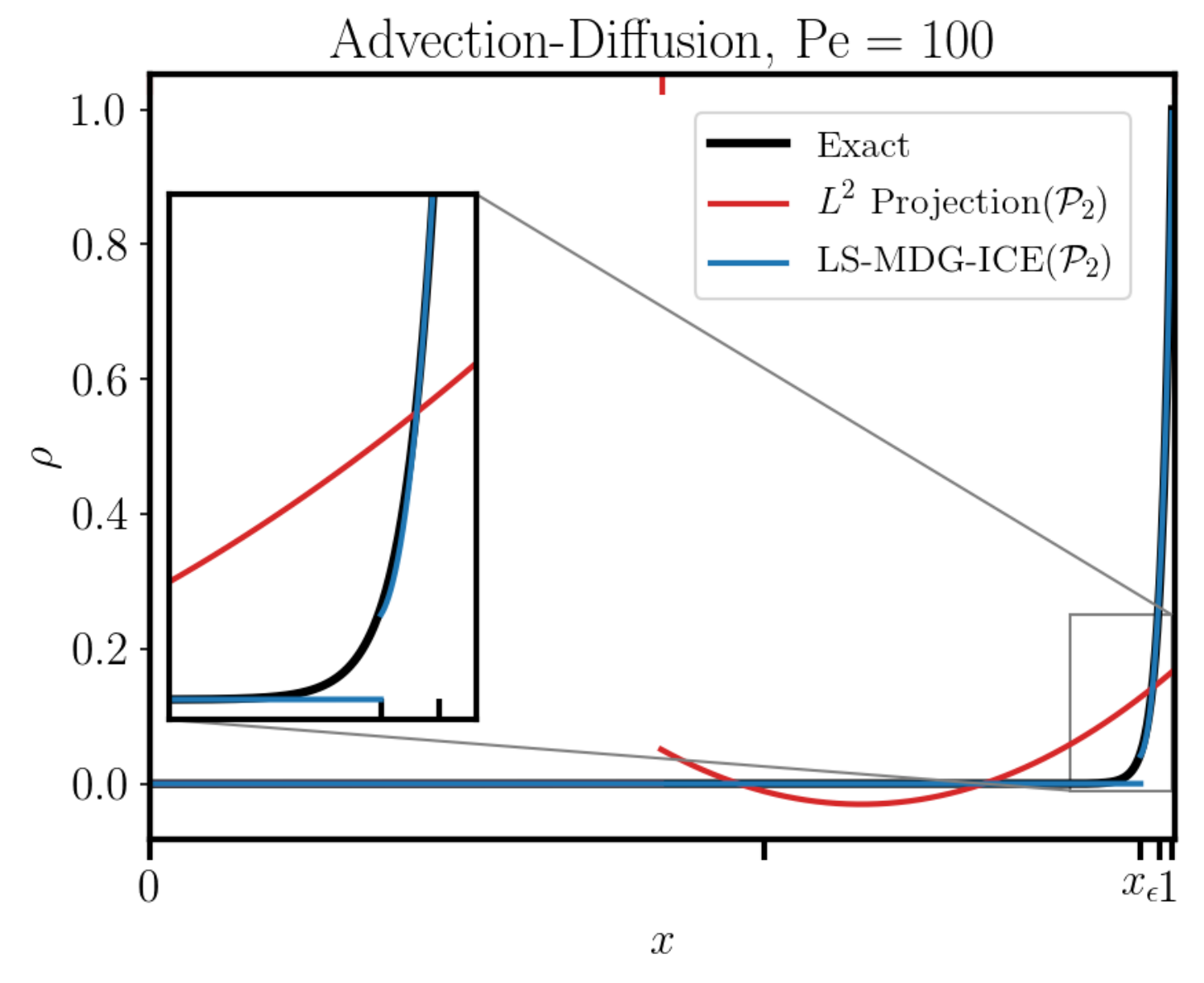} & \includegraphics[width=0.4\linewidth]{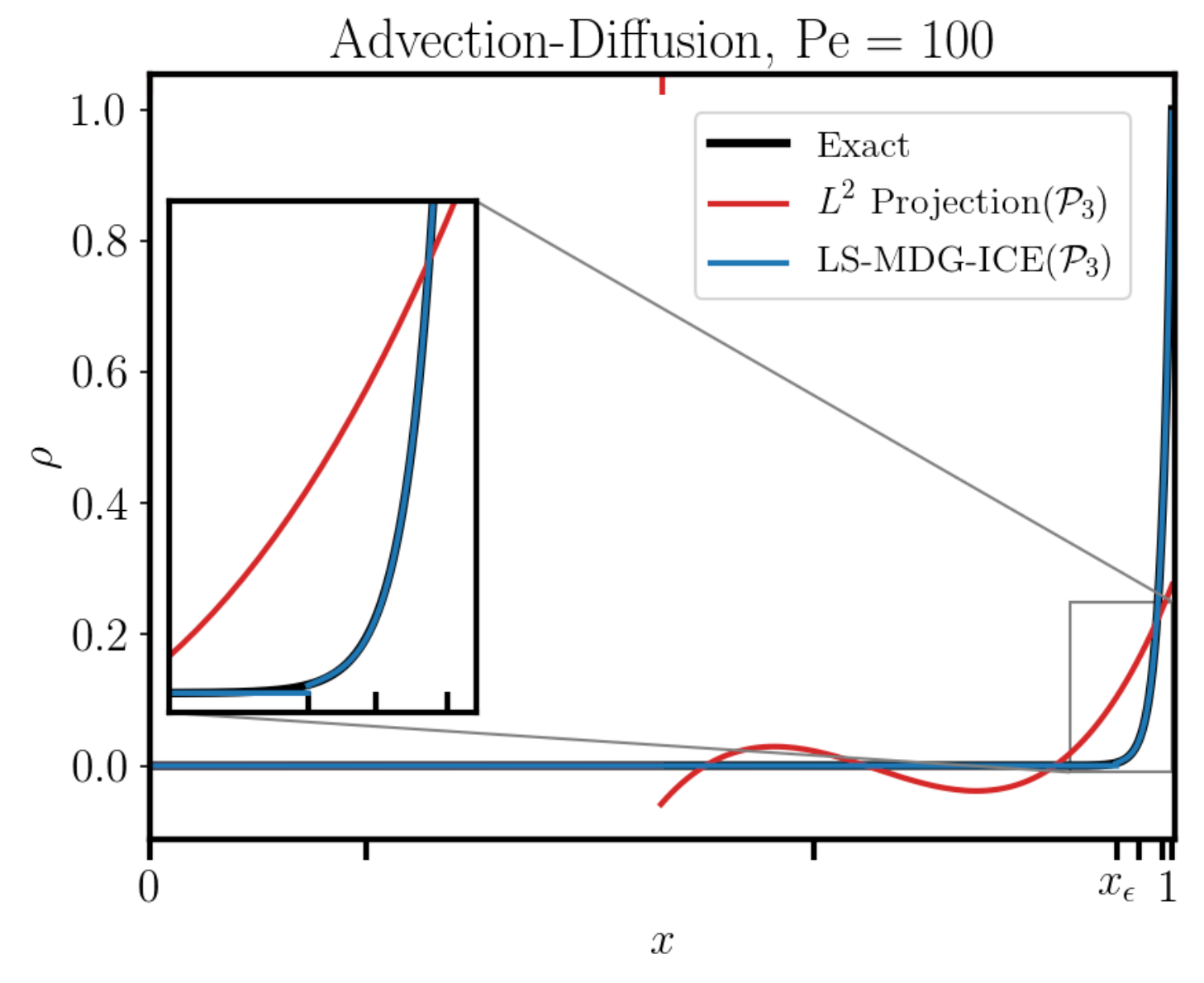}\tabularnewline
\includegraphics[width=0.4\linewidth]{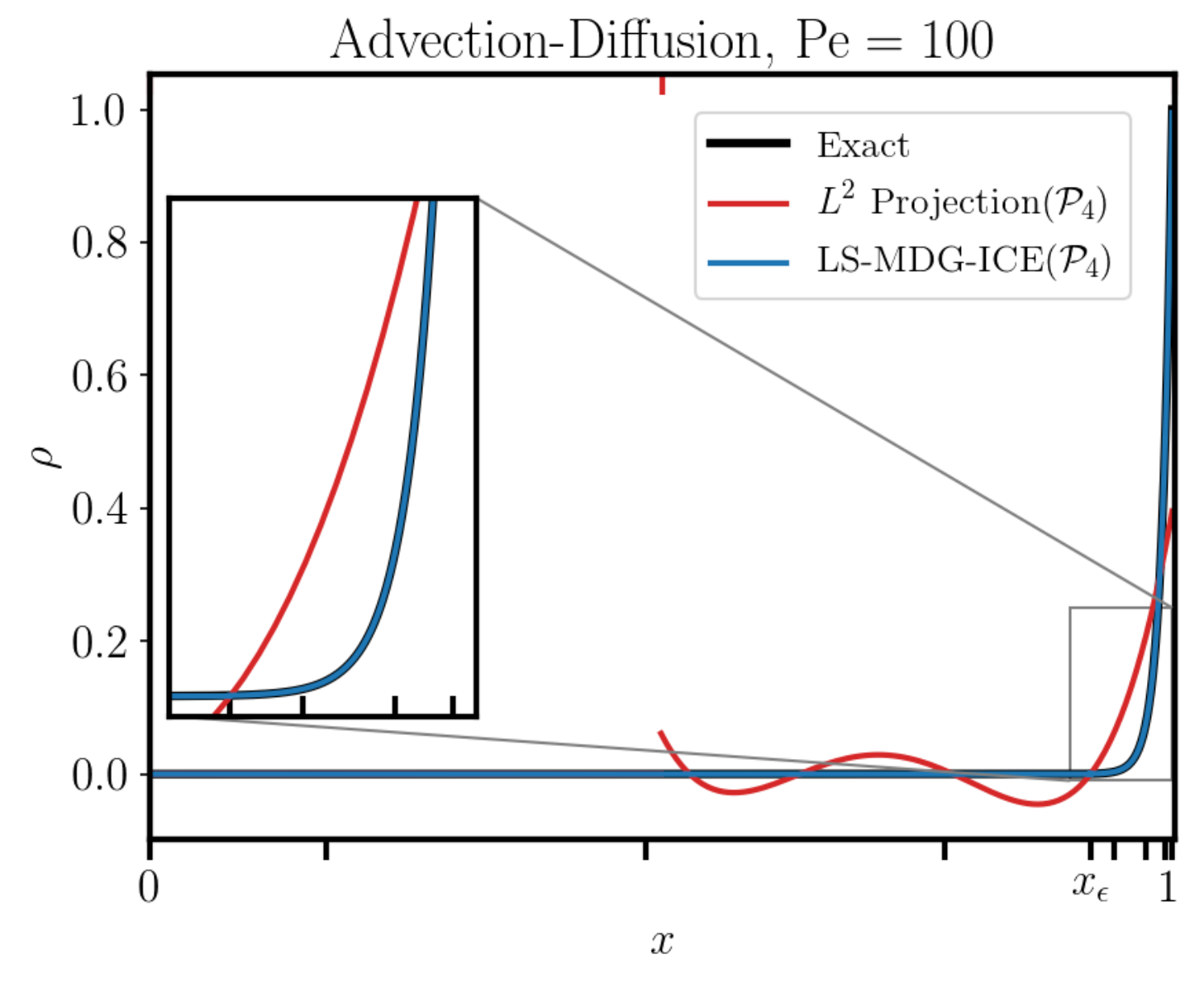} & \includegraphics[width=0.4\linewidth]{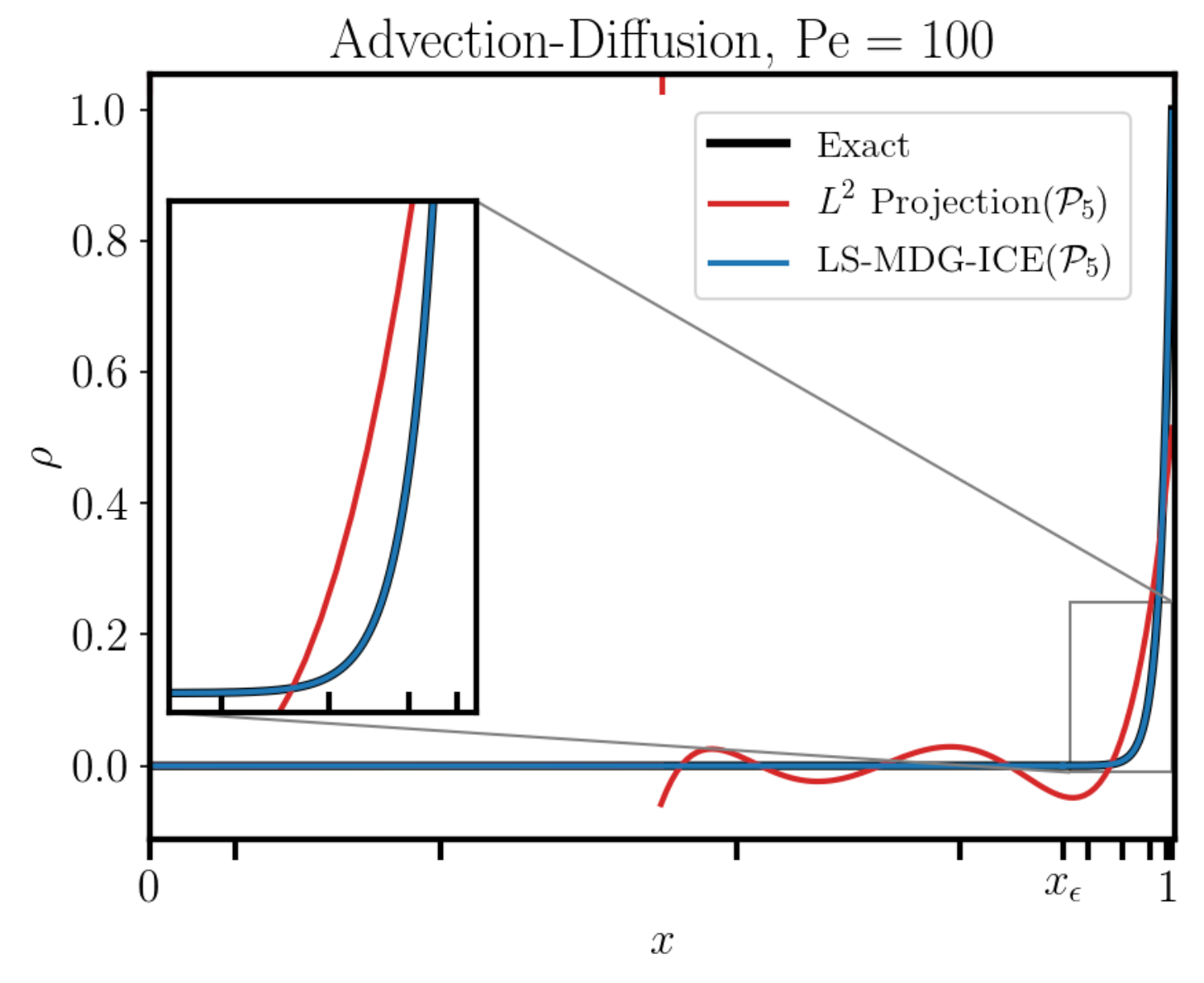}\tabularnewline
\end{tabular}\caption{\label{fig:steady-boundary-layer-Pe_0100}Grid adapted LS-MDG-ICE
solutions for a one-dimensional boundary layer type profile for $\mathrm{Pe}=100$
on a grid consisting of two isoparametric line cells are compared
to the exact solution and the $L^{2}$ projection of the exact solution
onto a uniform grid. The grid markers correspond to the adapted grid
and interior cell points given by the LS-MDG-ICE solutions. The interior
cell interface is denoted $x_{\epsilon}$. The initial uniform grid
is indicated with red tick marks.}
\end{figure}
\begin{figure}
\centering{}%
\begin{tabular}{cc}
\includegraphics[width=0.4\linewidth]{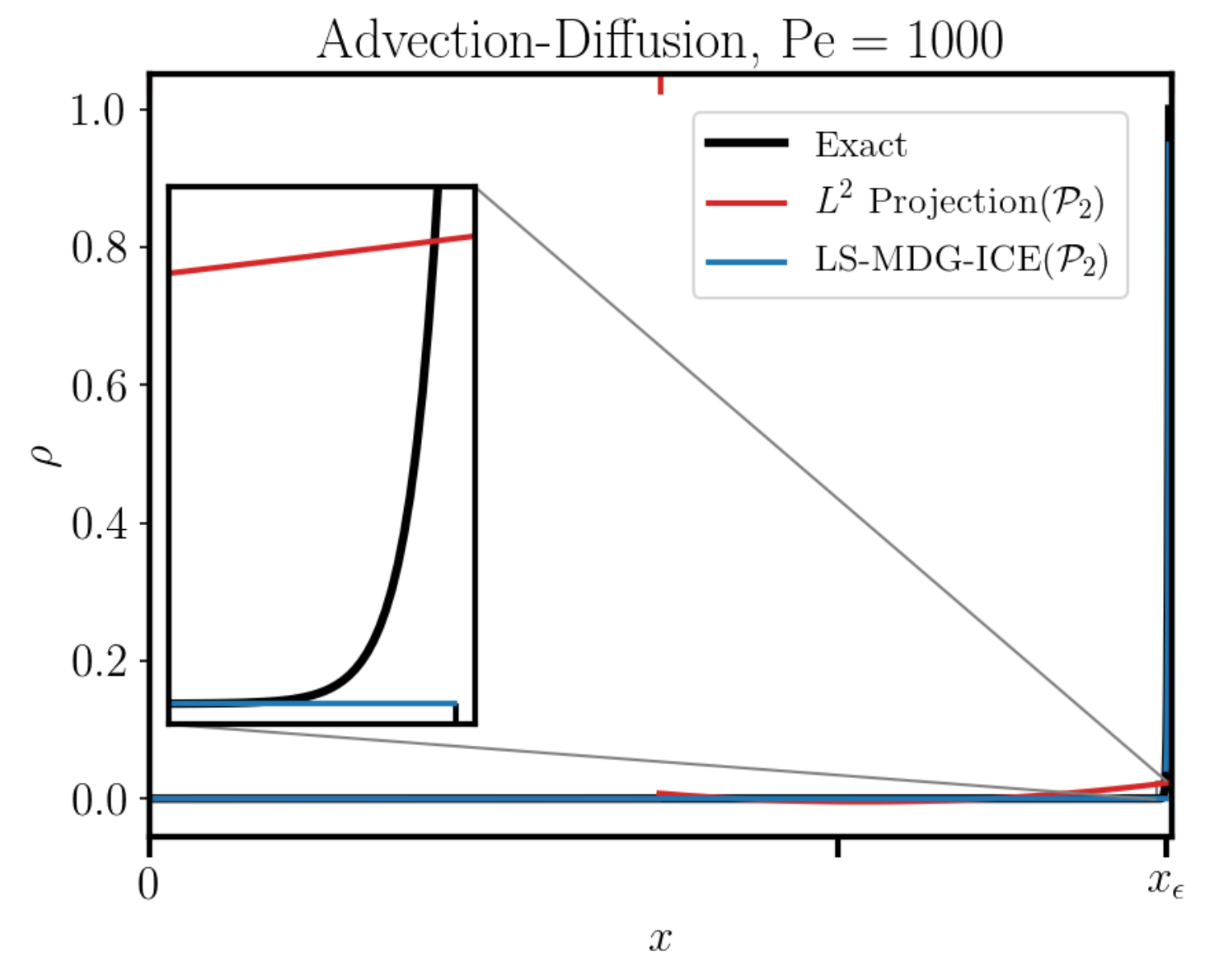} & \includegraphics[width=0.4\linewidth]{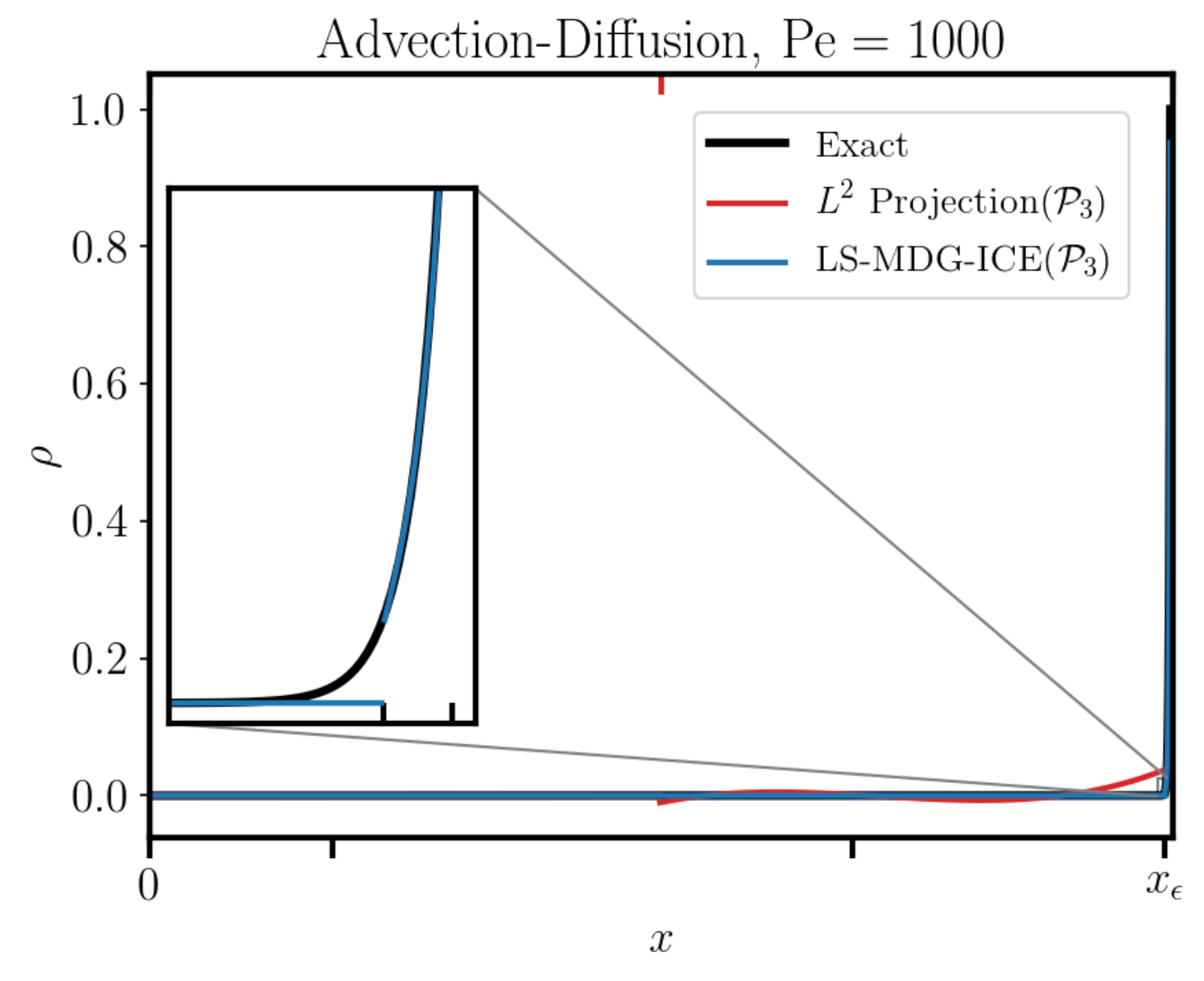}\tabularnewline
\includegraphics[width=0.4\linewidth]{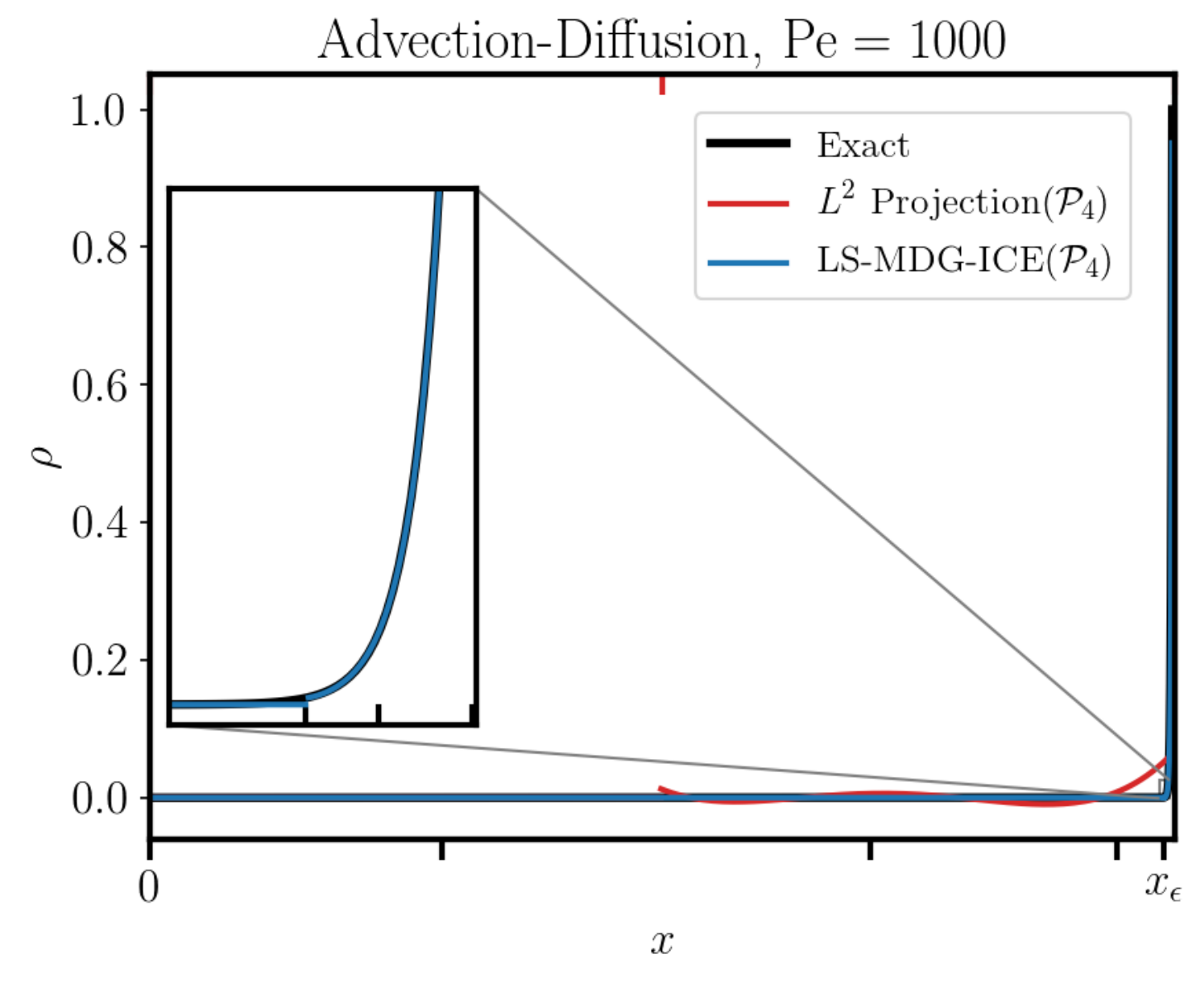} & \includegraphics[width=0.4\linewidth]{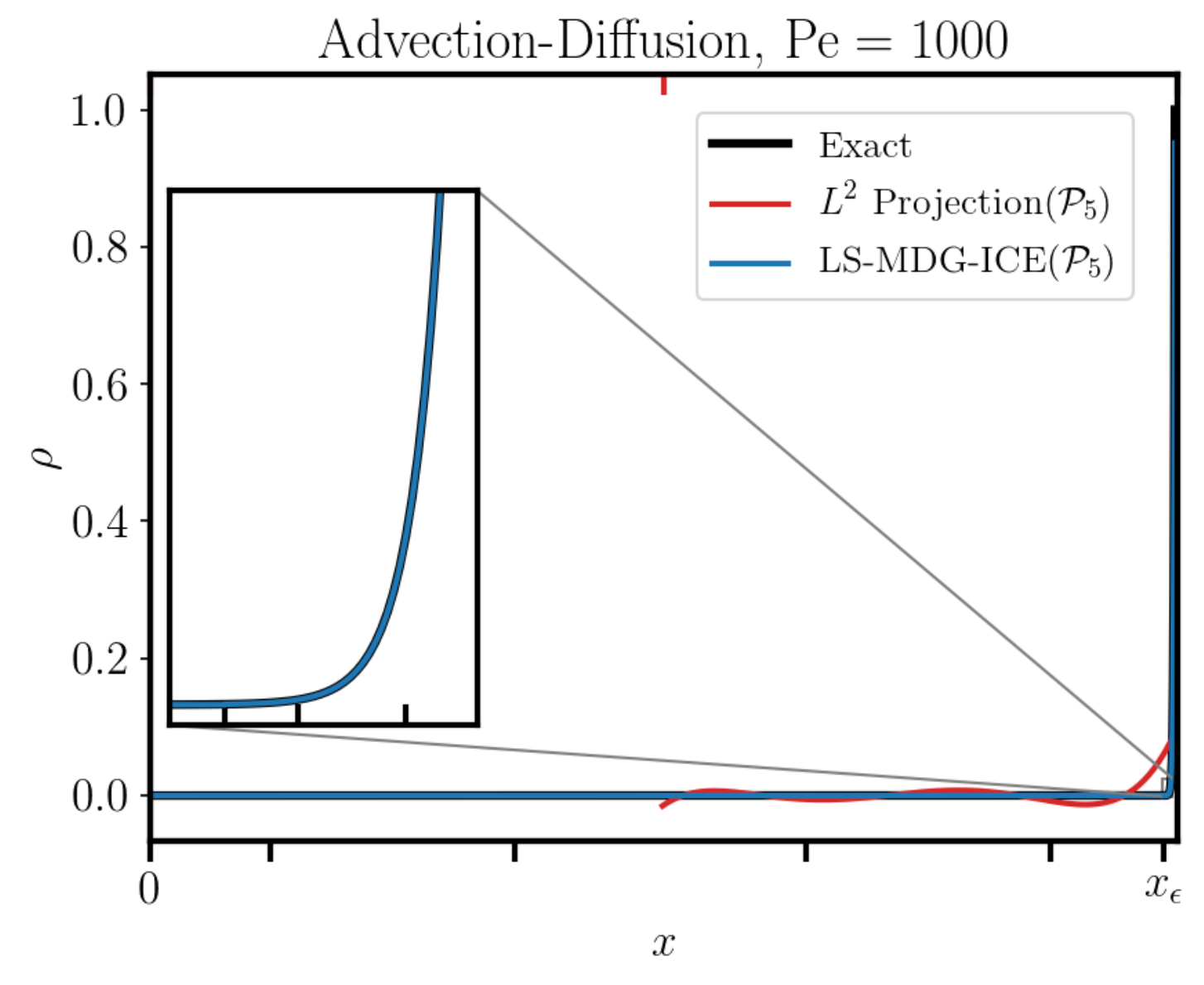}\tabularnewline
\end{tabular}\caption{\label{fig:steady-boundary-layer-Pe_1000}Grid adapted LS-MDG-ICE
solutions for a one-dimensional boundary layer type profile for $\mathrm{Pe}=1000$
on a grid consisting of two isoparametric line cells are compared
to the exact solution and the $L^{2}$ projection of the exact solution
onto a uniform grid. The grid markers correspond to the adapted grid
and interior cell points of the left cell given by the LS-MDG-ICE
solutions. For clarity, the interior points corresponding to the right
cell are not included. The interior cell interface is denoted $x_{\epsilon}$.
The initial uniform grid is indicated with red tick marks.}
\end{figure}

Figure~\ref{fig:steady-boundary-layer-p2} compares the LS-MDG-ICE
solutions for $\mathrm{Pe}=10,100,1000$ to the exact solution, which
is numerically unstable due to round-off error for $\mathrm{Pe}=1000$,
and overlays the three solutions LS-MDG-ICE solutions on rescaled
domains. As expected, the solutions are self-consistent. Furthermore,
the LS-MDG-ICE solutions do not exhibit oscillations in the presence
of sharp gradients.

Figure~\ref{fig:steady-boundary-layer-Pe_0100} and Figure~\ref{fig:steady-boundary-layer-Pe_1000}
present LS-MDG-ICE solutions on grids consisting of two isoparametric
line cells for polynomial degrees $\mathcal{P}_{2}$, $\mathcal{P}_{3}$,
$\mathcal{P}_{4}$, and $\mathcal{P}_{5}$ and Péclet numbers of 100
and 1000 respectively. The solutions do not exhibit oscillations,
instead, under-resolved solutions are represented discontinuously,
enabling highly accurate approximations of the gradient, or flux,
near the boundary. Furthermore, this case demonstrates LS-MDG-ICE
innate ability to not only reposition the interior cell interface,
denoted $x_{\epsilon}$ in Figures~\ref{fig:steady-boundary-layer-p2},~\ref{fig:steady-boundary-layer-p2},
but also nonlinearly deform the cell shape representation to improve
the accuracy of the approximation. We refer to this behavior as anisotropic
curvilinear $r$-adaptivity and, for the one-dimensional cases considered
in this work, we will provide evidence that it enables LS-MDG-ICE
to achieve super-optimal convergence with respect to the exact solution.

\begin{table}
\centering{}%
\begin{tabular}{c|cccc}
$\mathrm{Pe}$ & $x_{\epsilon}\left(\mathcal{P}_{2}\right)$ & $x_{\epsilon}\left(\mathcal{P}_{3}\right)$ & $x_{\epsilon}\left(\mathcal{P}_{4}\right)$ & $x_{\epsilon}\left(\mathcal{P}_{5}\right)$\tabularnewline
\hline 
$10$ & $0.74756464998474681$ & $0.68852737337261127$ & $0.64875294047110343$ & $0.53548391757202485$\tabularnewline
$100$ & $0.96910269349294942$ & $0.94529226568428737$ & $0.91970922330845084$ & $0.89308633425846451$\tabularnewline
$1000$ & $0.99690998474116876$ & $0.99452868699898989$ & $0.99196943402826943$ & $0.98930724837060802$\tabularnewline
$10,000$ & $0.99969099755499946$ & $0.99945286890577256$ & $0.99919687528580259$ & $0.99892997586082066$\tabularnewline
$100,000$ & $0.99996909975446002$ & $0.99994528681641881$ & $0.9999196875816686$ & $0.99989299759774075$\tabularnewline
\end{tabular}\caption{\label{tab:steady-boundary-layer-two-cell-interface-positions} Position
of the interior interface, $x_{\epsilon}$, corresponding to LS-MDG-ICE
solutions on a grid consisting of two isoparametric line cells for
$\mathrm{P_{e}}=10,10^{2},10^{3},10^{4},10^{5}$. The LS-MDG-ICE solutions
for $\mathrm{Pe}=10,100,1000$ are shown in Figures~\ref{fig:steady-boundary-layer-Pe_0100}-\ref{fig:steady-boundary-layer-Pe_1000}
respectively.}
\end{table}

Finally, in order to verify that LS-MDG-ICE automatically solves for
grids that resolve the relevant physical scales, we tabulate the positions
of the interior interface, $x_{\epsilon}$, for LS-MDG-ICE solutions
on a grid consisting of two isoparametric line cells in Table~\ref{tab:steady-boundary-layer-two-cell-interface-positions}.
The solutions for $\mathrm{Pe}=100,1000$ are shown in Figures~\ref{fig:steady-boundary-layer-Pe_0100}-\ref{fig:steady-boundary-layer-Pe_1000}
respectively, while the solutions for $\mathrm{Pe>1000}$ are not
shown. In order to converge at the optimal rate, the size of the cell
in the boundary layer is required to be

\begin{equation}
x_{r}-x_{\epsilon}\left(\mathcal{P}_{p},\mathrm{Pe}=10^{n}\right)\propto\frac{1}{\mathrm{\mathrm{Pe}}},\label{eq:interior-interface-position}
\end{equation}
for Péclet numbers, $\mathrm{Pe}=\left\{ 10^{n}:n\in\mathbb{N}\:\mathrm{and}\:n>1\right\} $,
where $x_{r}=1$ is the position of the right boundary, i.e., $x_{r}=1$.
For the data presented in Table~\ref{tab:steady-boundary-layer-two-cell-interface-positions}
indicates that LS-MDG-ICE automatically solves for grids on the diffusive
scale, i.e., $\mathcal{O}\left(\varepsilon=\nicefrac{1}{\mathrm{Pe}}\right)$,
overcoming a longstanding issue in optimal grid design for singularly
perturbed problems.

\subsubsection{Convergence for the case of linear advection-diffusion in one dimension}

We present convergence results for the $L^{2}$ error of the LS-MDG-ICE
approximation with respect to the exact solution~(\ref{eq:steady-boundary-layer-exact-solution})
for $\mathrm{Pe}=10$.
\begin{figure}
\subfloat[\label{fig:steady-boundary-layer-convergence-static-mesh}Grid convergence
of LS-MDG-ICE using linear line elements with a fixed discrete geometry. ]{\begin{centering}
\includegraphics[width=0.45\linewidth]{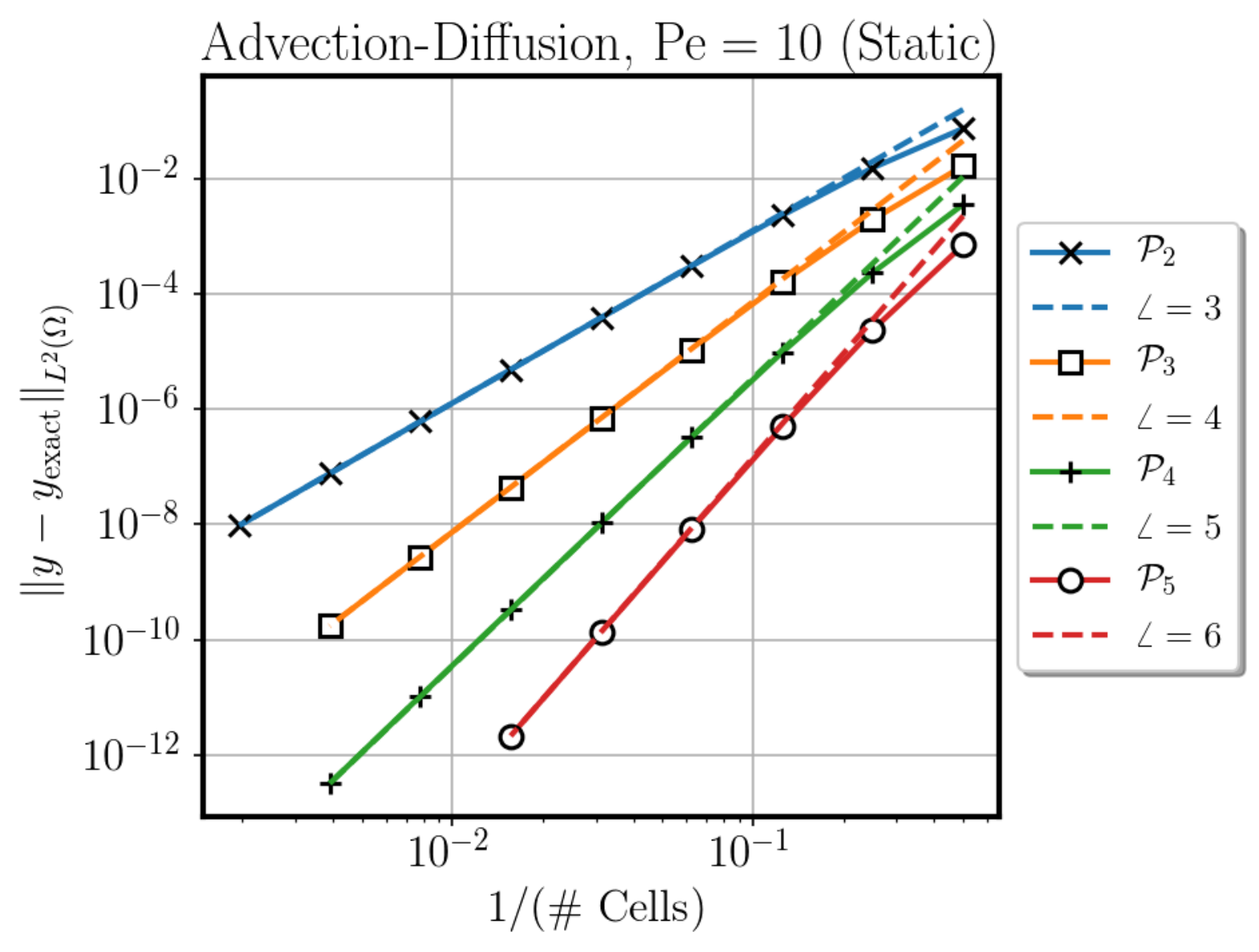}
\par\end{centering}
}\hfill{}
\begin{centering}
\subfloat[\label{fig:steady-boundary-layer-convergence-dynamic-mesh}Grid convergence
of LS-MDG-ICE using isoparametric line elements with a variable discrete
geometry.]{\begin{centering}
\includegraphics[width=0.45\linewidth]{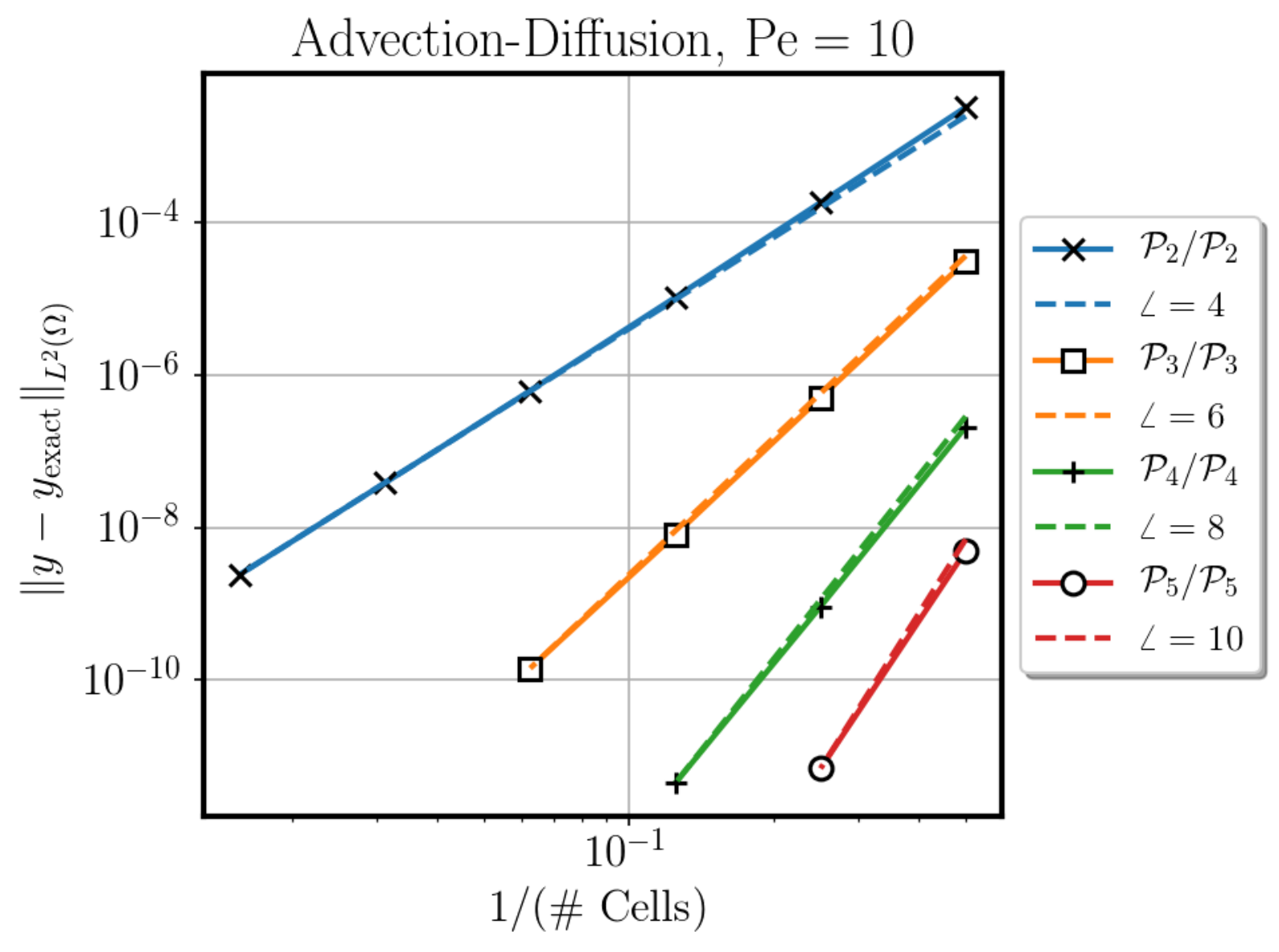}
\par\end{centering}
}
\par\end{centering}
\caption{\label{fig:steady-boundary-layer-convergence-static-dynamic-mesh}Convergence
plots for the linear advection-diffusion steady boundary layer type
solution with $\mathrm{Pe}=10$ for LS-MDG-ICE($\mathcal{P}_{2}$)
through LS-MDG-ICE($\mathcal{P}_{5}$) on both static and dynamic
grids. On static grids the LS-MDG-ICE discretization converges with
optimal ($p+1$) order. On dynamic grids, with isoparametric elements,
the LS-MDG-ICE discretization converges with super-optimal ($2p$)
order.}
\end{figure}
\begin{figure}
\subfloat[\label{fig:steady-boundary-layer-convergence-p2}LS-MDG-ICE($\mathcal{P}_{2}$)
static and variable discrete geometry convergence. ]{\begin{centering}
\includegraphics[width=0.45\linewidth]{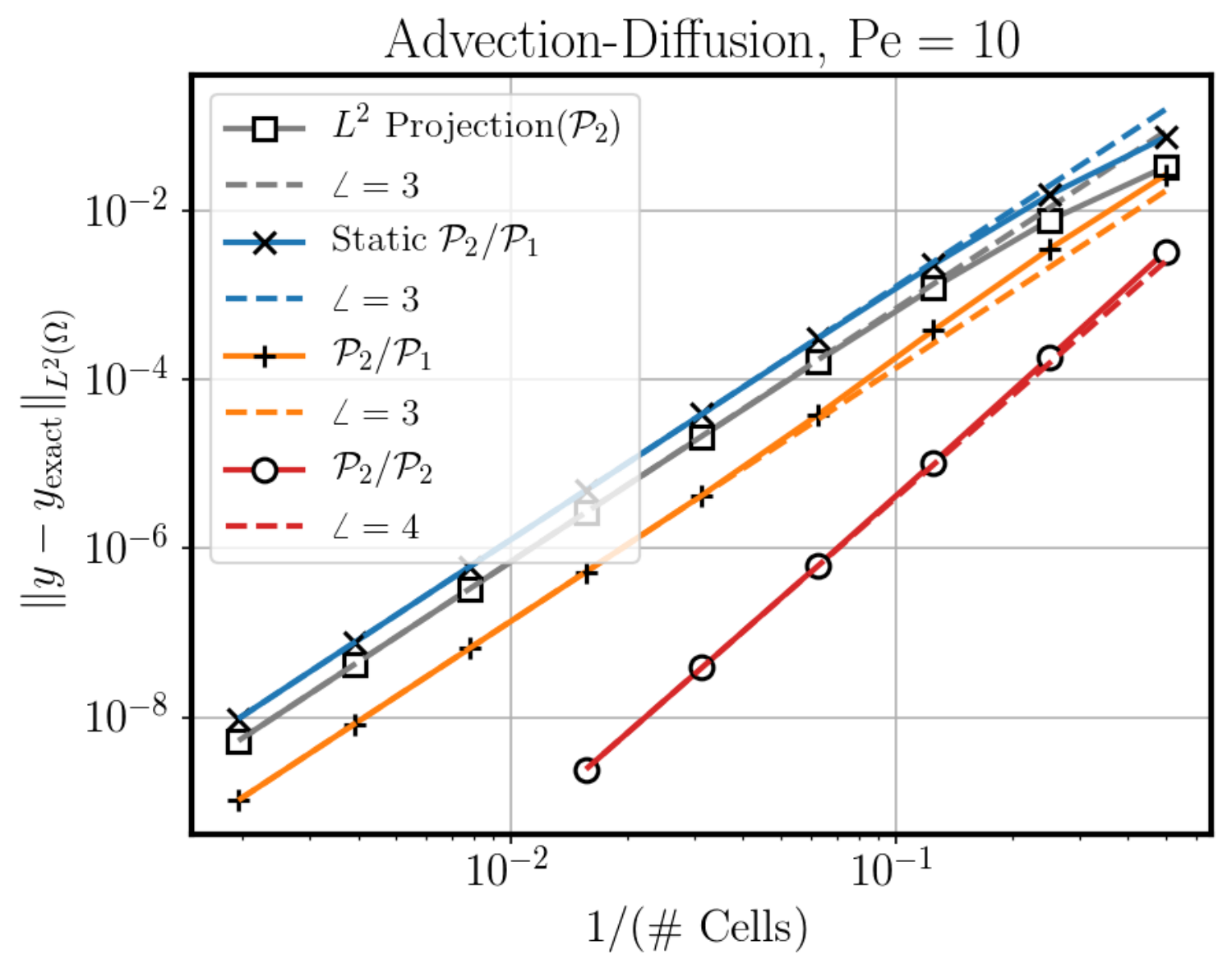}
\par\end{centering}
}\hfill{}
\begin{centering}
\subfloat[\label{fig:steady-boundary-layer-convergence-p3}LS-MDG-ICE($\mathcal{P}_{3}$)
variable discrete geometry convergence.]{\begin{centering}
\includegraphics[width=0.45\linewidth]{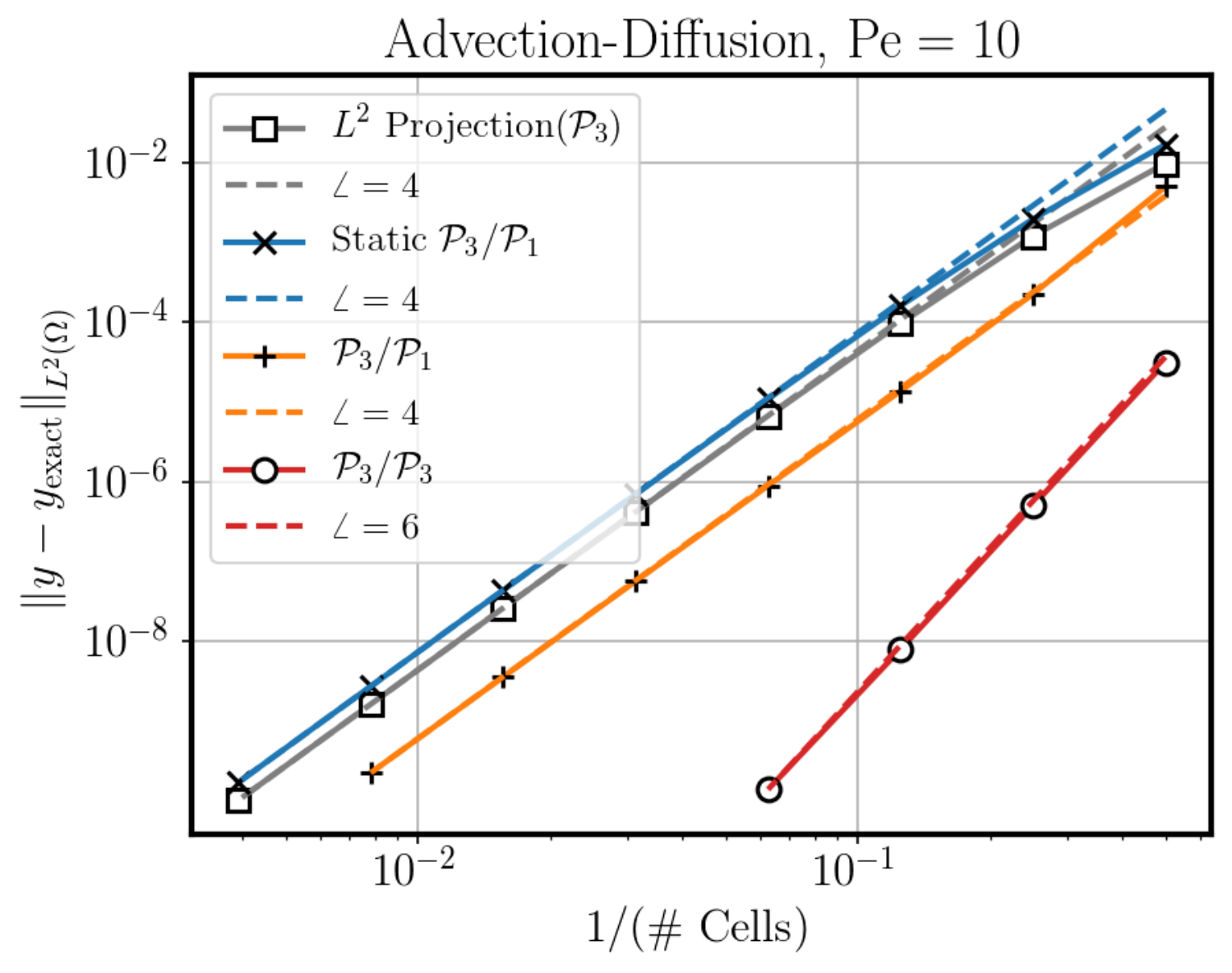}
\par\end{centering}
}
\par\end{centering}
\caption{\label{fig:steady-boundary-layer-convergence-p2-p3}Convergence plots
for the linear advection-diffusion steady boundary layer type solution
with $\mathrm{Pe}=10$ for LS-MDG-ICE($\mathcal{P}_{2}$) and LS-MDG-ICE($\mathcal{P}_{3}$).
Each of the coarsest grids consisted of $2$ line cells, while the
finest grid consisted of $512$, $512$, and $64$ lines cells for
the static LS-MDG-ICE($\mathcal{P}_{2}/\mathcal{P}_{1}$), LS-MDG-ICE($\mathcal{P}_{2}/\mathcal{P}_{1}$),
and LS-MDG-ICE($\mathcal{P}_{2}/\mathcal{P}_{2}$) discretizations
respectively and $256$, $128$, and $16$ lines cells for the static
LS-MDG-ICE($\mathcal{P}_{3}/\mathcal{P}_{1}$), dynamic LS-MDG-ICE($\mathcal{P}_{3}/\mathcal{P}_{1}$),
and dynamic LS-MDG-ICE($\mathcal{P}_{3}/\mathcal{P}_{3}$) discretizations
respectively.}
\end{figure}
\begin{figure}
\centering{}\includegraphics[width=0.48\linewidth]{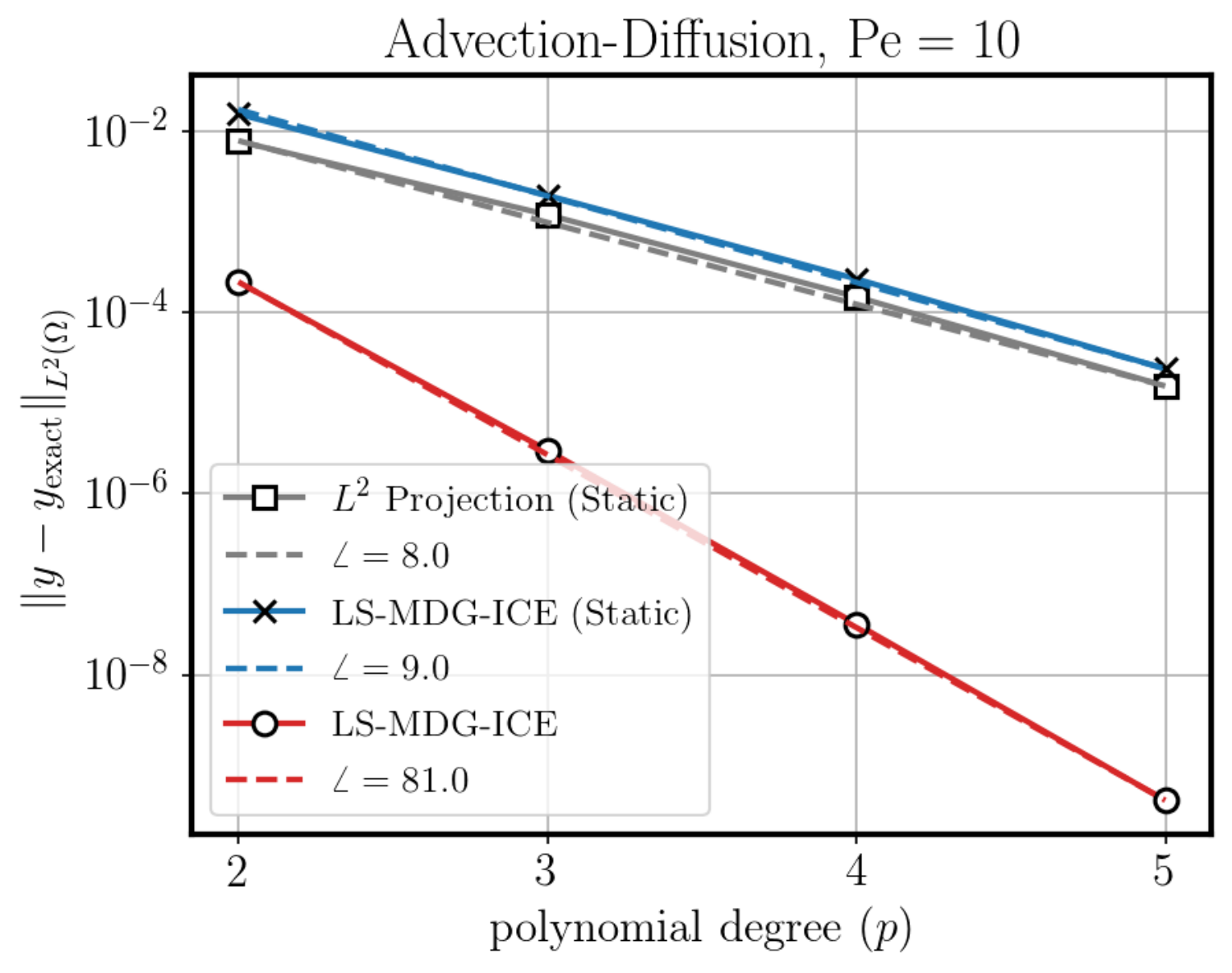}\caption{\label{fig:steady-boundary-layer-p-convergence-Pe-0010}Comparison
of convergence under polynomial refinement for the linear advection-diffusion
steady boundary layer type solution with $\mathrm{Pe}=10$ for LS-MDG-ICE
with a static and variable discrete geometry on a grid consisting
of four isoparametric line elements.}
\end{figure}
Figure~\ref{fig:steady-boundary-layer-convergence-static-dynamic-mesh}
presents convergence rates of LS-MDG-ICE with both a fixed and variable
discrete geometry for polynomial degrees $p=2,\ldots,5$. The coarsest
grid consists of $2$ line cells, while the finest grid consists of
$256$ line cells. For a fixed discrete geometry, LS-MDG-ICE exhibits
optimal $\left(p+1\right)$ order convergence. For a variable discrete
geometry, LS-MDG-ICE exhibits super-optimal $\left(2p\right)$ order
convergence.

Figure~\ref{fig:steady-boundary-layer-convergence-p2-p3} directly
compares the convergence rates of LS-MDG-ICE with a fixed and variable
discrete geometry for $p=2,3$. For reference, we have also included
the error associated with the $L^{2}$ projection of the exact solution
onto a uniform grid, which provides an upper bound on the accuracy
attainable by methods based on a static grid since it minimizes the
error in the $L^{2}$ norm. For the case of a variable discrete geometry,
we consider both linear and isoparametric line elements in order.
In the case of linear line elements LS-MDG-ICE exhibits the same optimal
$\left(p+1\right)$ order convergence as the case of a fixed discrete
geometry. However, the solution is almost 8 times more accurate. The
solution is most accurate for the case of isoparametric line cells.
The LS-MDG-ICE($\mathcal{P}_{3}$) requires $256$, $128$, and $16$
line cells to achieve an error on the order of $10^{-9}$ for the
case of a static discrete geometry, variable linear discrete geometry,
and variable isoparametric discrete geometry. 

Figure~\ref{fig:burgers-steady-viscous-shock-convergence-p2} compares
convergence of LS-MDG-ICE solutions under polynomial refinement for
$\mathrm{Pe}=10$ on grids consisting of four line cells for the case
of static discrete geometry and variable isoparametric discrete geometry.
Again, we include the error associated with the $L^{2}$ projection
of the exact solution onto a uniform grid for reference. The $L^{2}$
projection, static LS-MDG-ICE, and LS-MDG-ICE errors are plotted on
a log-linear plot, with a reference slope of $8$, $9$, and $81$,
respectively.

\subsection{Burgers steady viscous shock\label{sec:burgers-steady-viscous-shock}}

We consider the solution of the nonlinear conservation law in one
dimension corresponding to the viscous Burgers equation~\citep{Hop50,Col51}.
The exact solution is given by

\begin{equation}
y\left(x\right)=y_{R}+\frac{1}{2}\left(y_{L}-y_{R}\right)\left(1-\tanh\left(\frac{\left(y_{L}-y_{R}\right)\left(x-v_{s}t\right)}{4\epsilon}\right)\right)\label{eq:viscous-burgers-exact}
\end{equation}
where

\begin{equation}
v_{s}=\frac{1}{2}\left(y_{L}+y_{R}\right)\label{eq:burgers-shock-speed}
\end{equation}
is the shock speed, $\epsilon$ is a viscosity coefficient, and $t$
is the time. The boundary conditions at $x=\pm\nicefrac{1}{2}$ are
given by Equation~(\ref{eq:viscous-burgers-exact}), which results
in a stationary viscous shock, i.e., $v_{s}=0$. The solutions is
initialized with a piecewise constant profile
\begin{equation}
y\left(x\right)=\begin{cases}
y\left(\nicefrac{-1}{2}\right) & \textup{ if }x_{c}\le0\\
y\left(\nicefrac{+1}{2}\right) & \textup{ if }x_{c}>0
\end{cases}\label{eq:viscous-burgers-initial-condition}
\end{equation}
where $x_{c}$ is the cell centroid. For certain parameter regimes,
$\epsilon\le10^{-3}$, the solution to the problem is not unique since
perturbations to the boundary data are smaller than machine precision~\citep{Cha10}.
\begin{figure}
\subfloat[\label{fig:burgers-steady-viscous-shock-p2-Re-0100}LS-MDG-ICE$\left(\mathcal{P}_{2}\right)$ ]{\begin{centering}
\includegraphics[width=0.45\linewidth]{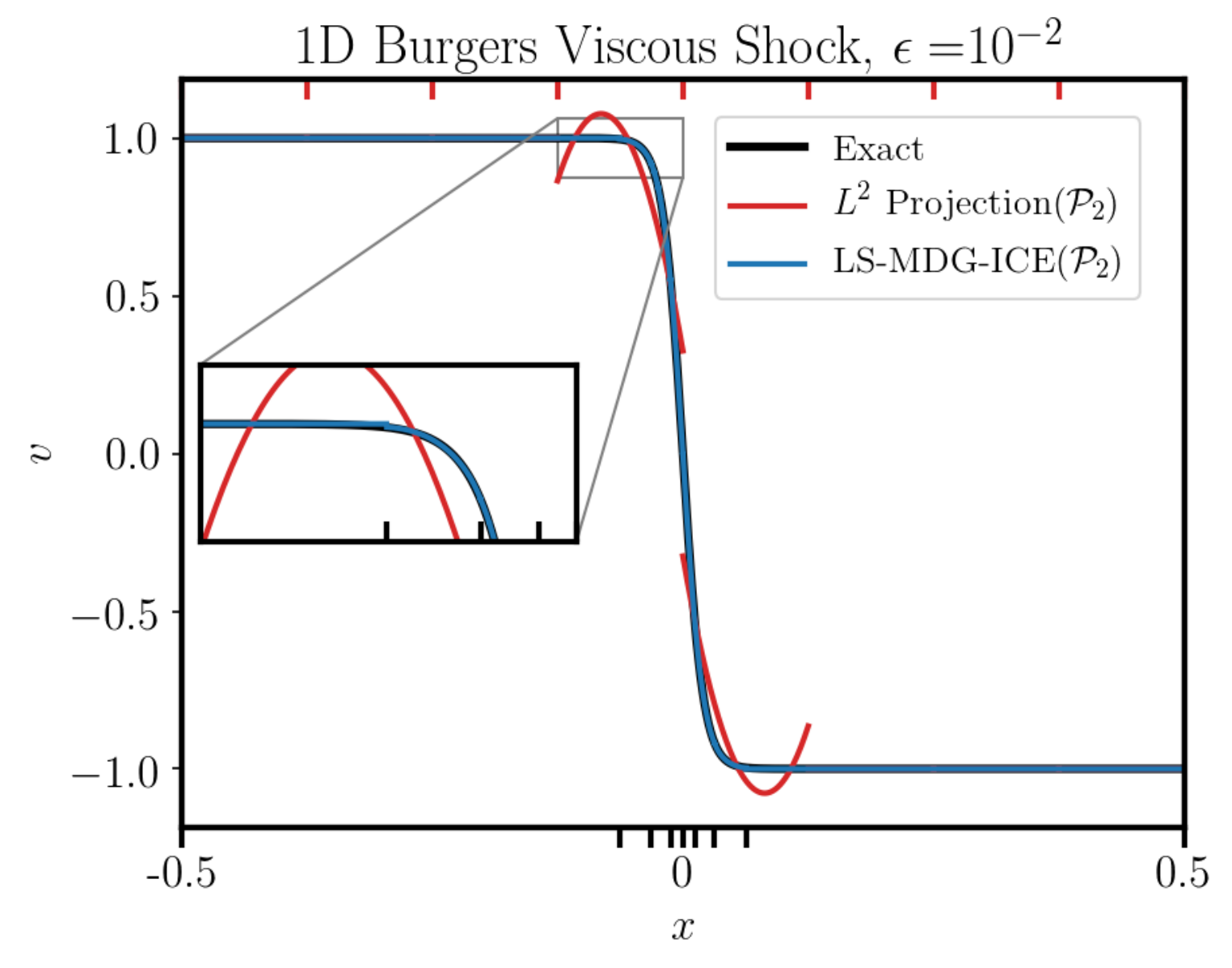}
\par\end{centering}
}\hfill{}
\begin{centering}
\subfloat[\label{fig:burgers-steady-viscous-shock-p4-Re-0100}LS-MDG-ICE$\left(\mathcal{P}_{4}\right)$ ]{\begin{centering}
\includegraphics[width=0.45\linewidth]{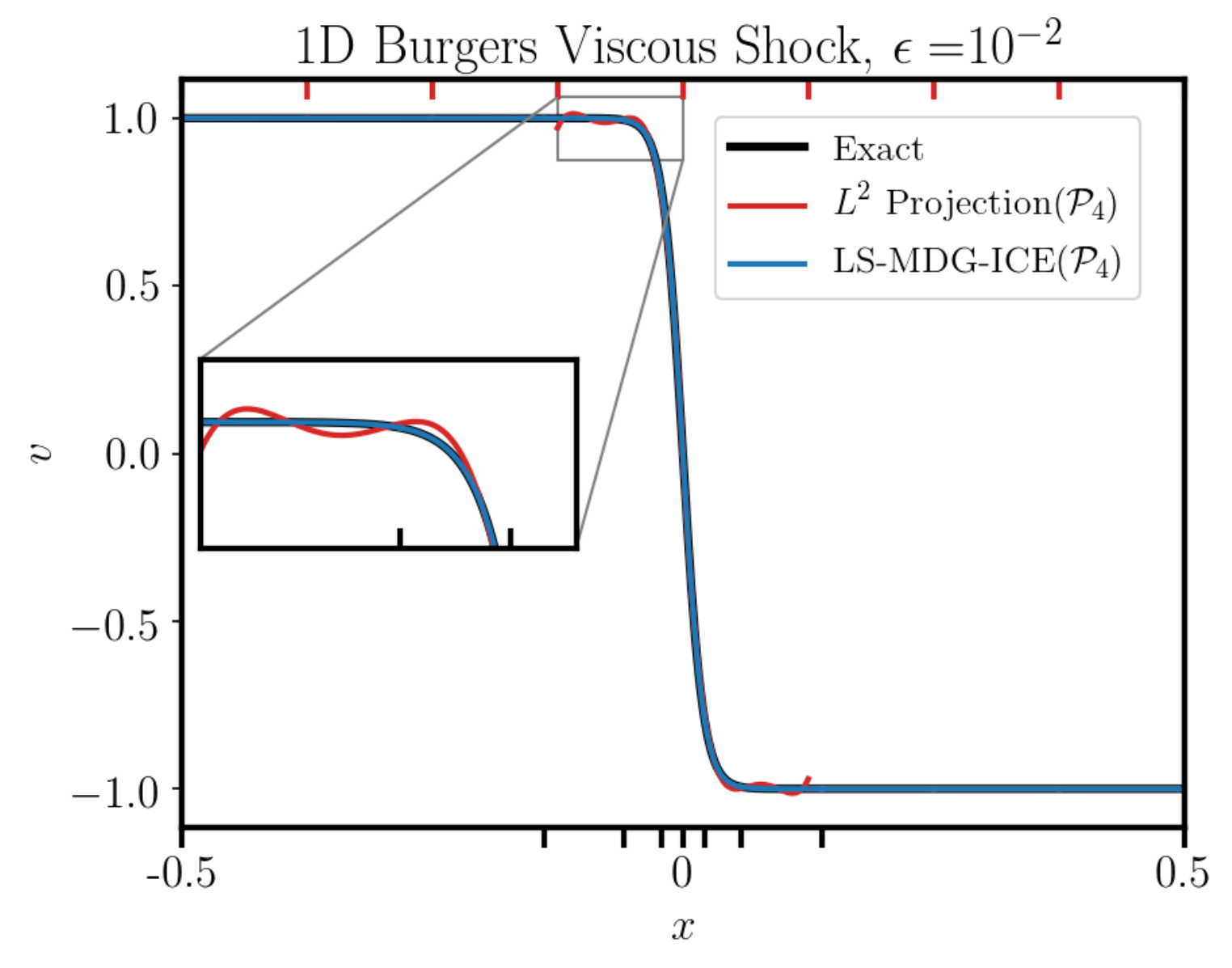}
\par\end{centering}
}
\par\end{centering}
\caption{\label{fig:burgers-steady-viscous-shock-Re-0100}Burgers steady viscous
shock solutions for $\epsilon=10^{-2}$ computed with LS-MDG-ICE$\left(\mathcal{P}_{2}\right)$
and LS-MDG-ICE$\left(\mathcal{P}_{4}\right)$ on eight isoparametric
line cells. The solutions were initialized on a uniform grid with
a piecewise constant profile given by Equation~(\ref{eq:viscous-burgers-initial-condition}).
The initial grid points are marked with red ticks on the top of the
plots and the adapted grid points given by the LS-MDG-ICE solutions
are marked with black ticks on the bottom of the plots.}
\end{figure}
\begin{figure}
\subfloat[\label{fig:burgers-steady-viscous-shock-p2-Re-1000}LS-MDG-ICE$\left(\mathcal{P}_{2}\right)$ ]{\begin{centering}
\includegraphics[width=0.45\linewidth]{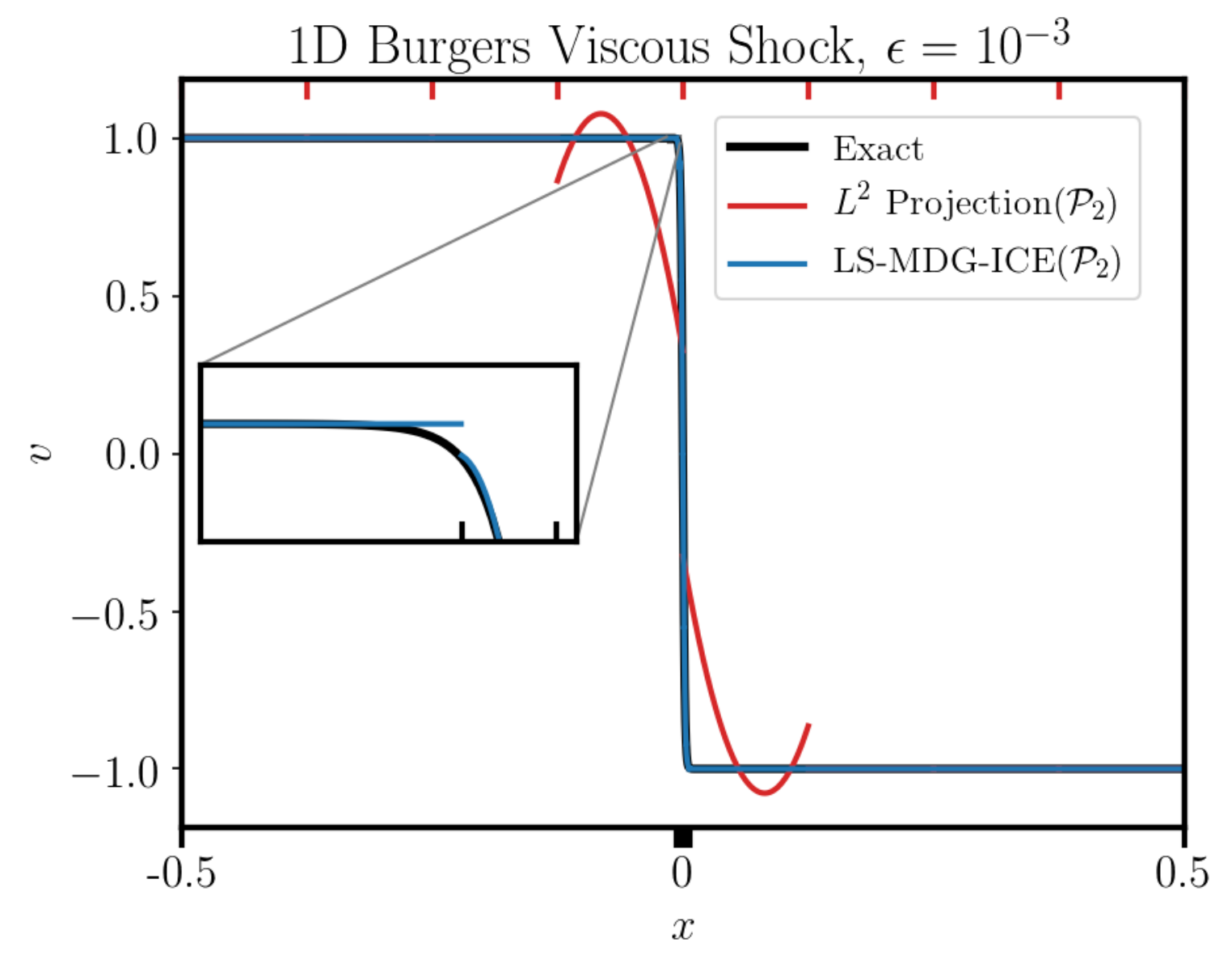}
\par\end{centering}
}\hfill{}
\begin{centering}
\subfloat[\label{fig:burgers-steady-viscous-shock-p4-Re-1000}LS-MDG-ICE$\left(\mathcal{P}_{4}\right)$ ]{\begin{centering}
\includegraphics[width=0.45\linewidth]{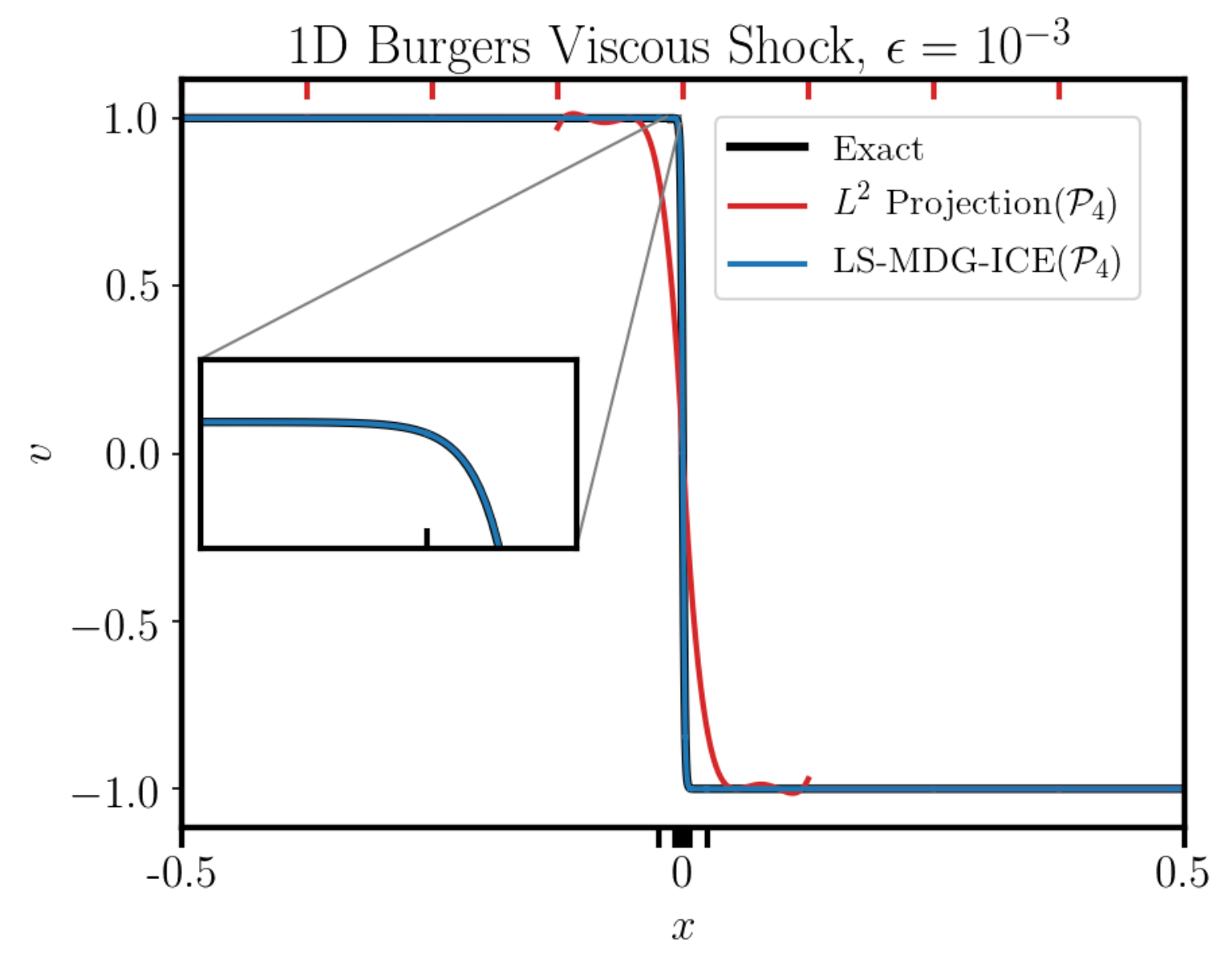}
\par\end{centering}
}
\par\end{centering}
\caption{\label{fig:burgers-steady-viscous-shock-Re-1000}Burgers steady viscous
solutions shock for $\epsilon=10^{-3}$ computed with LS-MDG-ICE$\left(\mathcal{P}_{2}\right)$
and LS-MDG-ICE$\left(\mathcal{P}_{4}\right)$ on eight isoparametric
line cells. The solutions were initialized on a uniform grid with
a piecewise constant profile given by Equation~(\ref{eq:viscous-burgers-initial-condition}).
The initial grid points are marked with red ticks on the top of the
plots and the adapted grid points given by the LS-MDG-ICE solutions
are marked with black ticks on the bottom of the plots.}
\end{figure}

Figure~\ref{fig:burgers-steady-viscous-shock-Re-0100} and Figure~\ref{fig:burgers-steady-viscous-shock-Re-1000}
present Burgers steady viscous shock solutions for $\epsilon=10^{-2}$
and $\epsilon=10^{-3}$ computed with LS-MDG-ICE$\left(\mathcal{P}_{2}\right)$
and LS-MDG-ICE$\left(\mathcal{P}_{4}\right)$ on eight isoparametric
line cells. The $L^{2}$ projection of the exact solution~(\ref{eq:viscous-burgers-exact})
onto a uniform grid is plotted as well and the initial uniform grid
is marked with red ticks. The LS-MDG-ICE solutions do not contain
oscillations, a problematic feature of under-resolved polynomial approximations,
which can be seen in the $L^{2}$ projection of the exact solution.
The discontinuities present in the LS-MDG-ICE$\left(\mathcal{P}_{2}\right)$
solutions shown in Figure~\ref{fig:burgers-steady-viscous-shock-p4-Re-0100}
and Figure~\ref{fig:burgers-steady-viscous-shock-p4-Re-1000} provide
a mechanism for accurately representing under-resolved features, which
localizes the error, and prevents it from polluting the solution upstream.
Furthermore, as the approximation order is increased, the unphysical
discontinuities are resolved and the LS-MDG-ICE solution becomes indistinguishable
from the exact solution. Figure~\ref{fig:burgers-steady-viscous-shock-p4-Re-0100}
and Figure~\ref{fig:burgers-steady-viscous-shock-p4-Re-1000} present
the LS-MDG-ICE$\left(\mathcal{P}_{4}\right)$ solution computed on
eight isoparametric line cells for $\epsilon=10^{-2}$ and $\epsilon=10^{-3}$
respectively.

\subsubsection{Convergence for the case of viscous Burgers flow in one dimension}

\begin{figure}
\centering{}\includegraphics[width=0.48\linewidth]{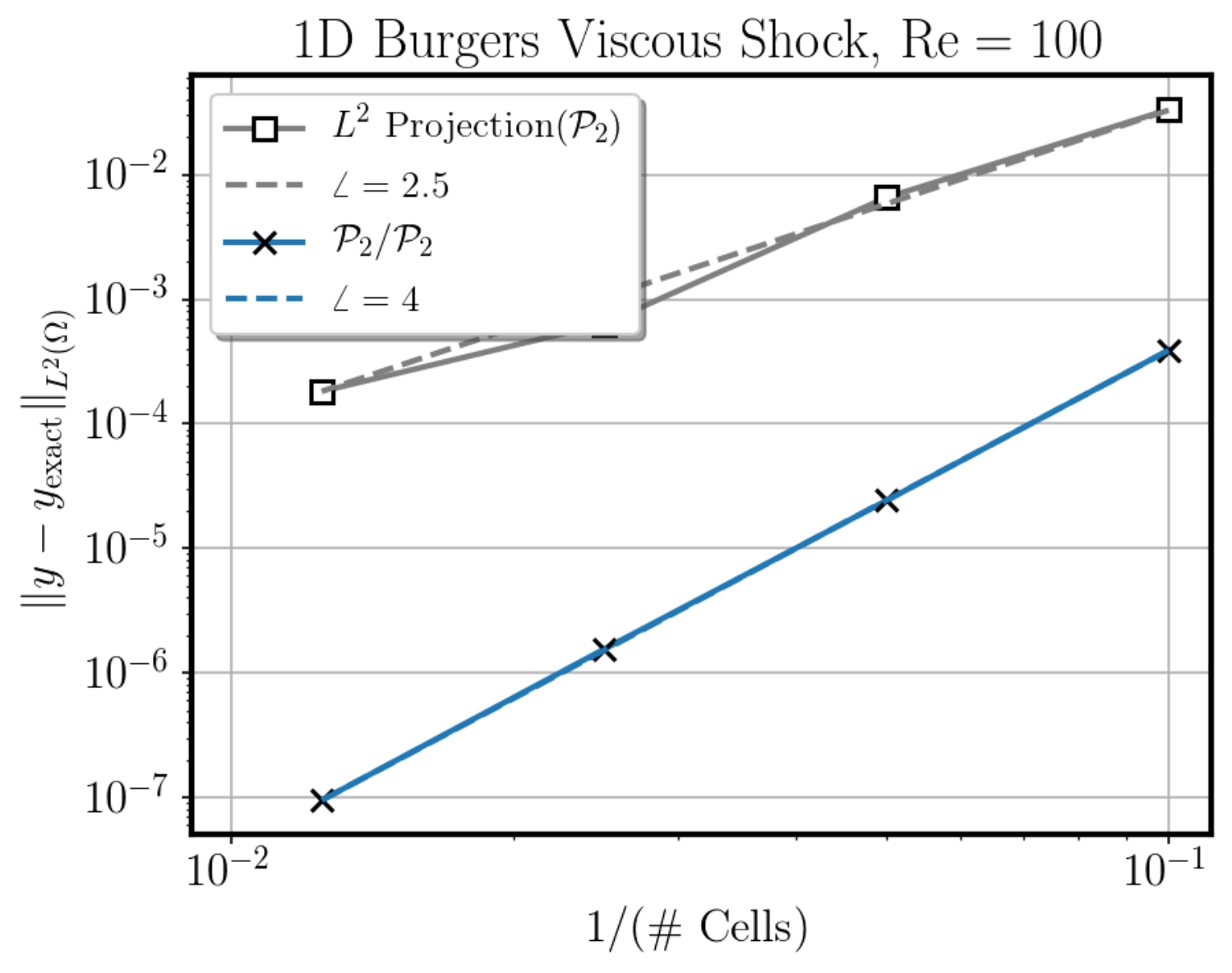}\caption{\label{fig:burgers-steady-viscous-shock-convergence-p2}Convergence
plot for Burgers steady viscous shock with $\epsilon=10^{-2}$ using
LS-MDG-ICE($\mathcal{P}_{2}$) isoparametric elements. The coarsest
grid consists of $10$ isoparametric line elements and the finest
grid consists of $80$ isoparametric line elements. }
\end{figure}
Figure~\ref{fig:burgers-steady-viscous-shock-convergence-p2} presents
convergence results under grid refinement for Burgers steady viscous
shock at $\epsilon=10^{-2}$ using LS-MDG-ICE($\mathcal{P}_{2}$).
The coarsest grid consists of $10$ isoparametric line elements and
the finest grid consists of $80$ isoparametric line elements. At
these moderate conditions, an excessively fine grid resolution is
required to obtain the optimal convergence rate under uniform refinement.
This lack of resolution is manifested in the error associated with
the $L^{2}$ projection of the exact solution, which converges sub-optimally
with a rate of $2.5$. In contrast, LS-MDG-ICE adapts the grid to
resolve the relevant physical scales and therefore overcomes the challenges
associated with traditional methods based on static grids. In fact,
LS-MDG-ICE again achieves super-optimal convergence at the rate of
$2p$ in the case of variable isoparametric discrete geometry. For
static grid methods, the use of a layer adapted grid, where the grid
is non-uniformly refined, or graded, to produce higher resolution
in regions with steep gradients, could, at best, restore optimal convergence
rates.

\subsection{Compressible Navier-Stokes steady viscous shock\label{sec:CNS-Viscous-Shock}}

We present the solution to a one-dimensional viscous shock for a compressible
Navier-Stokes flow in order to verify the nonlinear viscous flux formulation
and implementation. The initial left and right states, as well as
the boundary conditions, i.e., $y\left(x=-1\right)$ and $y\left(x=-1\right)$,
are given as
\begin{equation}
\left(\rho,v,T\right)\left(x\right)=\begin{cases}
\left(1,3.5,1\right) & x<0\\
\left(4.260869565217392,0.8214285714285714,3.3150510204081627\right) & x>0
\end{cases},\label{eq:CNS-Viscous-Shock-Initial-Conditions}
\end{equation}
where the right state was derived such that the jump conditions given
by Equation~(\ref{eq:interface-condition-strong-viscous}) are satisfied.
In one dimension, the problem reduces to a system of ordinary differential
equations that can be solved numerically, see~\citep{Cha10,Mas13}
for details.

\begin{figure}
\centering{}%
\begin{tabular}{ccc}
\includegraphics[width=0.3\linewidth]{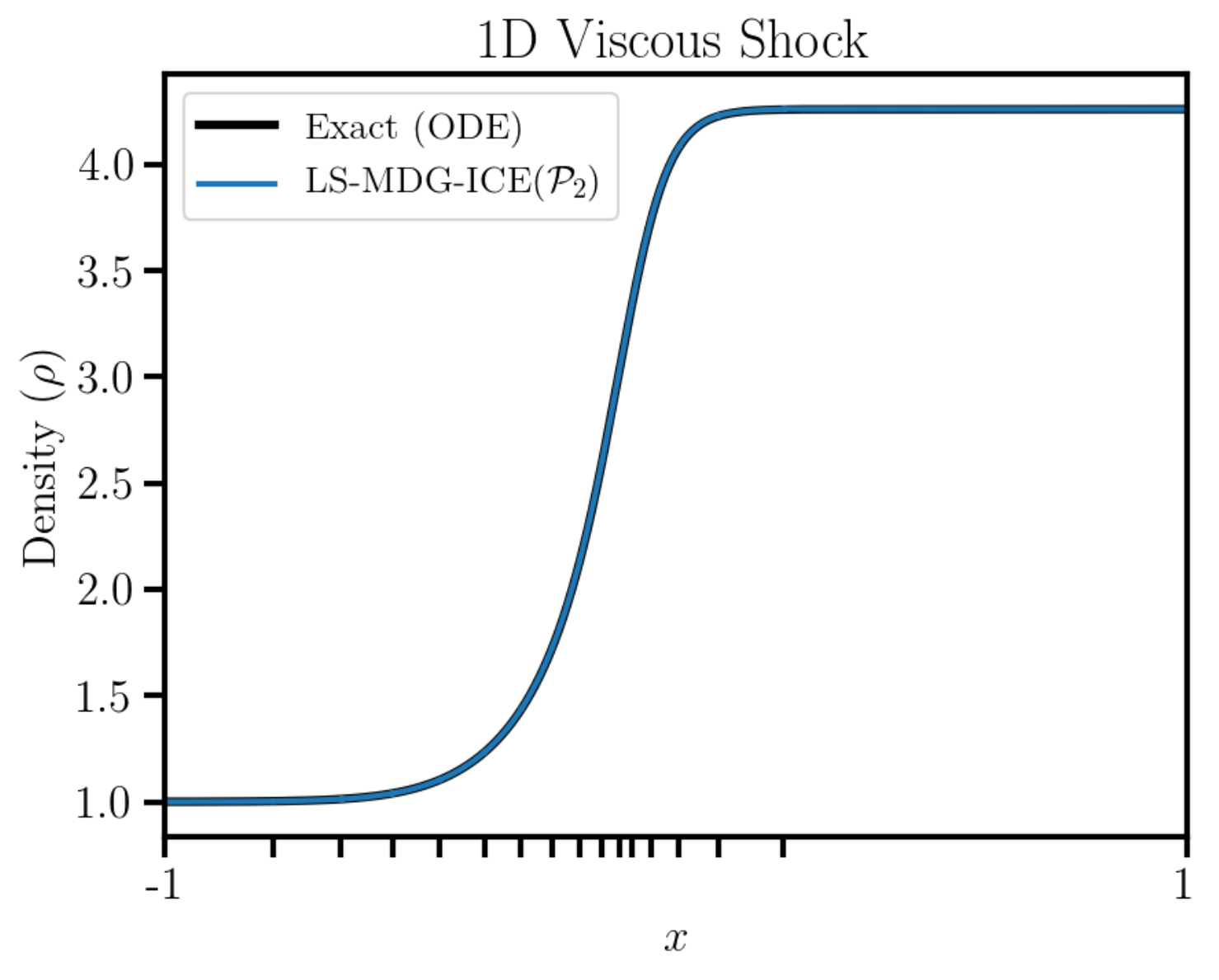} & \includegraphics[width=0.3\linewidth]{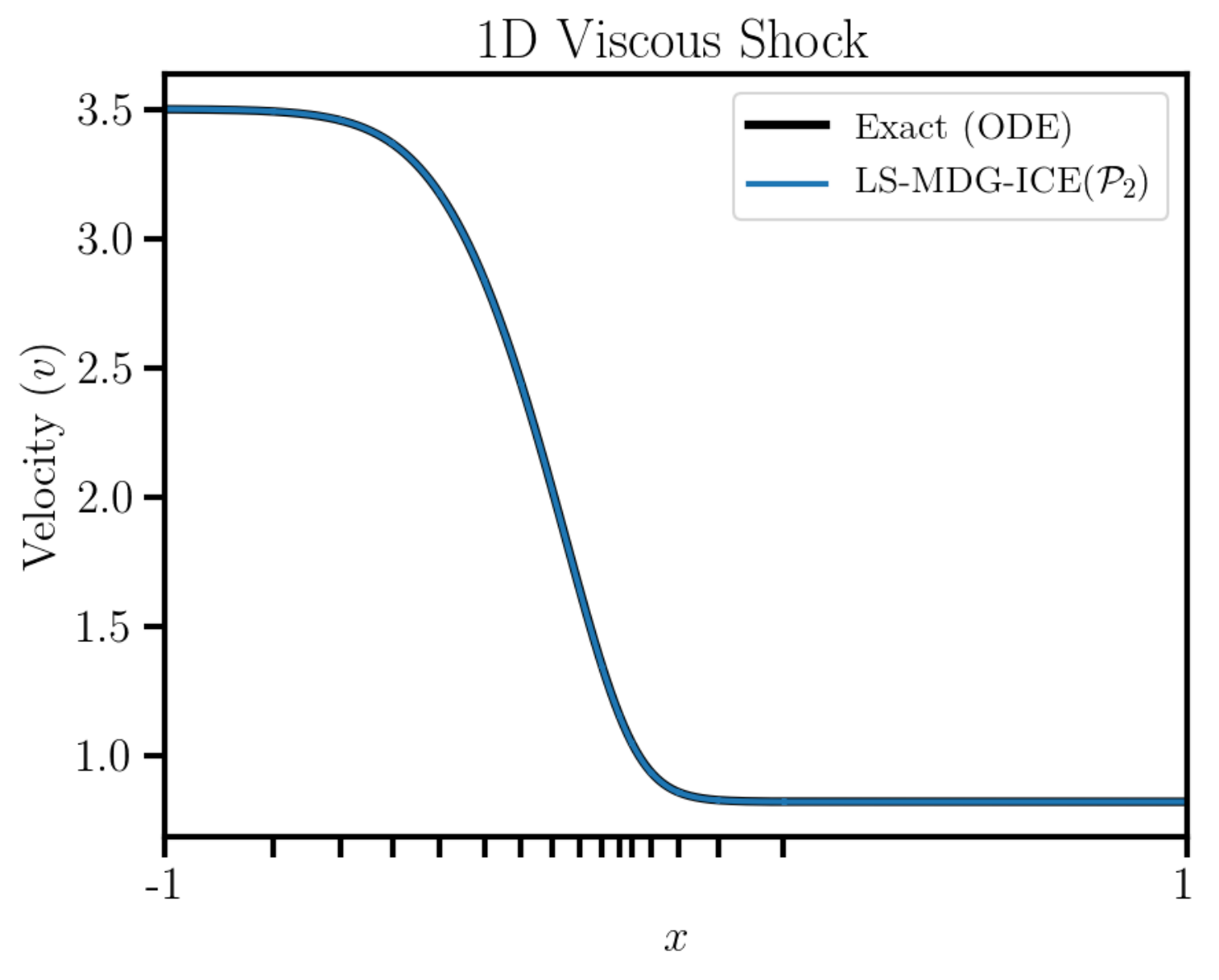} & \includegraphics[width=0.3\linewidth]{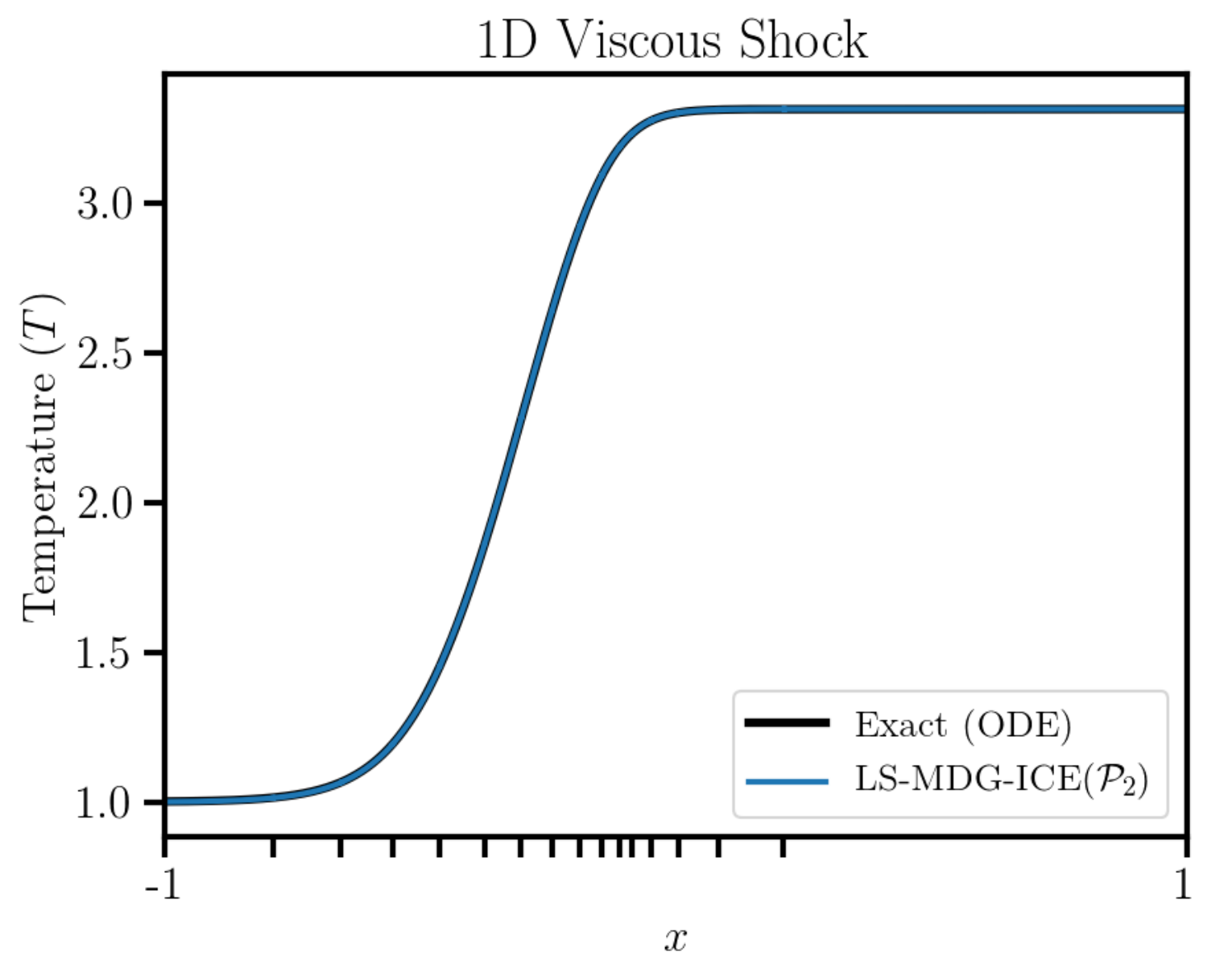}\tabularnewline
\includegraphics[width=0.3\linewidth]{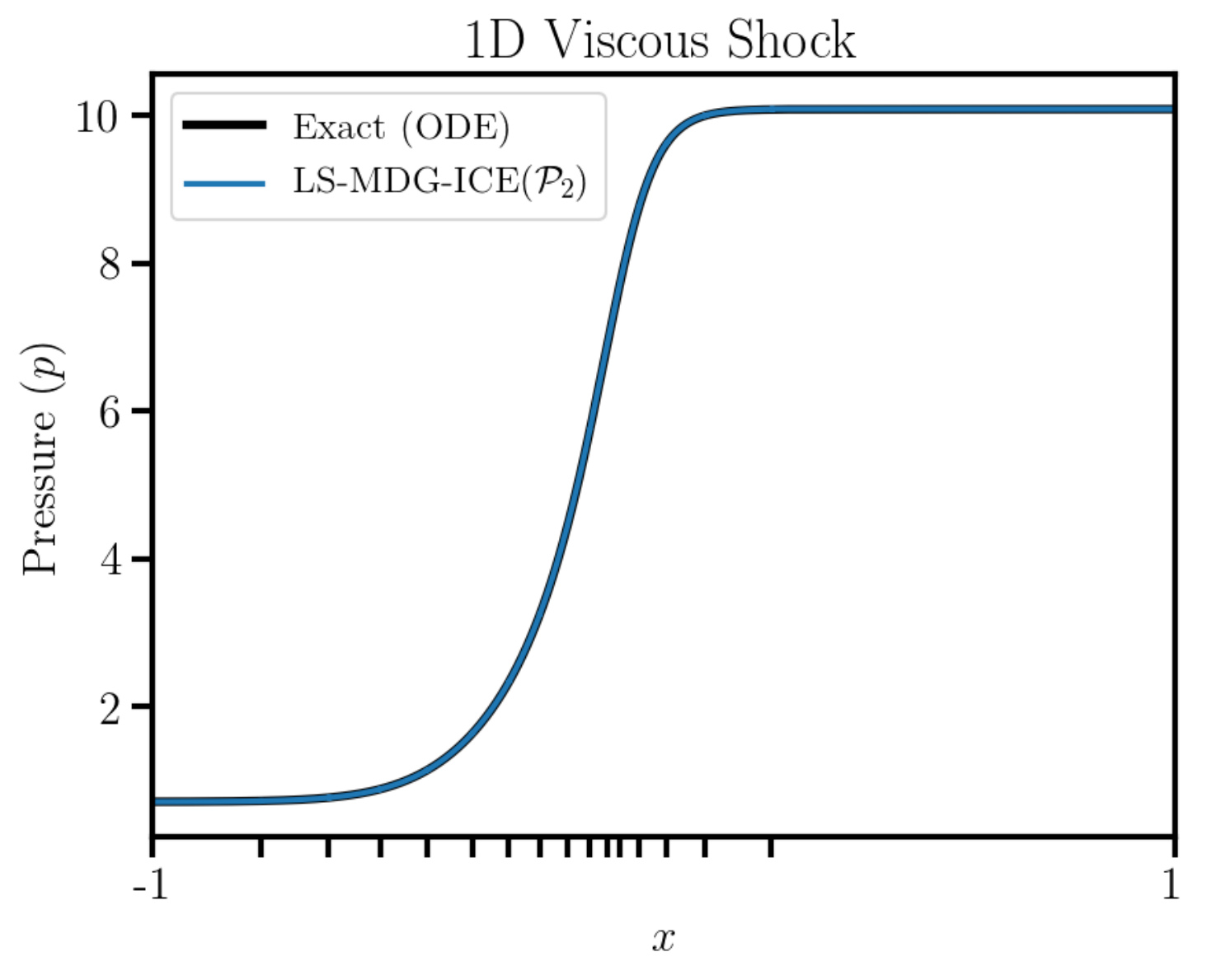} & \includegraphics[width=0.3\linewidth]{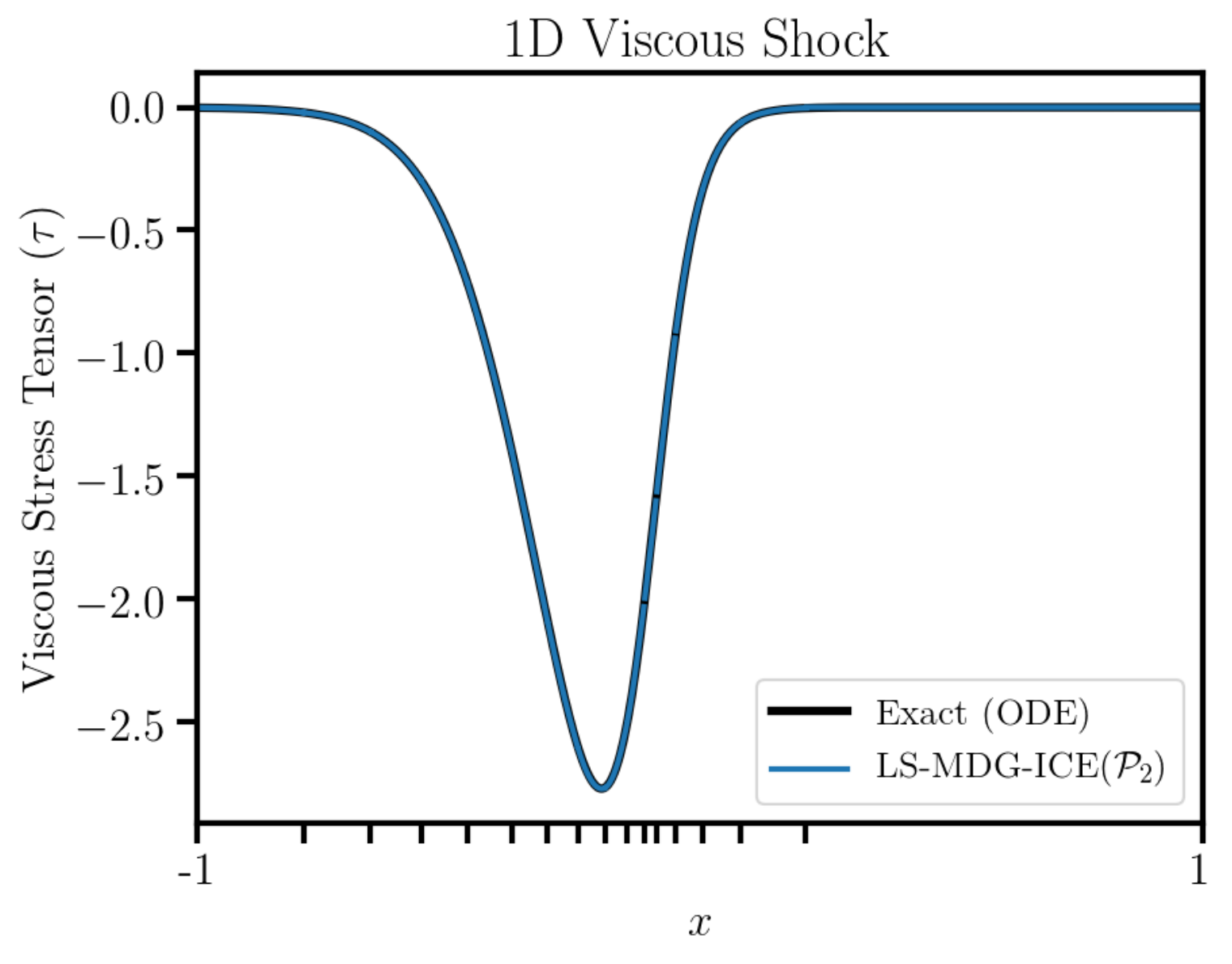} & \includegraphics[width=0.3\linewidth]{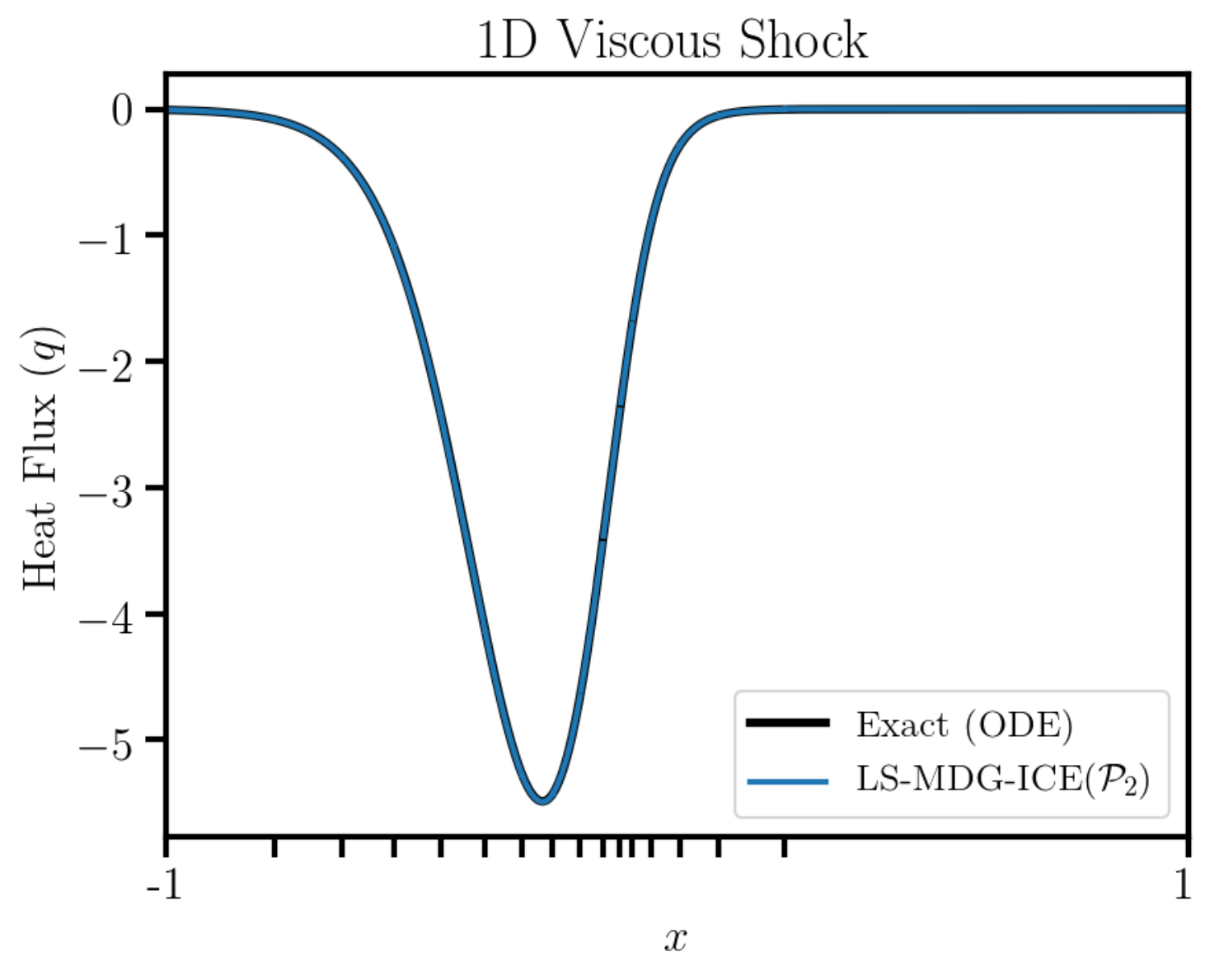}\tabularnewline
\end{tabular}\caption{\label{fig:CNS-Viscous-Shock} The LS-MDG-ICE($\mathcal{P}_{2}$)
computed on 16 isoparametric line cells for steady viscous Mach 3.5
shock at $\mathrm{Re}=25$, $\mathrm{Pr}=0.72$, and $T_{\infty}=293.15$
.}
\end{figure}
Figure~\ref{fig:CNS-Viscous-Shock} presents the LS-MDG-ICE($\mathcal{P}_{2}$)
solution computed on 16 isoparametric line cells and the ODE approximation
of the exact solution computed using fourth-order ODE integration
with a step size of $10^{-4}$ for a steady viscous Mach 3.5 shock
at $\mathrm{Re}=25$, $\mathrm{Pr}=0.72$, and $T_{\infty}=293.15$.
The location of the viscous shock corresponding to the initial conditions
given by Equation~(\ref{eq:CNS-Viscous-Shock-Initial-Conditions})
is not uniquely determined. The numerical approximation to the set
of ODEs describing the exact solution was shifted by $x=-0.14$ and
extrapolated to the right boundary, $x=1$, for direct comparison
with the LS-MDG-ICE solution. Figure~\ref{fig:CNS-Viscous-Shock}
also shows the final grid projected to the space $\mathcal{P}_{1}$,
which adapted to resolve multiple profiles simultaneously, providing
the highest resolution at the transition from the viscous shock to
the right state where oscillations generated by traditional methods
based on polynomial approximations are expected to be the greatest.

\subsection{Sinusoidal Waves\label{subsec:sinusoidal-waves-convergence}}

\begin{figure}
\subfloat[\label{fig:entropy-waves-convergence-velocity-weighted-0_1}Spatial
velocity $v_{x}=\frac{1}{10}$. Reproduction of the convergence results
originally presented in our previous work using a discrete least-squares
solver with an equal-order test space~\citep{Cor18} ]{\centering{}\includegraphics[width=0.48\linewidth]{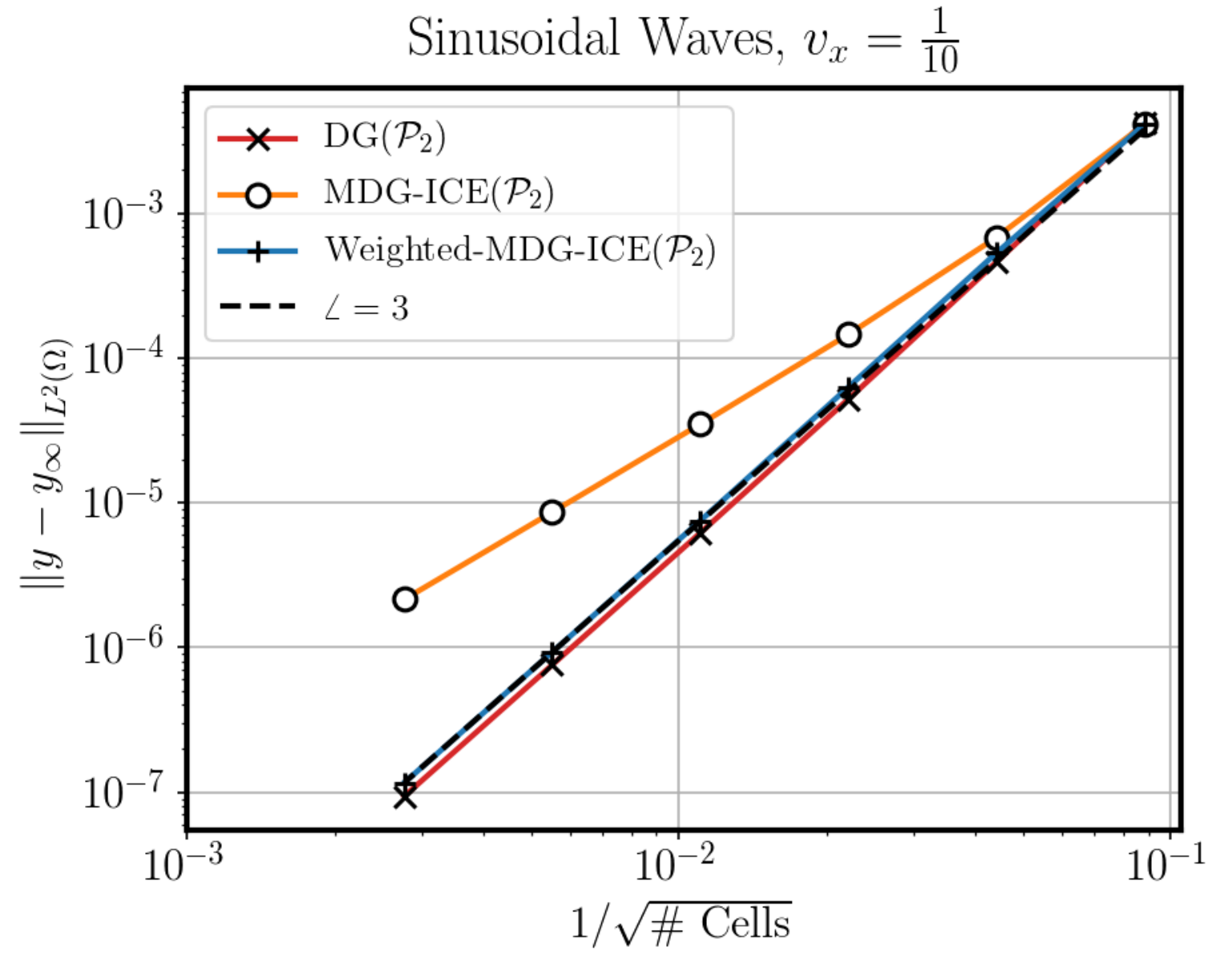}}\subfloat[\label{fig:entropy-waves-convergence-velocity-least-squares-0_1}Spatial
velocity $v_{x}=\frac{1}{10}$. Optimal-order convergence results
obtained using LS-MDG-ICE($\mathcal{P}_{2}$)]{\centering{}\includegraphics[width=0.48\linewidth]{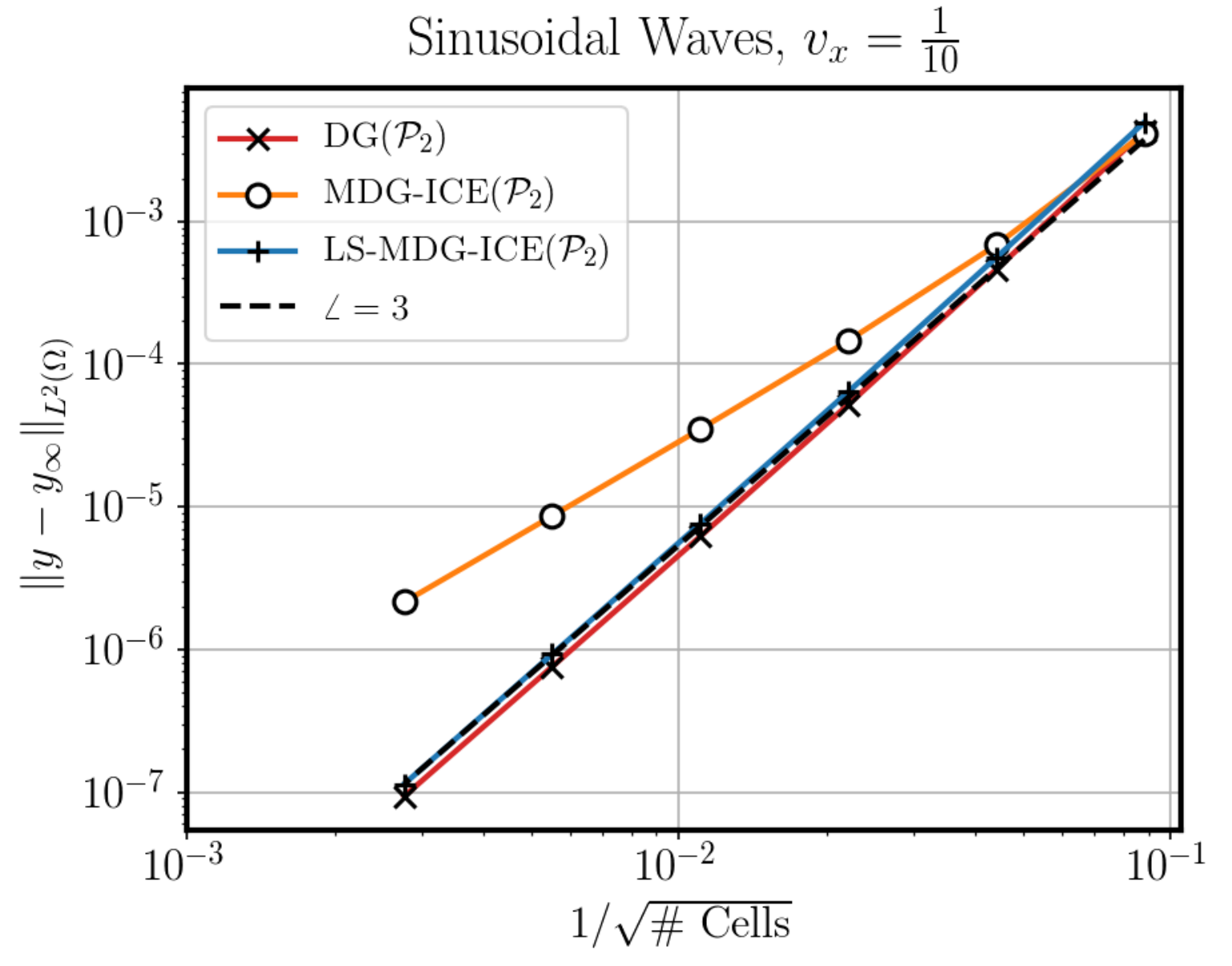}}

\caption{\label{fig:entropy-waves-convergence-velocity}Convergence plots for
linear advection of a smooth sinusoidal wave comparing a standard
DG($\mathcal{P}_{2}$) method with an MDG-ICE($\mathcal{P}_{2}$)
method using the discrete weak formulation given by Equation~(6.5)
in~\citep{Cor18}. The coarsest grid consisted of $128$ triangle
cells, while the finest grid consisted of $131,072$ triangle cells.
For this smooth flow the discrete geometry $u$ is fixed. For the
case of a spatial velocity $v_{x}=\frac{1}{10}$, unweighted MDG-ICE($\mathcal{P}_{2}$)
converges at a sub-optimal rate, while DG($\mathcal{P}_{2}$), weighted
MDG-ICE($\mathcal{P}_{2}$), and LS-MDG-ICE($\mathcal{P}_{2}$) achieve
optimal third-order convergence.}
\end{figure}

We revisit two-dimensional linear advection of a sinusoidal wave presented
in~\citep{Cor18} and demonstrate that LS-MDG-ICE achieves optimal-order
convergence with respect to grid refinement on a static, two-dimensional
grid. This problem involves a smooth flow defined over the space-time
domain $\Omega=\left(0,2\right)\times\left(0,2\right)$, governed
by the linear advection equation, cf. Section~\ref{subsec:Linear-advection-diffusion},
with a spatial velocity $v_{x}>0$ and diffusion coefficient $\epsilon=0$.
The exact solution

\begin{equation}
y\left(x,t\right)=\frac{7}{5}\left(1+\frac{1}{10}\sin\left(2\pi\left(x-v_{x}t\right)\right)\right),\label{eq:sinusoidal-waves-exact-solution}
\end{equation}
is used to define temporal inflow (initial) and spatial inflow boundary
conditions at $t=0$ and $x=0$.

Figure~\ref{fig:entropy-waves-convergence-velocity} presents convergence
results with respect to grid refinement for the case of $v_{x}=0.1$.
The initial coarse grid consisted of $128$ triangle cells and was
refined successively up to $131,072$ triangle cells. Figure~\ref{fig:entropy-waves-convergence-velocity-weighted-0_1}
reproduces the convergence results corresponding to DG($\mathcal{P}_{2}$),
MDG-ICE($\mathcal{P}_{2}$), and weighted MDG-ICE($\mathcal{P}_{2}$)
originally presented in our earlier work~\citep{Cor18}. In this
case, unweighted MDG-ICE($\mathcal{P}_{2}$) converges at a sub-optimal
rate, while both DG($\mathcal{P}_{2}$) and weighted MDG-ICE($\mathcal{P}_{2}$)
achieve optimal third-order convergence. Optimal third-order convergence
is recovered for MDG-ICE if the interface term is scaled by a coefficient
$\alpha_{h}\in\mathbb{R}$, which, for this particular problem, is
defined as
\begin{equation}
\alpha_{h}=2^{i},\label{eq:interface-residual-weight}
\end{equation}
for $i=0,\cdots,n$, where $n$ is the number of times the solution
was refined. Thus, $\alpha_{h}$ is one on the coarsest level and
increases by a factor of two upon each level of refinement. This formulation
is referred to as weighted-MDG-ICE($\mathcal{P}_{2}$) in Figure~\ref{fig:entropy-waves-convergence-velocity-weighted-0_1}.

The interface weight $\alpha_{h}$ that will recover optimal-order
convergence is unknown in general. In contrast, the least-squares
formulation uses optimal test functions that are generated automatically
and do not require an ad hoc scaling. Figure~\ref{fig:entropy-waves-convergence-velocity-least-squares-0_1}
replaces weighted MDG-ICE with LS-MDG-ICE. As expected, LS-MDG-ICE
converges at the optimal rate and is comparable to DG($\mathcal{P}_{2}$)
in terms of accuracy. The error with respect to the exact solution
on the finest grid is $1.14\times10^{-7}$ for LS-MDG-ICE and $9.36\times10^{-8}$
DG($\mathcal{P}_{2}$).

\subsection{Inviscid Mach 3 bow shock\label{sec:inviscid-bow-shock}}

\begin{figure}
\begin{centering}
\subfloat[\label{fig:Inviscid-Bow-Shock-Mesh}392 isoparametric $\mathcal{P}_{2}$
triangle cells]{\begin{centering}
\includegraphics[width=0.9\linewidth]{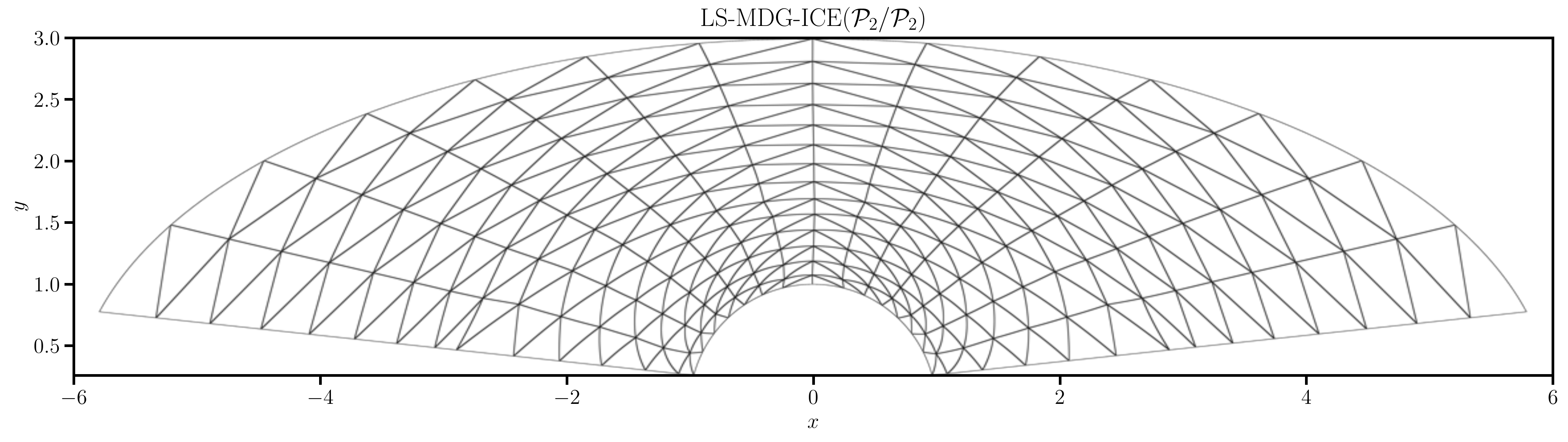} 
\par\end{centering}
}
\par\end{centering}
\begin{centering}
\subfloat[\label{fig:Inviscid-Bow-Shock-Mach}392 isoparametric $\mathcal{P}_{2}$
triangle cells]{\begin{centering}
\includegraphics[width=0.9\linewidth]{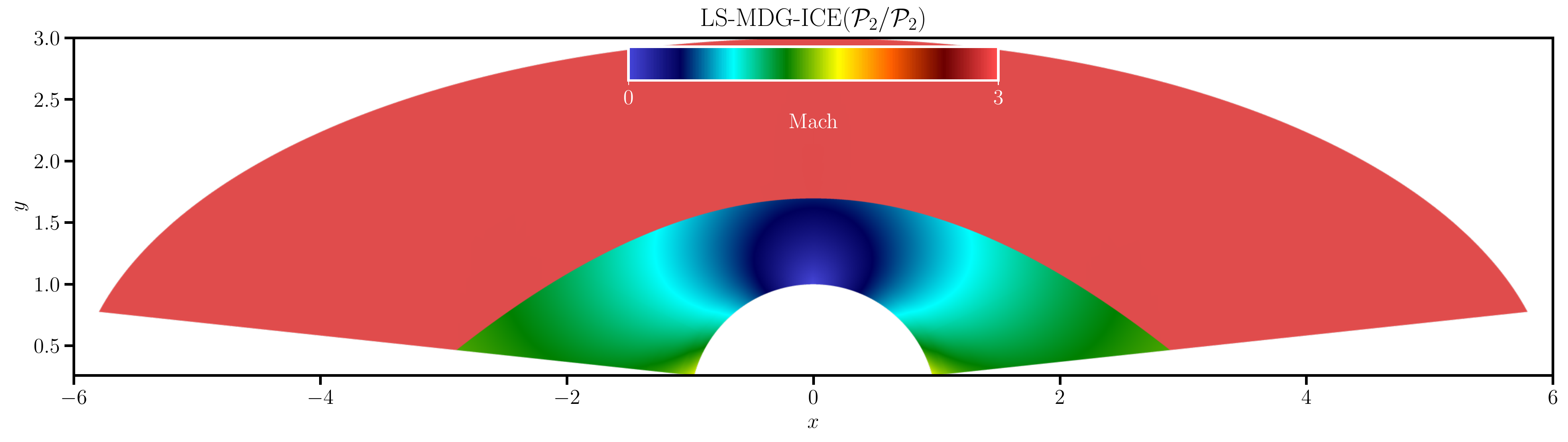}
\par\end{centering}
}
\par\end{centering}
\begin{centering}
\subfloat[\label{fig:Inviscid-Bow-Shock-Pressure}392 isoparametric $\mathcal{P}_{2}$
triangle cells]{\begin{centering}
\includegraphics[width=0.9\linewidth]{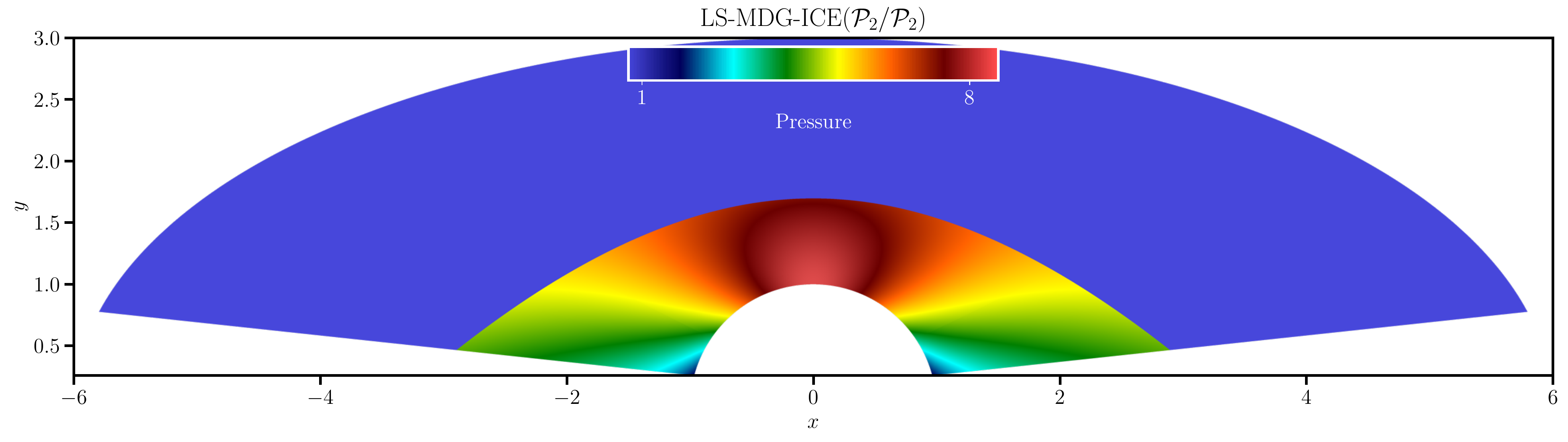}
\par\end{centering}
}
\par\end{centering}
\begin{centering}
\subfloat[\label{fig:Inviscid-Bow-Shock-EnthalpyStagnation}392 isoparametric
$\mathcal{P}_{2}$ triangle cells]{\begin{centering}
\includegraphics[width=0.9\linewidth]{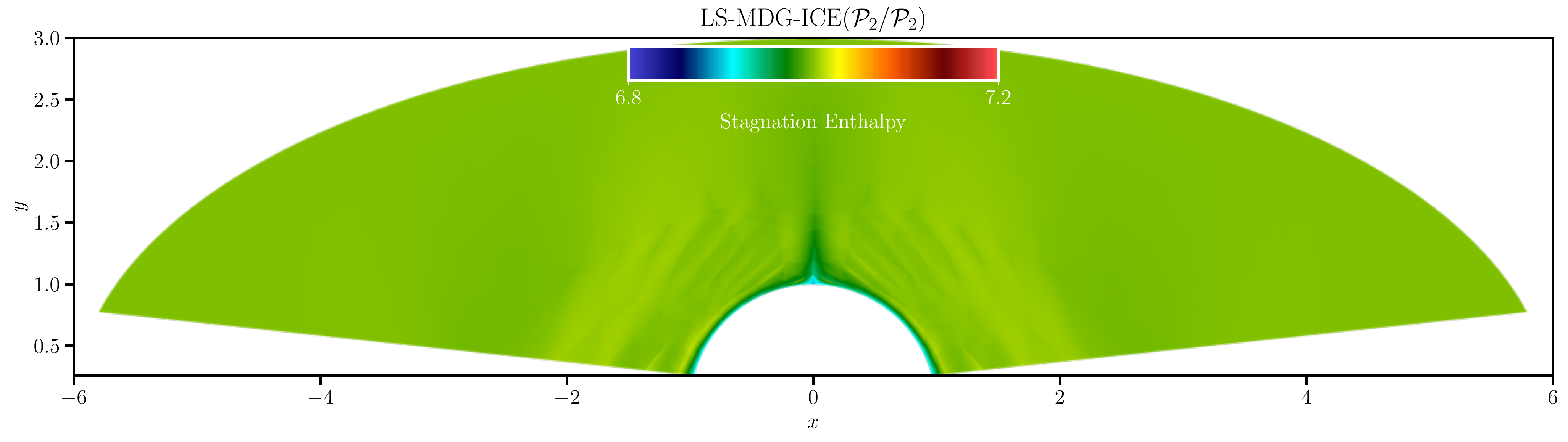}
\par\end{centering}
}
\par\end{centering}
\centering{}\caption{\label{fig:Inviscid-Bow-Shock-MDG-ICE}Inviscid Mach 3 bow shock:
The LS-MDG-ICE solution computed using $\mathcal{P}_{2}$ isoparametric
triangle elements. The LS-MDG-ICE field variables were initialized
by projecting the cell averaged DG($\mathcal{P}_{2}$) solution. The
location of the shock along the line $x=0$ was computed as $y=1.698$
for a stand-off distance of $0.698$. The error of the stagnation
enthalpy, $\left\Vert H-H_{\infty}\right\Vert _{L^{2}\left(\Omega\right)}$,
with respect to the exact solution, $H_{\infty}=7$, is $2.2736\times10^{-2}$.}
\end{figure}
\begin{figure}
\subfloat[\label{fig:Inviscid-Bow-Shock-1D-Pressure}The pressure sampled along
$x=0$. The exact pressure at the stagnation point, $p\approx8.614975$,
is marked with the symbol $\times$. The computed pressure at the
stagnation point on is $8.585034$.]{\begin{centering}
\includegraphics[width=0.45\linewidth]{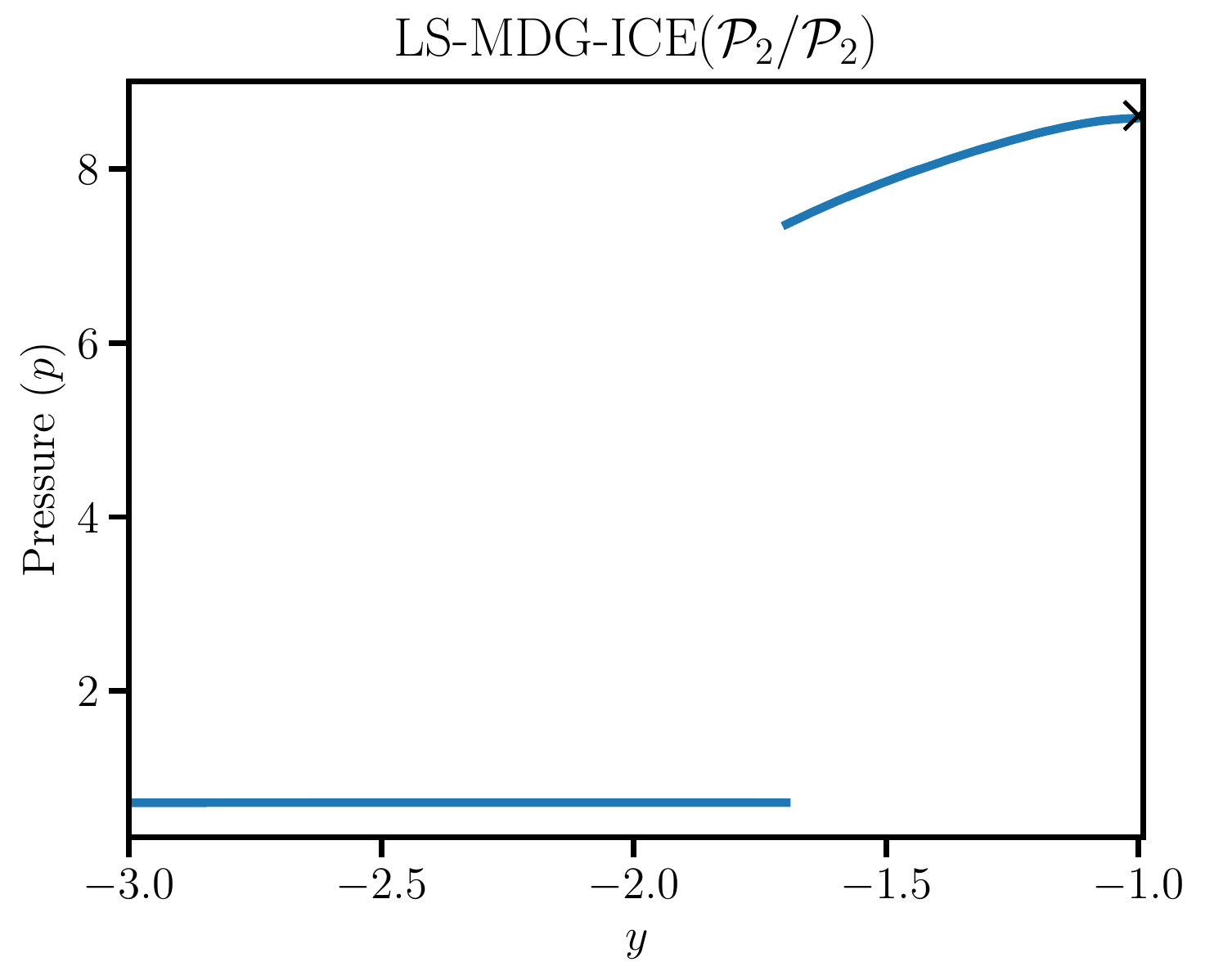}
\par\end{centering}
}\hfill{}\subfloat[\label{fig:Inviscid-Bow-Shock-1D-Mach}The Mach number sampled along
$x=0$. The exact Mach number at the stagnation point, $M=0$, is
marked with the symbol $\times$. The computed Mach number at the
stagnation point is $6.86725\times10^{-3}$.]{\begin{centering}
\includegraphics[width=0.45\linewidth]{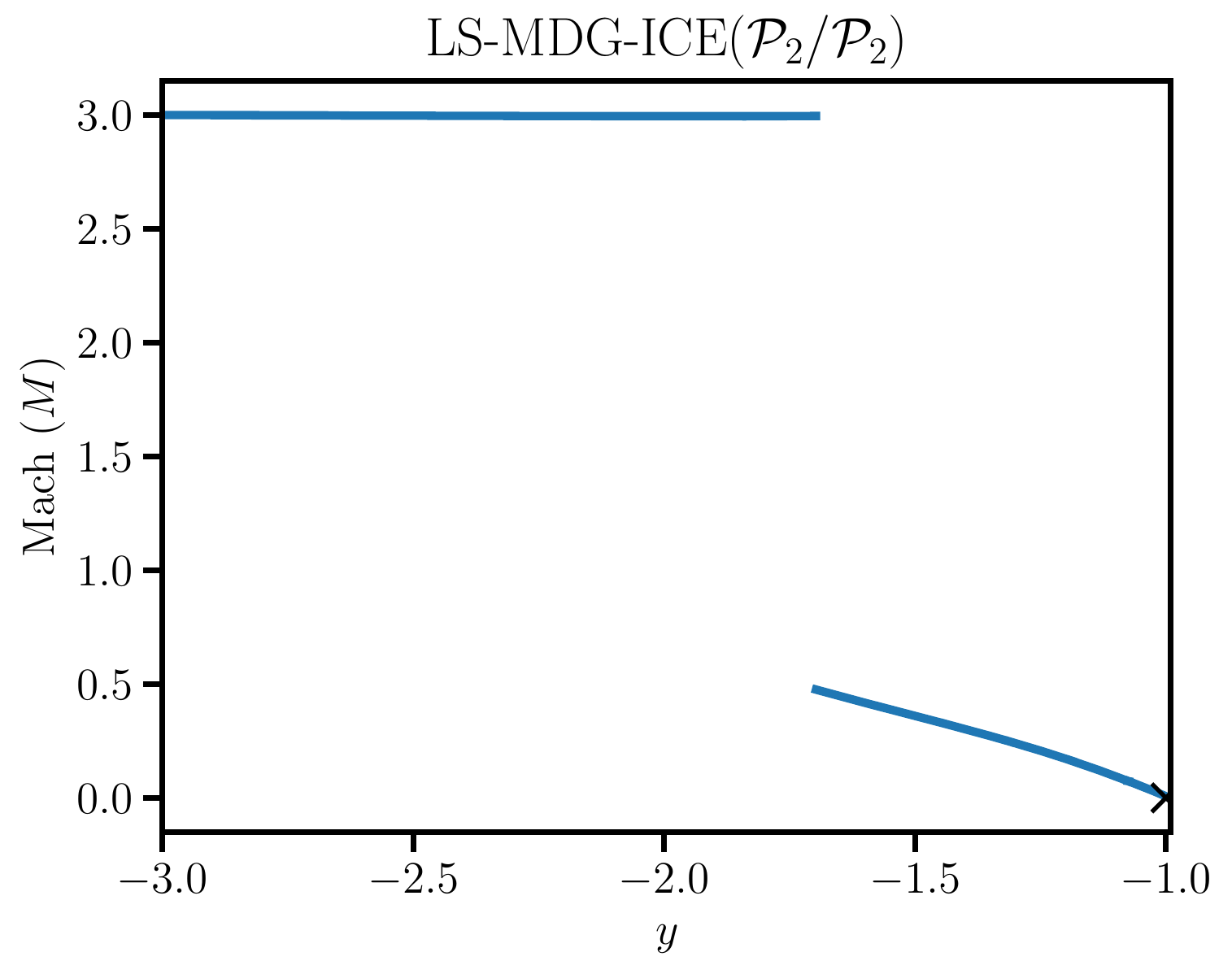}
\par\end{centering}
}

\caption{\label{fig:Inviscid-Bow-Shock-1D}Inviscid Mach 3 bow shock: Computed
with LS-MDG-ICE using $\mathcal{P}_{2}$ isoparametric triangle elements.
The location of the shock was computed as $y=1.698$ for a stand-off
distance of $0.698$.}
\end{figure}
The least-squares formulation was applied to the steady inviscid Mach
3 bow shock. Figure~\ref{fig:Inviscid-Bow-Shock-MDG-ICE} presents
the final mesh, temperature, pressure, and stagnation enthalpy fields
corresponding to the MDG-ICE($\mathcal{P}_{2}$) solution computed
using $392$ $\mathcal{P}_{2}$ isoparametric triangle elements. The
field variables were initialized by projecting the cell averaged DG($\mathcal{P}_{2}$)
solution. The error of the stagnation enthalpy, $\left\Vert H-H_{\infty}\right\Vert _{L^{2}\left(\Omega\right)}$,
with respect to the exact solution, $H_{\infty}=7$, was reduced from
$4.0932\times10^{-1}$ to $2.2736\times10^{-2}$. The stand off distance
was estimated to be to $0.698$, which is in agreement with the estimate
computed using the discrete least-squares formulation reported previously~\citep{Cor18}.

Figure~\ref{fig:Inviscid-Bow-Shock-1D} presents the pressure and
Mach number sampled along the centerline, $x=0$. The computed pressure
at the stagnation point is $8.585034$\footnote{The discrepancy in the pressure scale, in particular, the pressure
at the stagnation point, $8.585034$, reported in this work and the
value reported in~\citep{Cor18}, $12.0559064$, is due to the nondimensionalization
of the freestream quantities. In this work, we nondimensionalize by
a term proportional to the freestream speed of sound squared that
results in a nondimensional freestream pressure of $p_{\infty}=\nicefrac{1}{\gamma}$.
In~\citep{Cor18}, the nondimensional freestream pressure was $p_{\infty}=1$.
Rescaling the computed pressure to match the scaling in~\citep{Cor18}
results in a pressure at the stagnation point of $12.019048$.}. The exact pressure at the stagnation point, $p\approx8.614975$,
is marked with the symbol $\times$ in Figure~\ref{fig:Inviscid-Bow-Shock-1D-Pressure}.
The computed Mach number at the stagnation point is $M=6.86725\times10^{-3}$.
The exact Mach number at the stagnation point, $M=0$, is marked with
the symbol $\times$ in Figure~\ref{fig:Inviscid-Bow-Shock-1D-Pressure}

\subsection{Viscous Mach 5 bow shock\label{sec:viscous-bow-shock}}

To verify the viscous LS-MDG-ICE discretization of compressible Navier-Stokes
flow, as described in Section~\ref{subsec:Compressible-Navier-Stokes-flow},
and demonstrate the ability of the LS-MDG-ICE to simultaneously resolve
both the interior viscous shock and viscous boundary layer via anisotropic
curvilinear $r$-adaptivity. We consider Mach 5 flows at a Reynolds
numbers of $10^{3}$ and $10^{4}$ and extend the one-dimensional
analysis from~\citep{Ker20} to the case of $\mathrm{Re}=10^{4}$.

\begin{figure}
\begin{centering}
\subfloat[\label{fig:Viscous-Bow-Shock-DG-Mesh}392 $\mathcal{P}_{1}$-shape
triangle cells]{\begin{centering}
\includegraphics[width=0.9\linewidth]{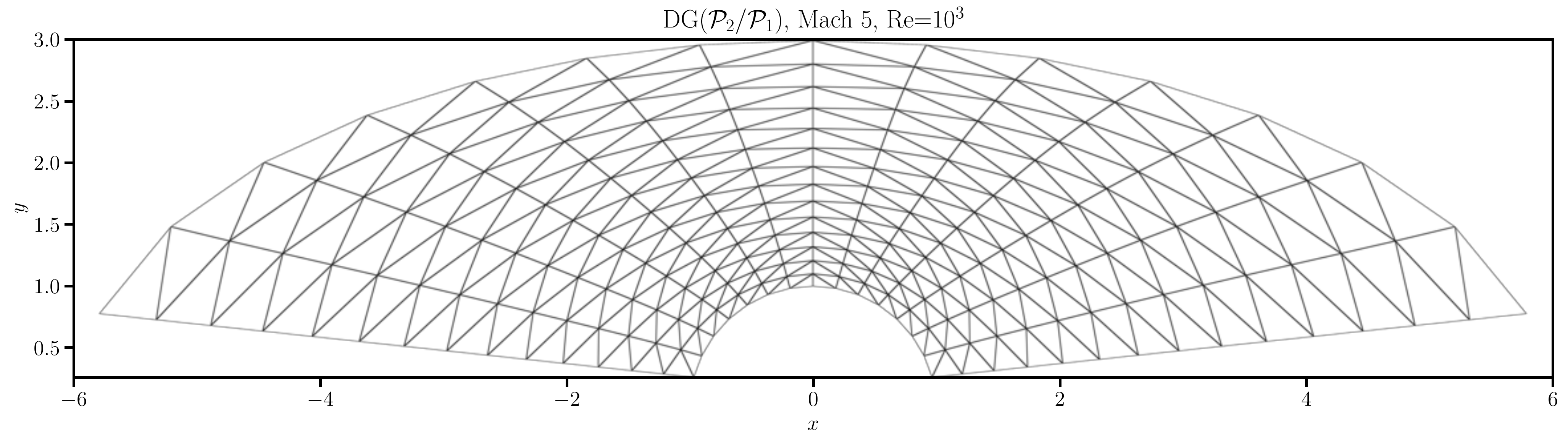}
\par\end{centering}
}
\par\end{centering}
\begin{centering}
\subfloat[\label{fig:Viscous-Bow-Shock-DG-Temperature}392 linear $\mathcal{P}_{2}$
triangle cells]{\begin{centering}
\includegraphics[width=0.9\linewidth]{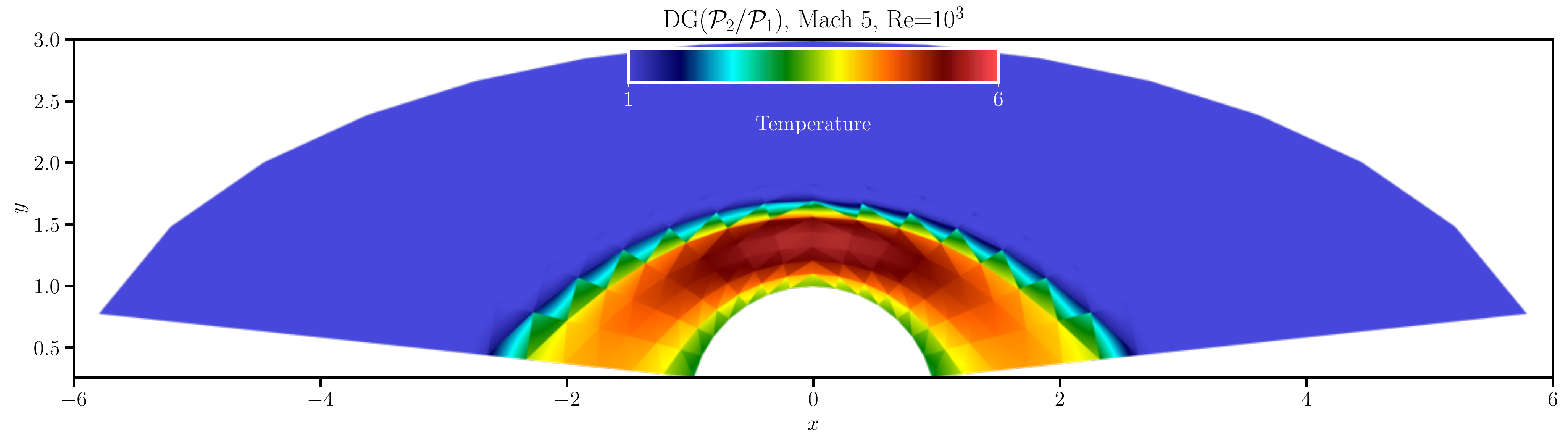}
\par\end{centering}
}
\par\end{centering}
\caption{\label{fig:Viscous-Bow-Shock-DG}The initial linear grid and temperature
field corresponding to a shock captured DG($\mathcal{P}_{2}/\mathcal{P}_{1}$)
solution for the viscous Mach 5 bow shock at $\mathrm{Re}=10^{3}$.}
\end{figure}

A supersonic viscous flow over a cylinder in two dimensions is characterized
by the Reynolds number $\mathrm{Re}$, and the freestream Mach number
$M_{\infty}$. The Reynolds number is defined as,

\begin{equation}
\mathrm{Re}=\frac{\rho vL}{\mu},\label{eq:reynolds-number}
\end{equation}
where $L$ is the characteristic length. The freestream Mach number
is defined $M_{\infty}=\frac{v_{\infty}}{c_{\infty}}$ 
\begin{equation}
M_{\infty}=\frac{v_{\infty}}{c_{\infty}},\label{eq:mach-number}
\end{equation}
where $v_{\infty}$ is the freestream velocity, $c_{\infty}=\sqrt{\nicefrac{\gamma P_{\infty}}{\rho_{\infty}}}$
is the freestream speed of sound, $P_{\infty}$ is the freestream
pressure, and $\rho_{\infty}$ is the freestream density. In this
work, we assume the flow is traveling in the negative $y$-direction,
i.e., from top to bottom in Figure~\ref{fig:Viscous-Bow-Shock-DG}
and Figure~\ref{fig:Viscous-Bow-Shock-Re-10000-MDG-ICE}. Supersonic
inflow and outflow boundary conditions are applied on the ellipse
boundary, which is defined by $\left(\nicefrac{x}{3}\right)^{2}+\left(\nicefrac{y}{6}\right)^{2}=1$,
and outflow planes respectively. An isothermal no-slip wall is specified
at the surface of the cylinder of radius $r=1$, centered at the origin.
The temperature at the isothermal wall is given as $T_{\mathrm{wall}}=2.5T_{\infty}$,
where $T_{\infty}$ is the freestream temperature.

We follow the approach described in~\citep{Ker20} and use local
edge refinement to adaptively split highly anisotropic elements within
the viscous structures as they were resolved by LS-MDG-ICE. This simple
approach is highly effective for maintaining the original grid resolution
downstream of the viscous shock, which becomes less resolved as the
solver moves the grid to anisotropically resolve the shock, as well
as for inserting additional degrees of freedom in order to better
resolve the curved shock geometry.

Furthermore, we follow the approach of Zahr et al.~\citep{Zah19}
and incorporate the local cell size into the regularization by scaling
$\lambda_{u}$ of~(\ref{eq:regularization-operator}) by a factor
proportional to the inverse of the element volume. This allowed for
more aggressive grid motion, as highly anisotropic elements were less
likely to degenerate as the grid points were repositioned to resolve
the sharp viscous features. The result is a coarser grid in the region
upstream of the viscous shock where the flow is constant and finer
grid resolution through the viscous shock as compared to the corresponding
solution shown in Figure~3.5~ of~\citep{Ker20}, which did not
weight $\lambda_{u}$ by the inverse cell volume.
\begin{figure}
\begin{centering}
\subfloat[\label{fig:Viscous-Bow-Shock-Mesh}502 sub-parametric $\mathcal{P}_{4}/\mathcal{P}_{3}$
triangle cells]{\begin{centering}
\includegraphics[width=0.9\linewidth]{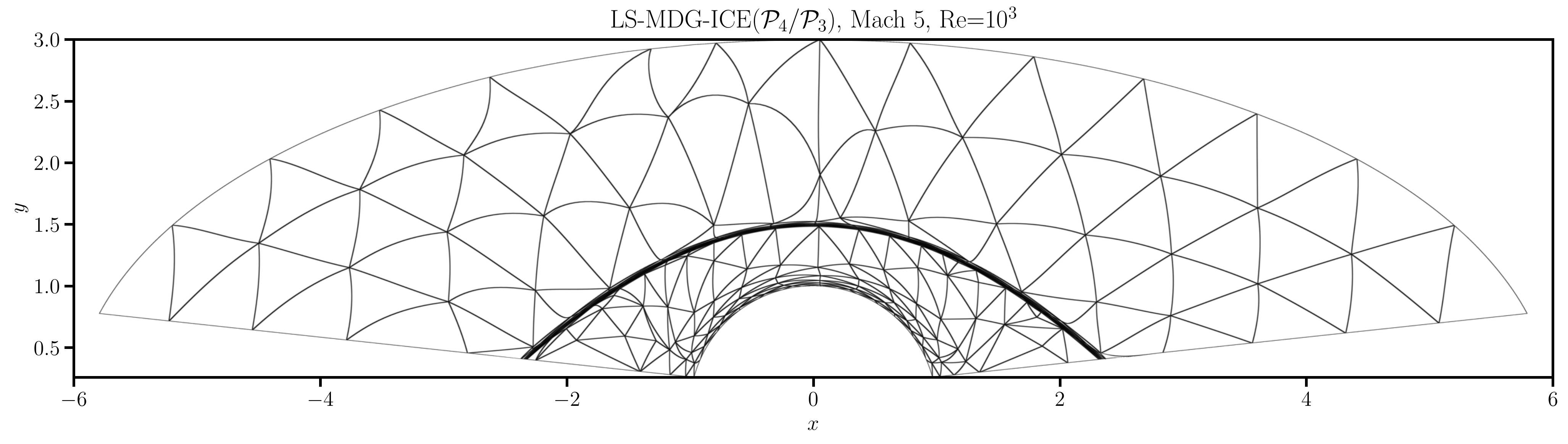}
\par\end{centering}
}
\par\end{centering}
\begin{centering}
\subfloat[\label{fig:Viscous-Bow-Shock-Temperature}502 sub-parametric $\mathcal{P}_{4}/\mathcal{P}_{3}$
triangle cells]{\begin{centering}
\includegraphics[width=0.9\linewidth]{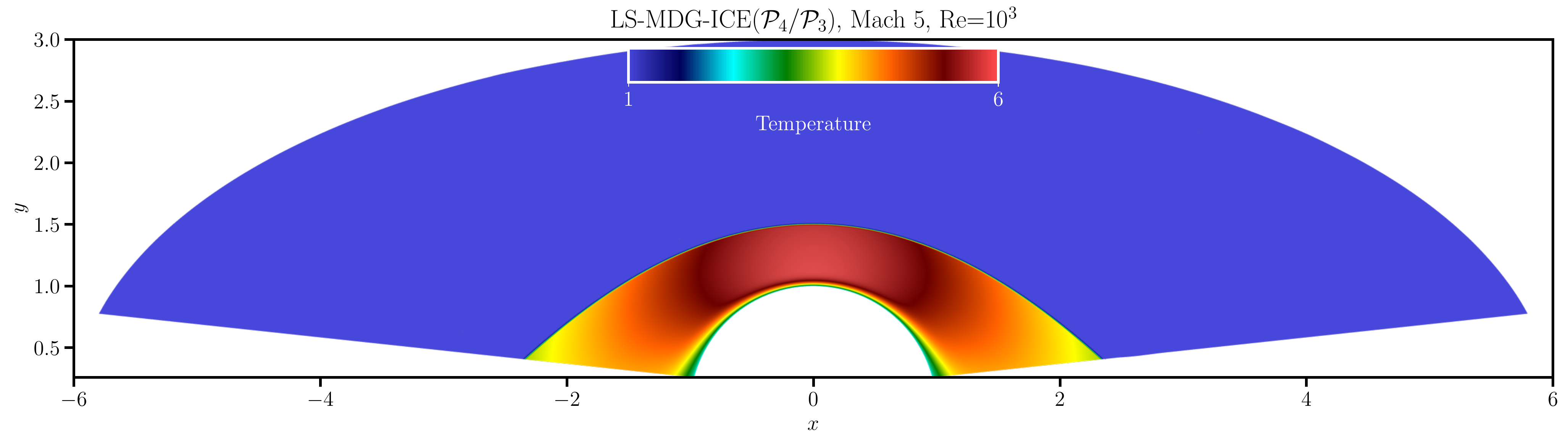}
\par\end{centering}
}
\par\end{centering}
\begin{centering}
\subfloat[\label{fig:Viscous-Bow-Shock-Mach}502 sub-parametric $\mathcal{P}_{4}/\mathcal{P}_{3}$
triangle cells]{\begin{centering}
\includegraphics[width=0.9\linewidth]{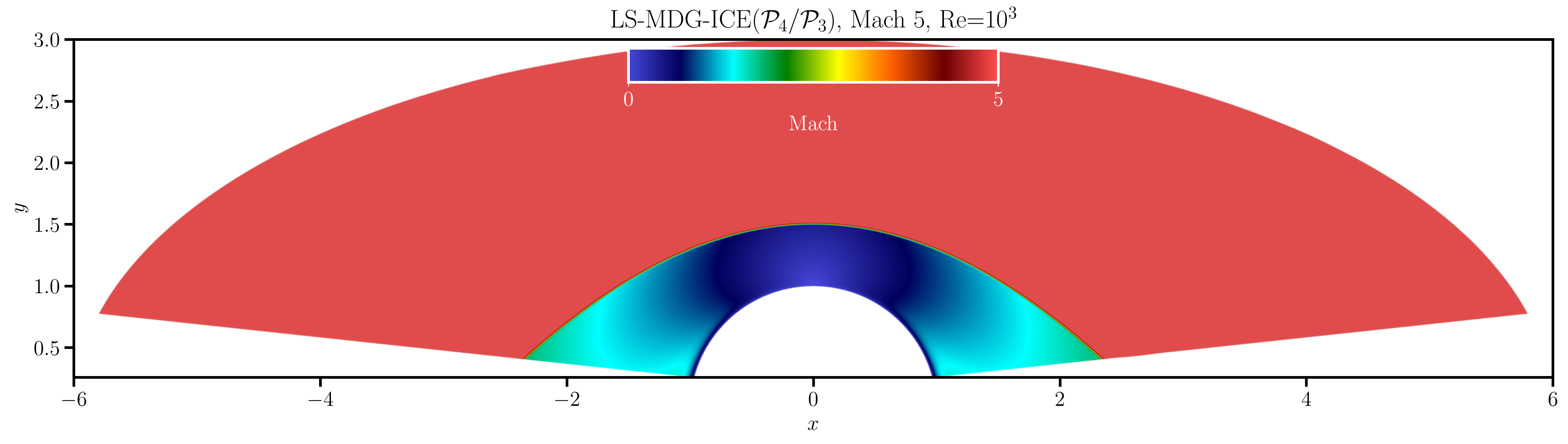}
\par\end{centering}
}
\par\end{centering}
\begin{centering}
\subfloat[\label{fig:Viscous-Bow-Shock-Pressure}502 sub-parametric $\mathcal{P}_{4}/\mathcal{P}_{3}$
triangle cells]{\begin{centering}
\includegraphics[width=0.9\linewidth]{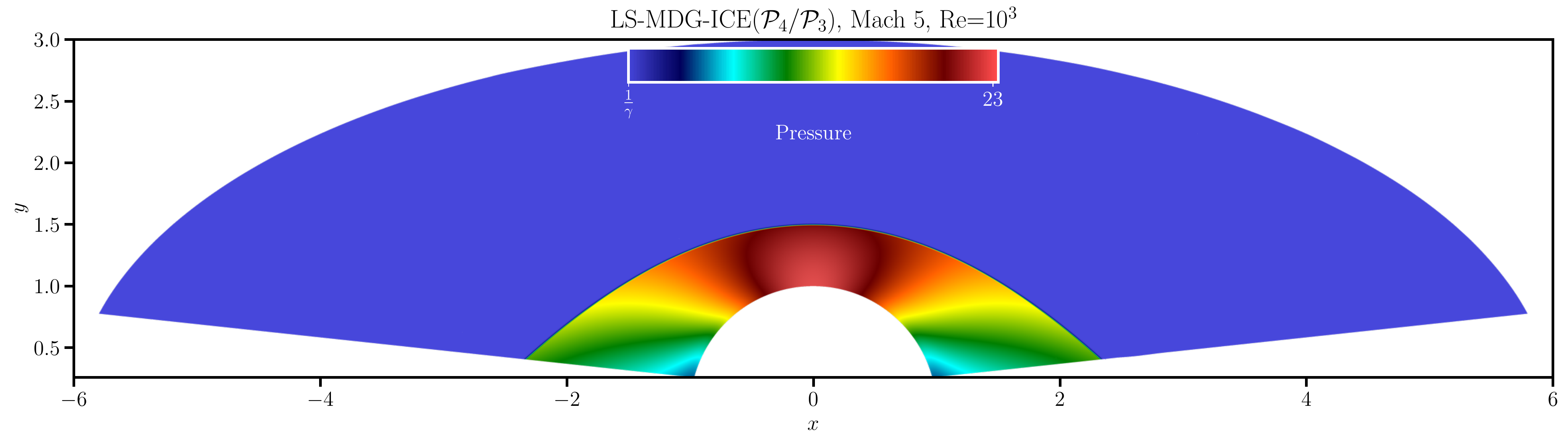}
\par\end{centering}
}
\par\end{centering}
\centering{}\caption{\label{fig:Viscous-Bow-Shock-MDG-ICE}The LS-MDG-ICE solution computed
using $\mathcal{P}_{4}/\mathcal{P}_{3}$ sub-parametric triangle elements
for the viscous Mach 5 bow shock at $\mathrm{Re}=10^{3}$. The LS-MDG-ICE
grid was initialized by projecting the linear triangle grid shown
in Figure~\ref{fig:Viscous-Bow-Shock-DG-Mesh}. The LS-MDG-ICE field
variables were initialized by projecting the cell-averaged DG($\mathcal{P}_{2}/\mathcal{P}_{1}$)
solution shown in Figure~\ref{fig:Viscous-Bow-Shock-DG-Temperature}.
The LS-MDG-ICE flux variables were initialized to zero for consistency
with the initial piecewise constant field variables. The location
of the shock along the line $x=0$ was computed as $y=1.5000575$
for a stand-off distance of $0.5000575$.}
\end{figure}
\begin{figure}
\subfloat[\label{fig:Bow-Shock-1D-Temperature}The exact temperature at the
stagnation point, $T=2.5$, is marked with the symbol $\times$.]{\begin{centering}
\includegraphics[width=0.3\linewidth]{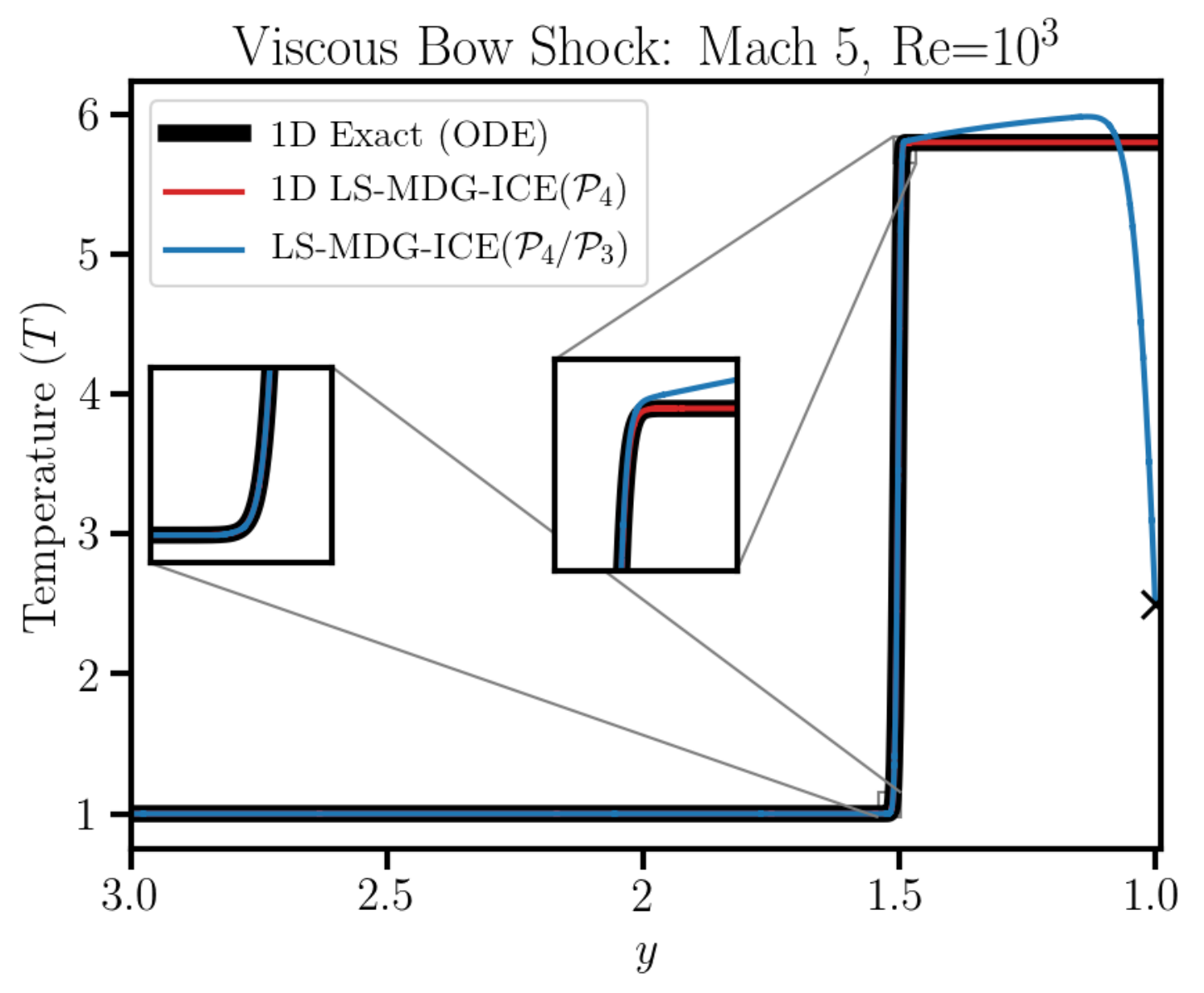}
\par\end{centering}
}\hfill{}\subfloat[\label{fig:Bow-Shock-1D-VelocityX}The exact normal velocity at the
stagnation point, $v_{n}=0$, is marked with the symbol $\times$.]{\begin{centering}
\includegraphics[width=0.3\linewidth]{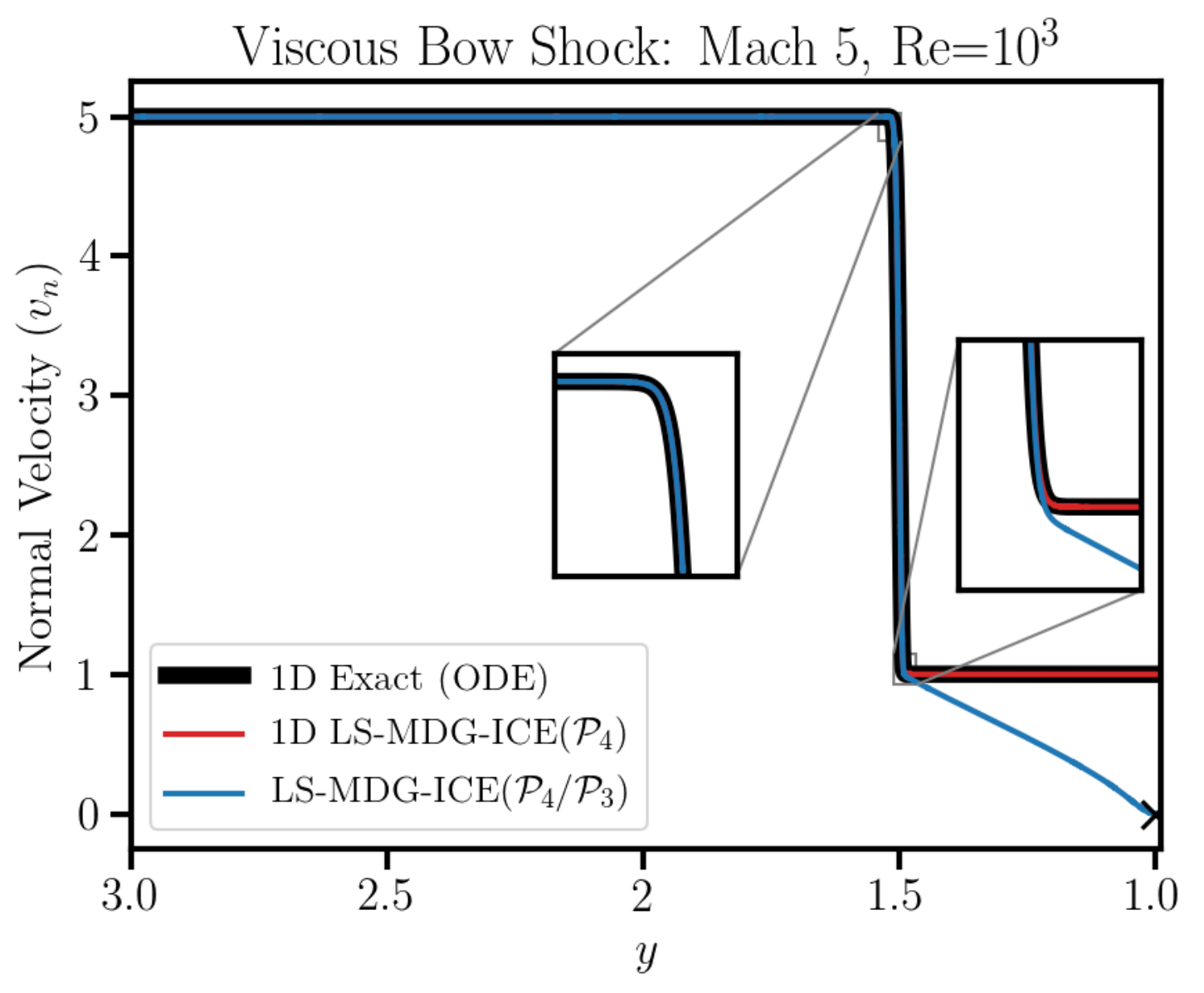}
\par\end{centering}
}\hfill{}\subfloat[\label{fig:Bow-Shock-1D-Pressure}The exact pressure at the stagnation
point for an inviscid flow, $p\approx23.324$, is marked with the
symbol $\times$.]{\begin{centering}
\includegraphics[width=0.3\linewidth]{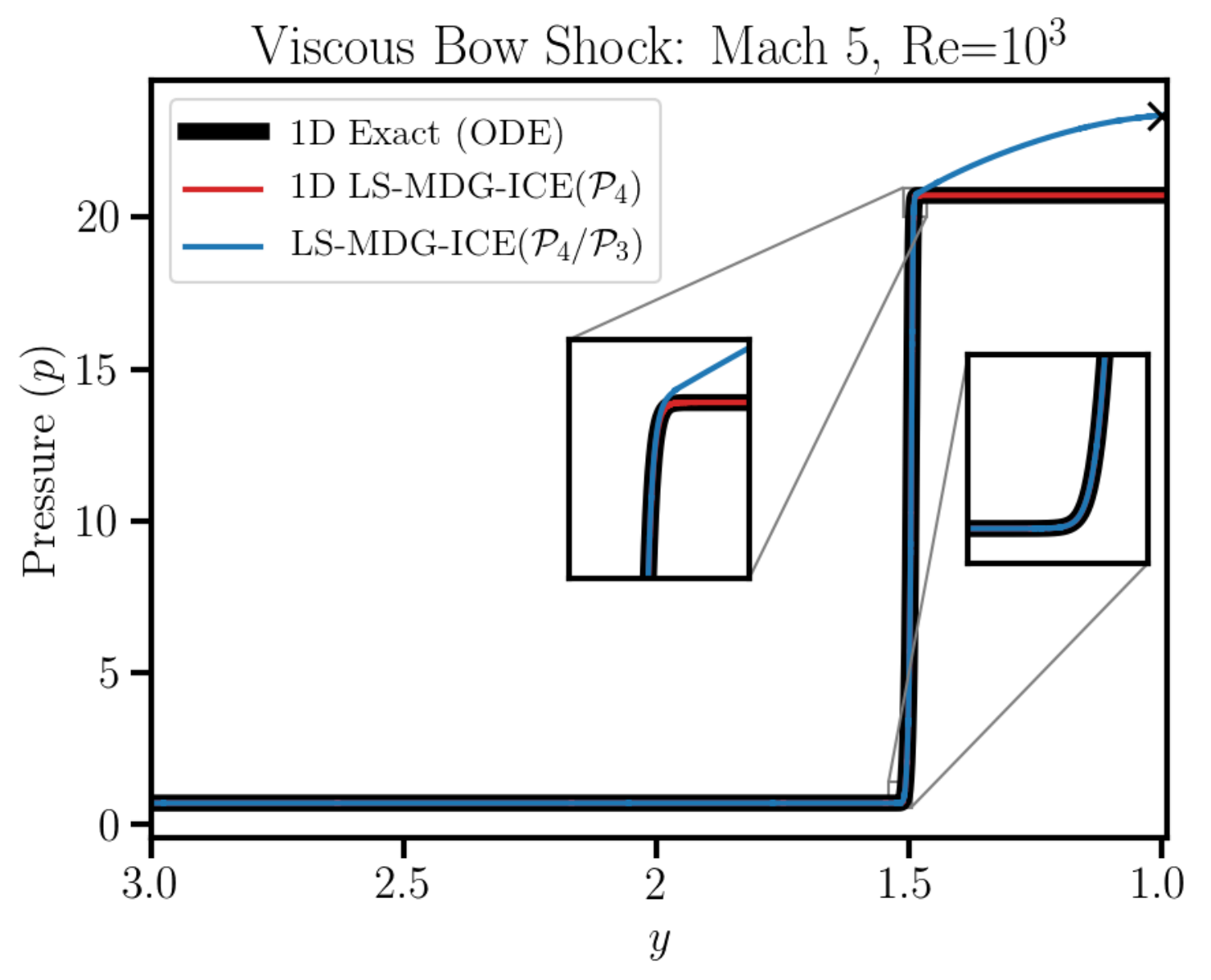}
\par\end{centering}
}

\subfloat[\label{fig:Bow-Shock-1D-Density}The density at the stagnation point,
computed using the stagnation pressure corresponding to an inviscid
flow, $\rho\approx13.0614$, is marked with the symbol $\times$.
]{\begin{centering}
\includegraphics[width=0.3\linewidth]{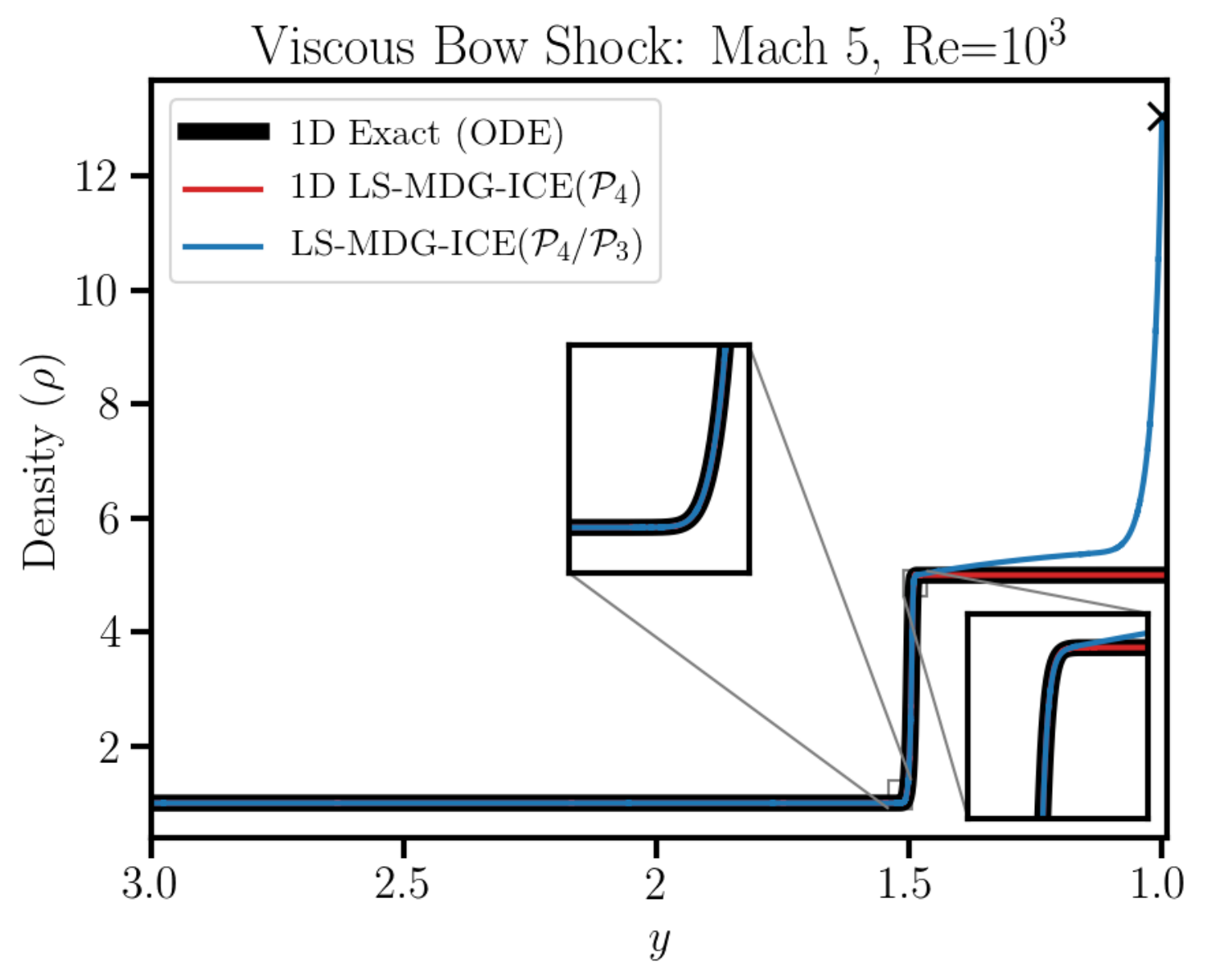}
\par\end{centering}
}\hfill{}\subfloat[\label{fig:Bow-Shock-1D-FluxXMomentumX}The normal component of the
normal viscous stress tensor, $\tau_{nn}$. ]{\begin{centering}
\includegraphics[width=0.3\linewidth]{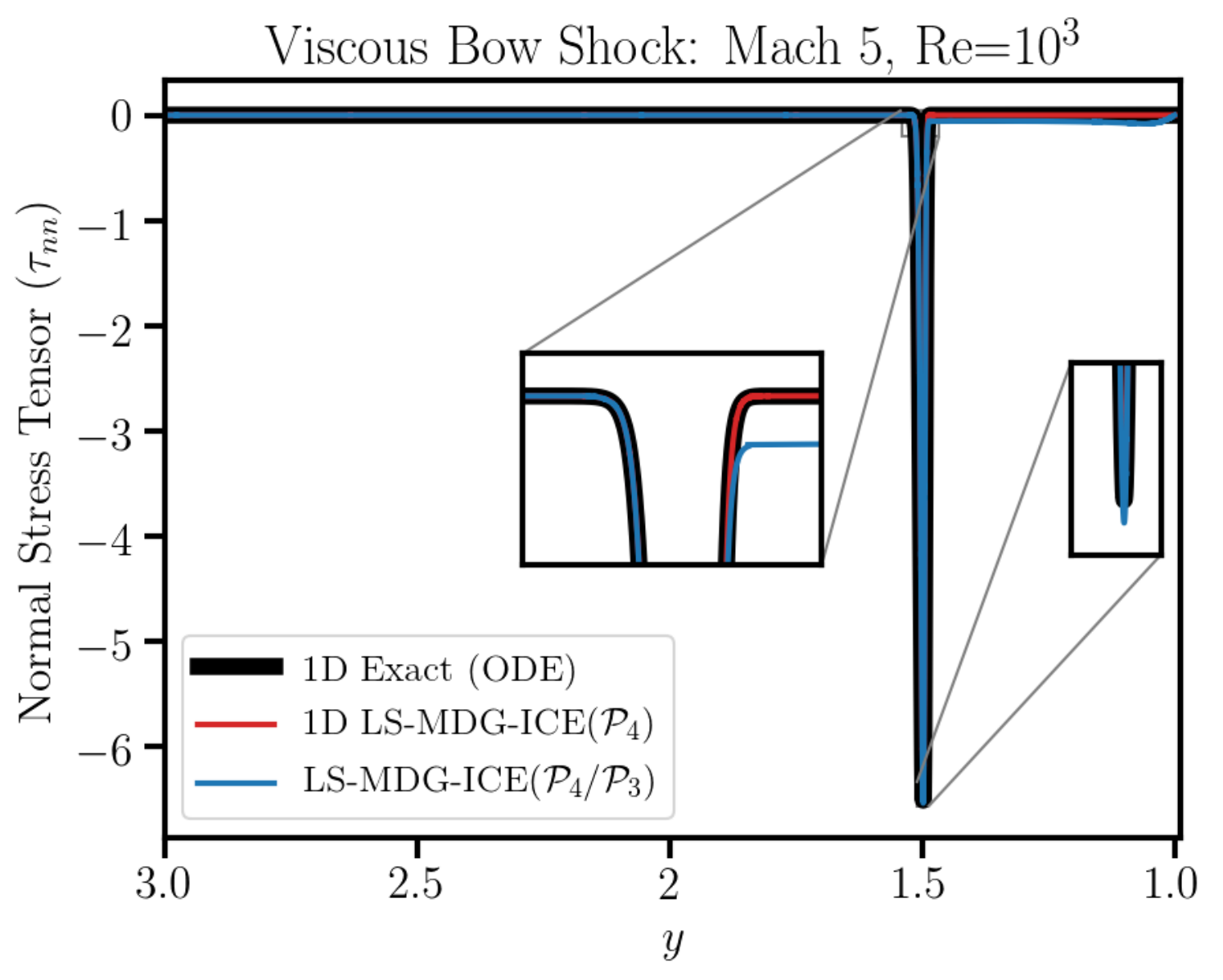}
\par\end{centering}
}\hfill{}\subfloat[\label{fig:Bow-Shock-1D-FluxXEnergyStagnationDensity}The normal thermal
heat flux, $q_{n}$. ]{\begin{centering}
\includegraphics[width=0.3\linewidth]{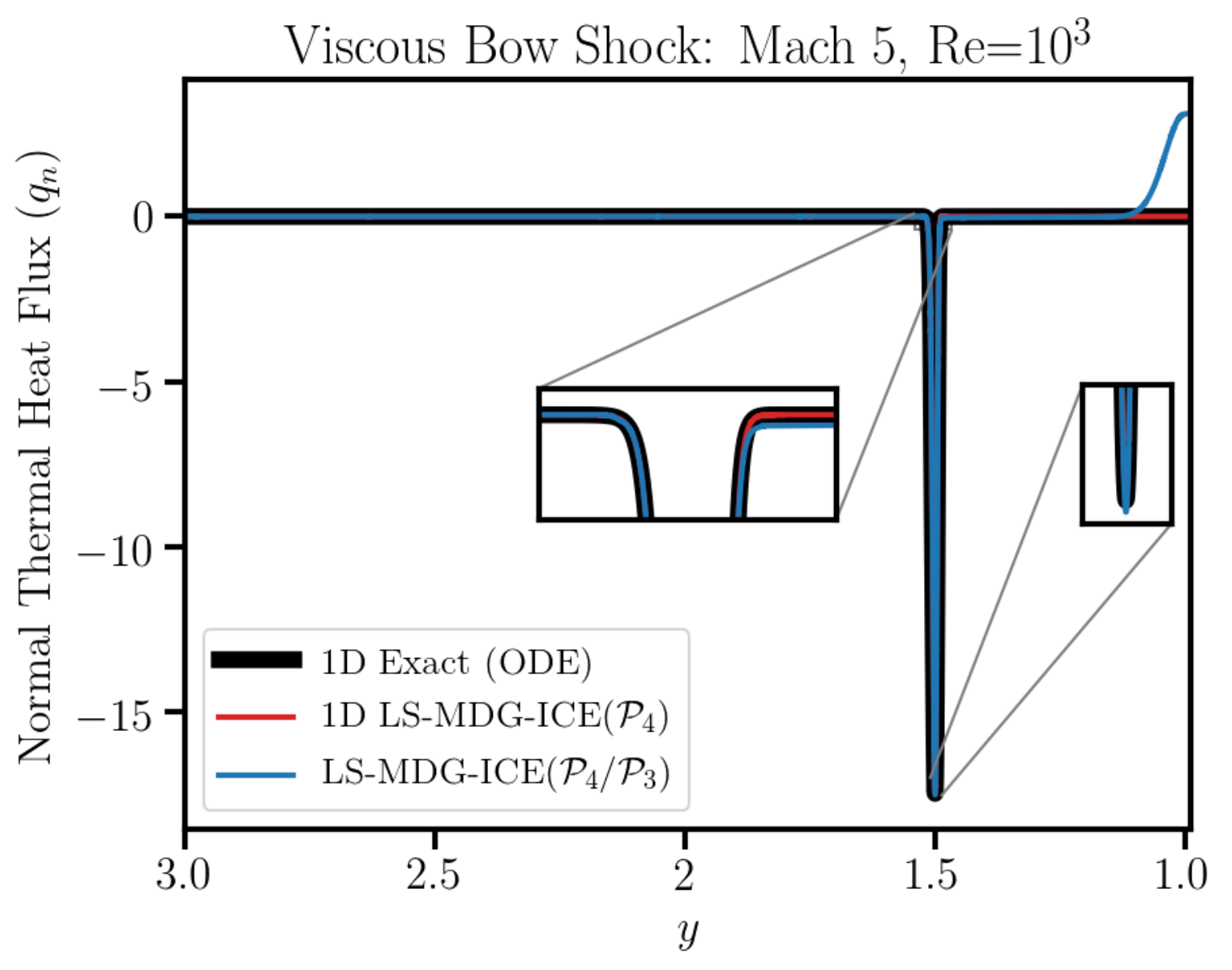}
\par\end{centering}
}

\caption{\label{fig:Bow-Shock-Re-1000-1D}Centerline profiles, sampled along
$x=0$, for the viscous Mach 5 bow shock at $\mathrm{Re}=10^{3}$
computed with LS-MDG-ICE($\mathcal{P}_{4}/\mathcal{P}_{3}$) compared
to ODE and MDG-ICE($\mathcal{P}_{4}$) approximations of the exact
solution for the corresponding one-dimensional viscous shock. The
one-dimensional LS-MDG-ICE($\mathcal{P}_{4}$) approximation was computed
using 16 isoparametric line cells. The location of the shock was computed
as as $y=1.5000575$ for a stand-off distance of $0.5000575$.}
\end{figure}
\begin{figure}
\begin{centering}
\subfloat[\label{fig:Viscous-Bow-Shock-Re-10000-Mesh}508 sub-parametric $\mathcal{P}_{4}/\mathcal{P}_{3}$
triangle cells]{\begin{centering}
\includegraphics[width=0.9\linewidth]{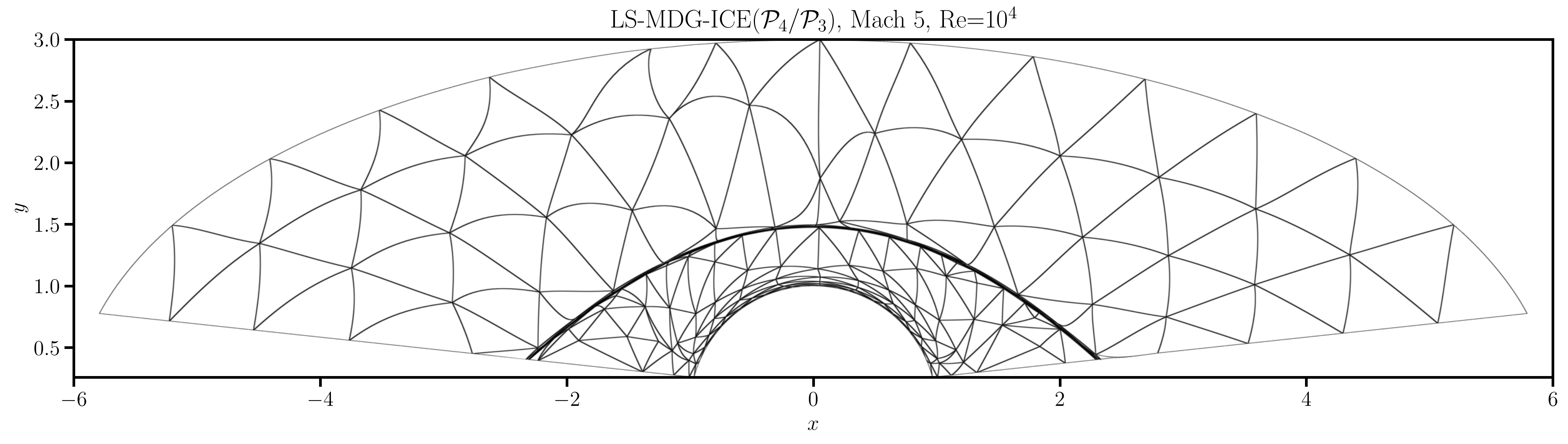}
\par\end{centering}
}
\par\end{centering}
\begin{centering}
\subfloat[\label{fig:Viscous-Bow-Shock-Re-10000-Temperature}508 sub-parametric
$\mathcal{P}_{4}/\mathcal{P}_{3}$ triangle cells]{\begin{centering}
\includegraphics[width=0.9\linewidth]{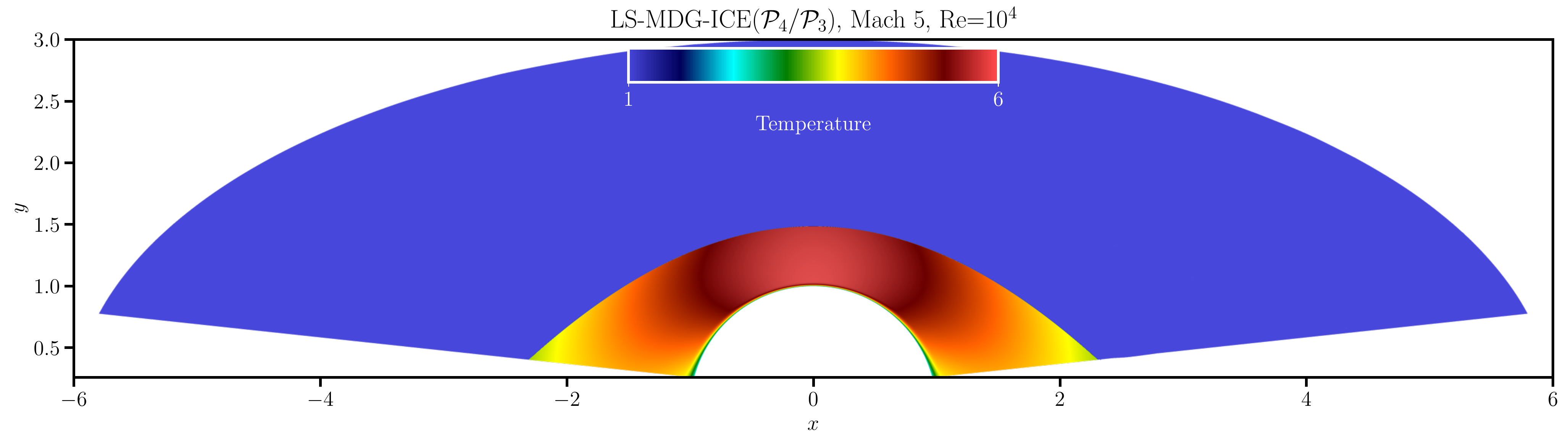}
\par\end{centering}
}
\par\end{centering}
\begin{centering}
\subfloat[\label{fig:Viscous-Bow-Shock-Re-10000-Mach}508 sub-parametric $\mathcal{P}_{4}/\mathcal{P}_{3}$
triangle cells]{\begin{centering}
\includegraphics[width=0.9\linewidth]{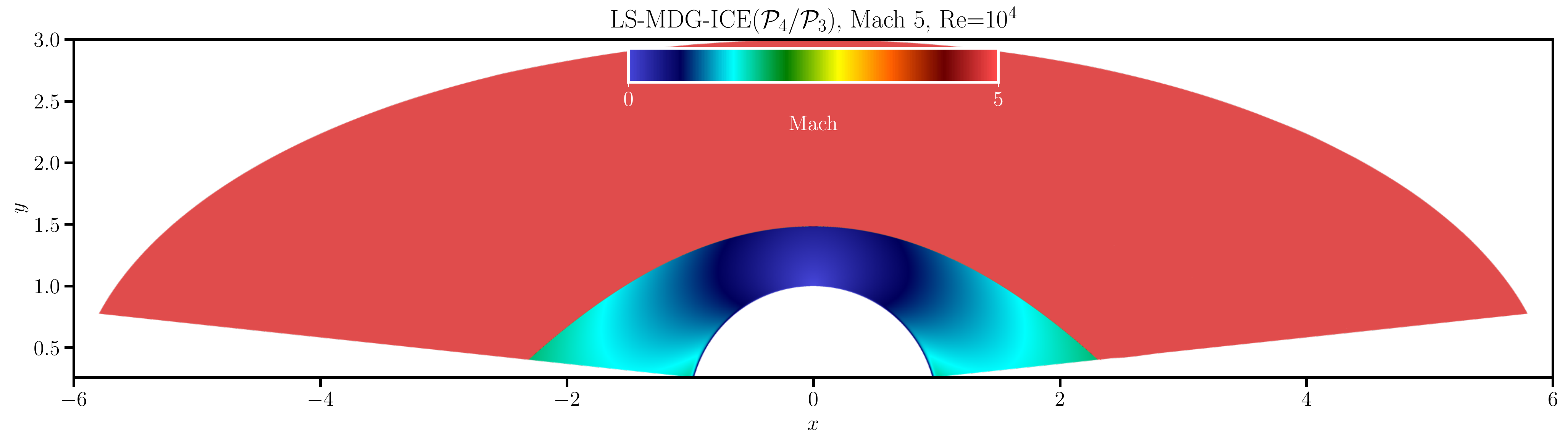}
\par\end{centering}
}
\par\end{centering}
\begin{centering}
\subfloat[\label{fig:Viscous-Bow-Shock-Re-10000-Pressure}508 sub-parametric
$\mathcal{P}_{4}/\mathcal{P}_{3}$ triangle cells]{\begin{centering}
\includegraphics[width=0.9\linewidth]{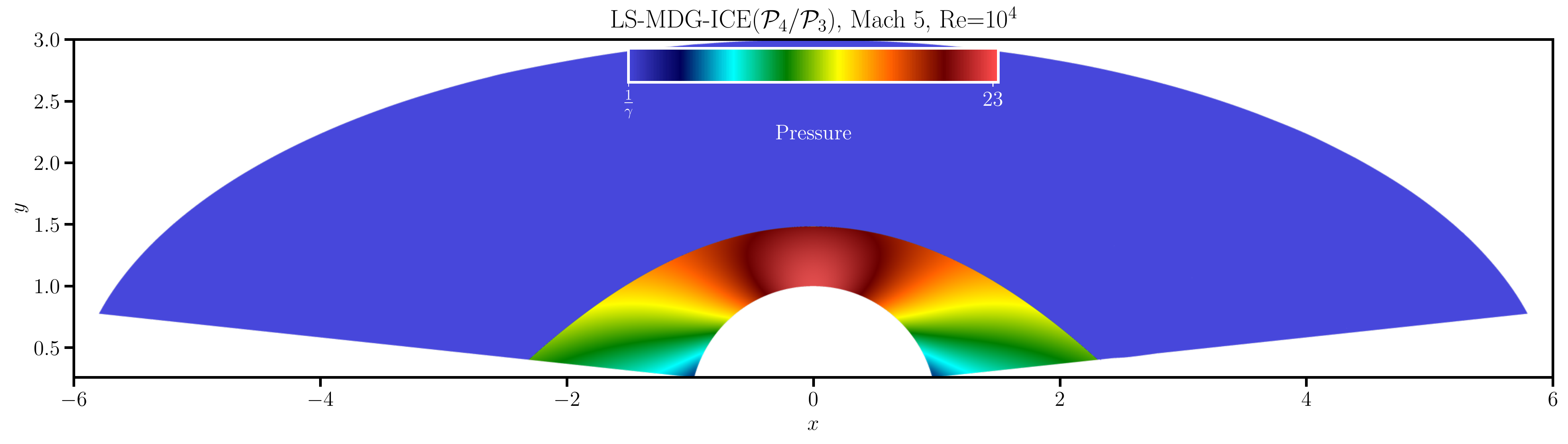}
\par\end{centering}
}
\par\end{centering}
\centering{}\caption{\label{fig:Viscous-Bow-Shock-Re-10000-MDG-ICE}The LS-MDG-ICE solution
computed using $\mathcal{P}_{4}/\mathcal{P}_{3}$ sub-parametric triangle
elements for the viscous Mach 5 bow shock at $\mathrm{Re}=10^{4}$.
The MDG-ICE solution was initialized from the LS-MDG-ICE solution
at $10^{3}$ $\mathrm{Re}$ shown in Figure~\ref{fig:Viscous-Bow-Shock-DG-Mesh}.
The location of the shock along the line $x=0$ was computed as $y=1.4847275$
for a stand-off distance of $0.4847275$.}
\end{figure}
\begin{figure}
\subfloat[\label{fig:Bow-Shock-1D-Re-10000-Temperature}The temperature sampled
along $x=0$. The exact temperature at the stagnation point, $T=2.5$,
is marked with the symbol $\times$. ]{\begin{centering}
\includegraphics[width=0.3\linewidth]{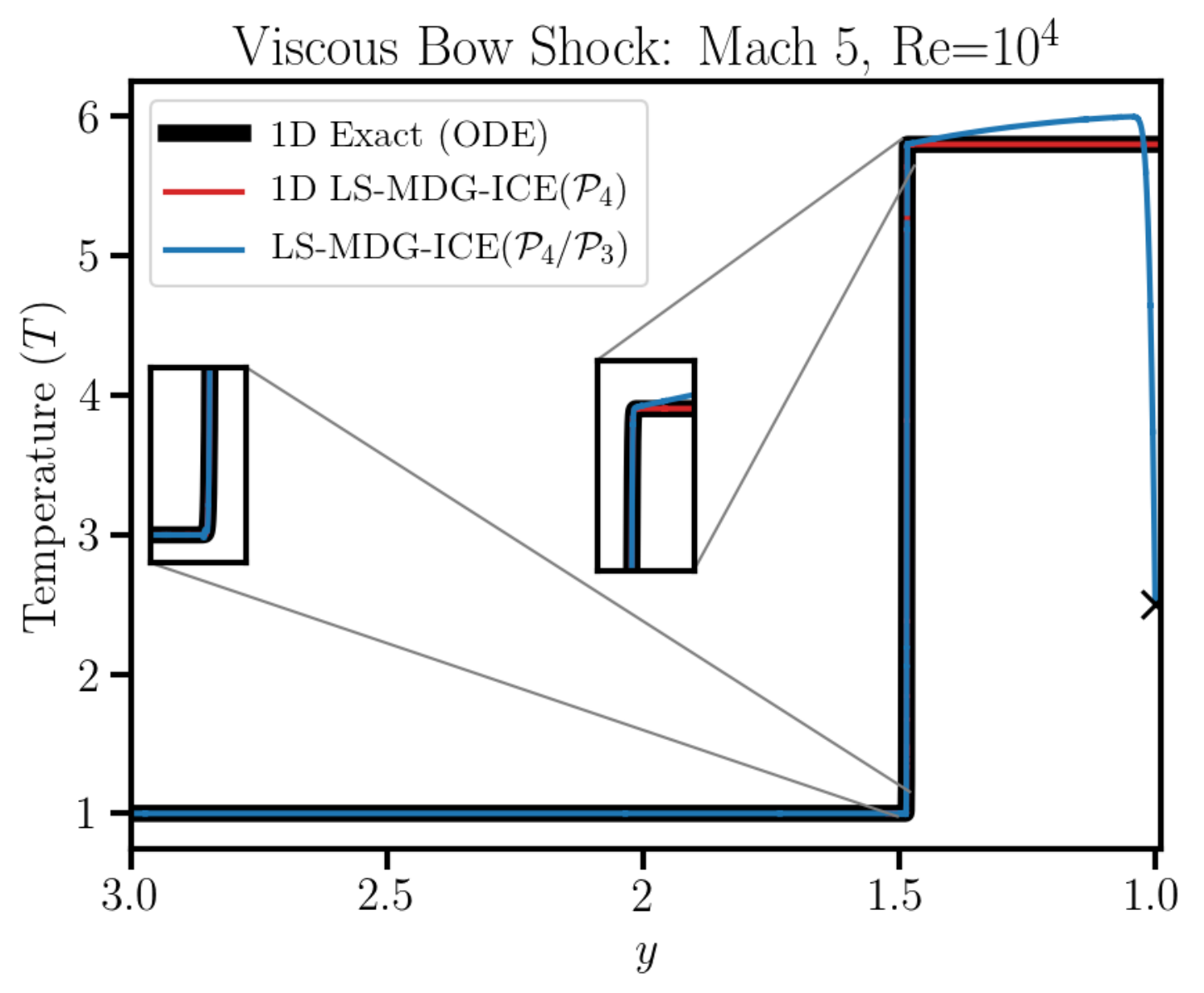}
\par\end{centering}
}\hfill{}\subfloat[\label{fig:Bow-Shock-1D-Re-10000-VelocityX}The normal velocity, $v_{n}=0$,
sampled along $x=0$. The exact normal velocity at the stagnation
point, $v_{n}=0$, is marked with the symbol $\times$. ]{\begin{centering}
\includegraphics[width=0.3\linewidth]{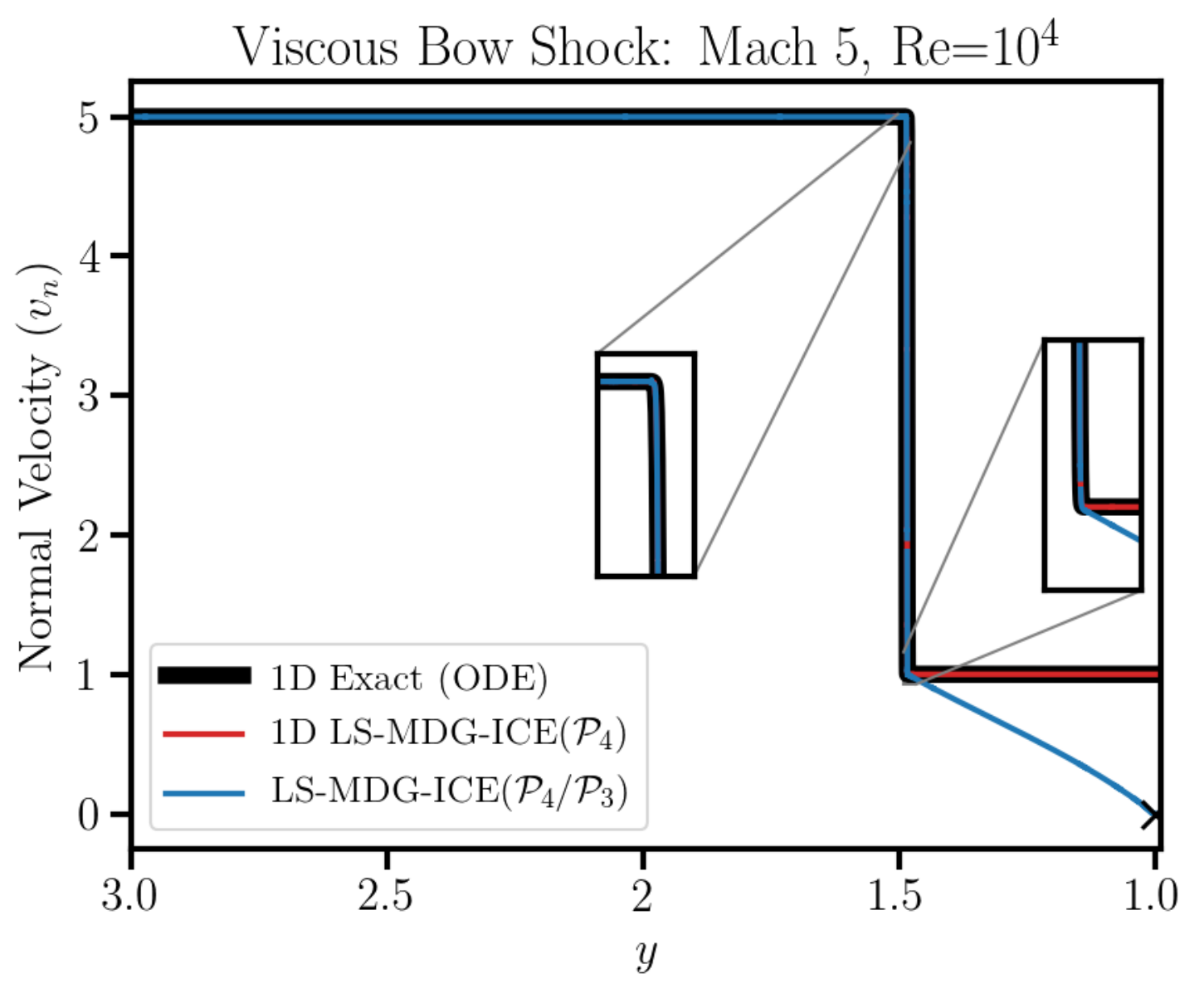}
\par\end{centering}
}\hfill{}\subfloat[\label{fig:Bow-Shock-1D-Re-10000-Pressure}The pressure, $p$, sampled
along $x=0$. The exact pressure at the stagnation point for an inviscid
flow, $p\approx23.324$, is marked with the symbol $\times$.]{\begin{centering}
\includegraphics[width=0.3\linewidth]{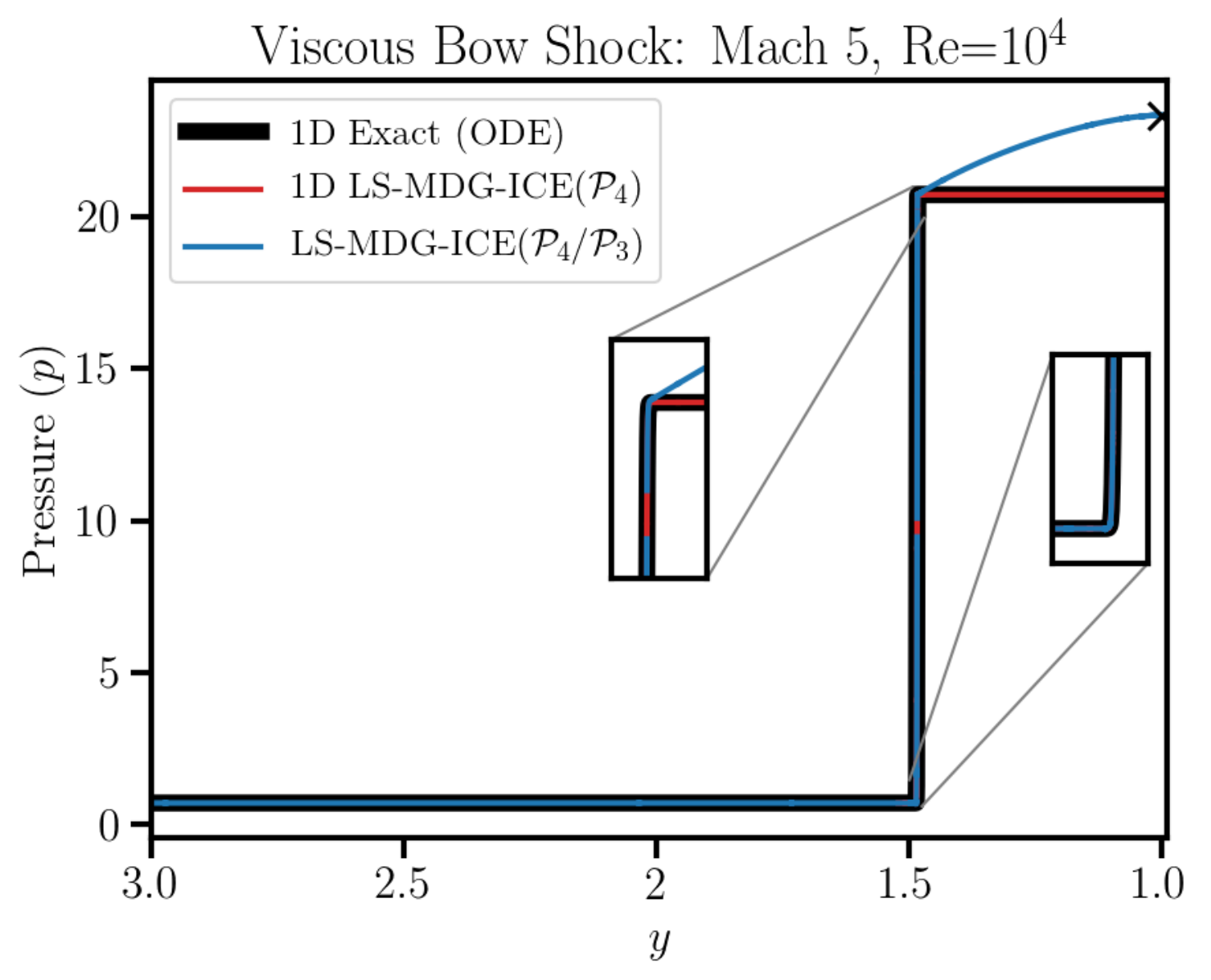}
\par\end{centering}
}

\subfloat[\label{fig:Bow-Shock-1D-Re-10000-Density}The density, $\rho$, sampled
along $x=0$. The density at the stagnation point, computed using
the stagnation pressure corresponding to an inviscid flow, $\rho\approx13.061$,
is marked with the symbol $\times$. ]{\begin{centering}
\includegraphics[width=0.3\linewidth]{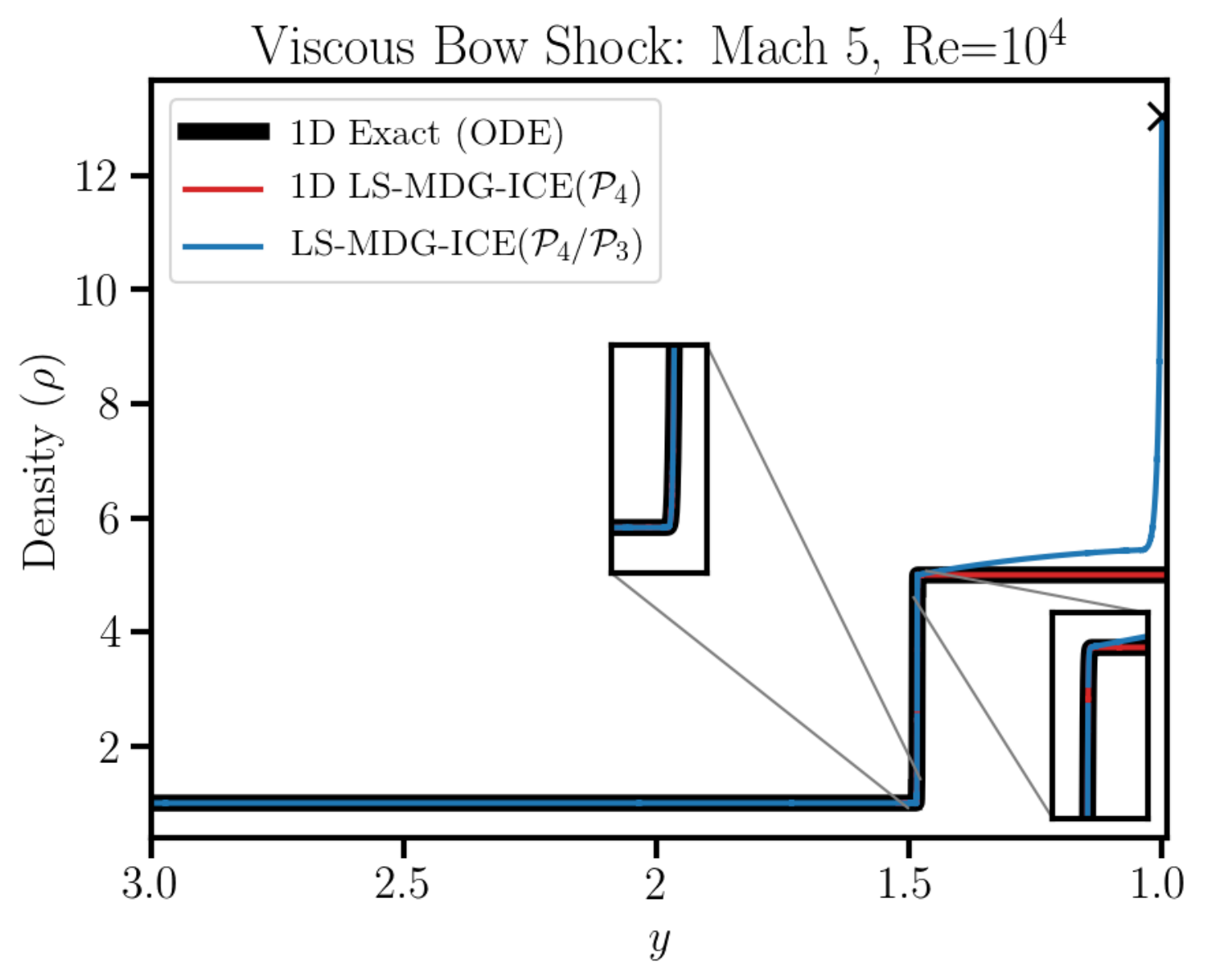}
\par\end{centering}
}\hfill{}\subfloat[\label{fig:Bow-Shock-1D-Re-10000-FluxXMomentumX}The normal component
of the normal viscous stress tensor, $\tau_{nn}$, sampled along $x=0$.
]{\begin{centering}
\includegraphics[width=0.3\linewidth]{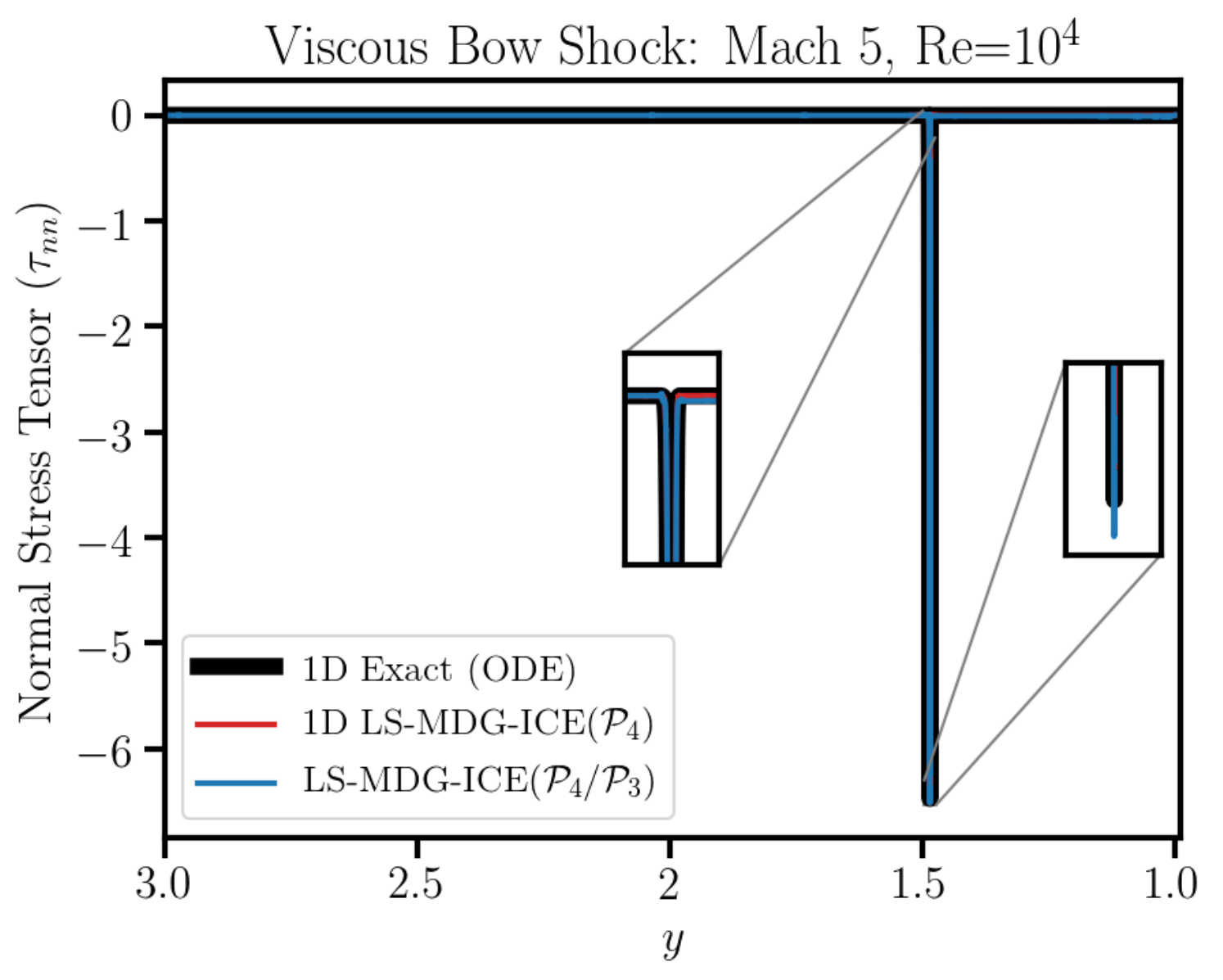}
\par\end{centering}
}\hfill{}\subfloat[\label{fig:Bow-Shock-1D-Re-10000-FluxXEnergyStagnationDensity}The
normal thermal heat flux, $q_{n}$, sampled along $x=0$. ]{\begin{centering}
\includegraphics[width=0.3\linewidth]{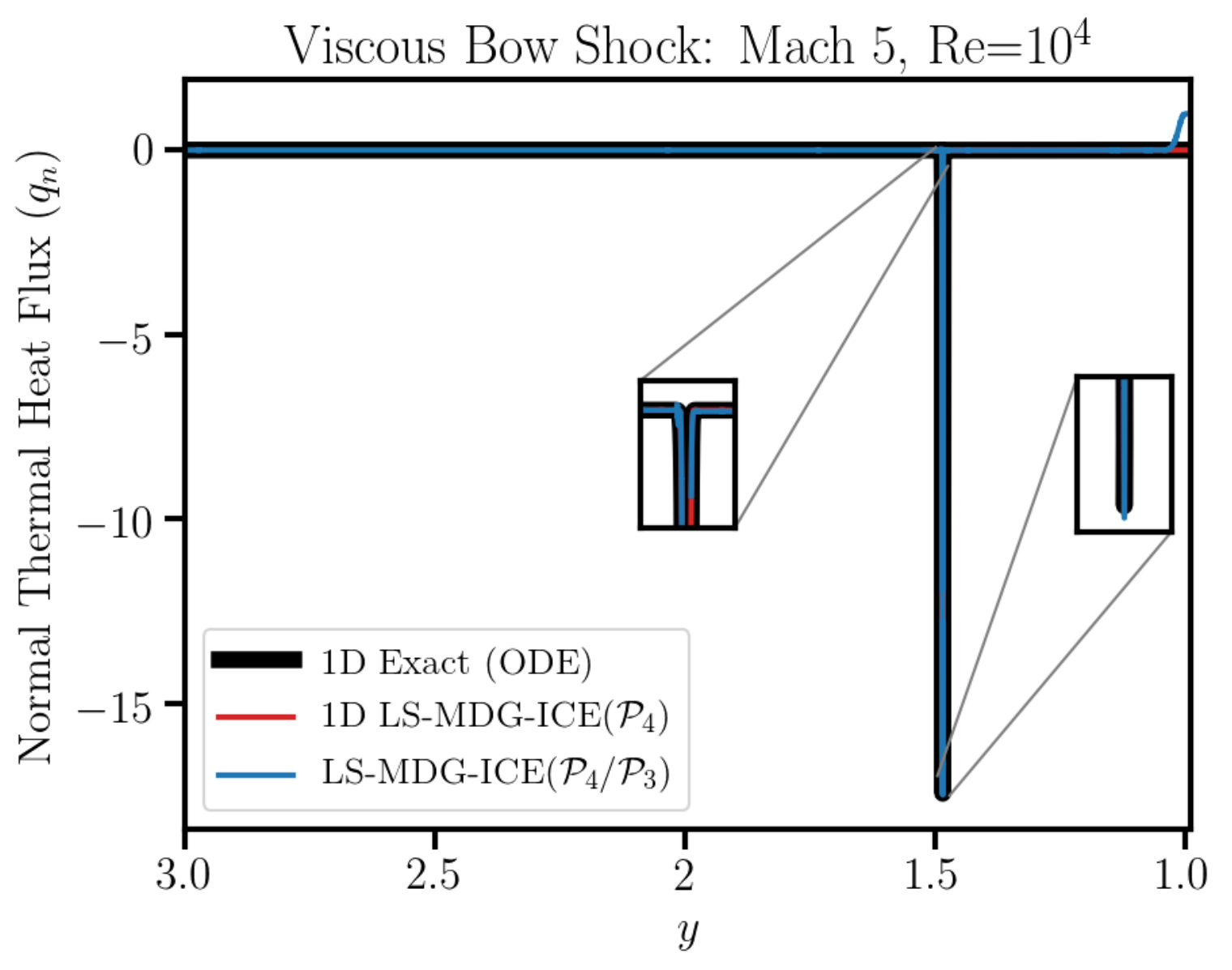}
\par\end{centering}
}

\caption{\label{fig:Bow-Shock-Re-10000-1D}Centerline profiles of temperature
and normal velocity for the viscous Mach 5 bow shock at $\mathrm{Re}=10^{4}$
computed with LS-MDG-ICE($\mathcal{P}_{4}/\mathcal{P}_{3}$) compared
to ODE and LS-MDG-ICE($\mathcal{P}_{4}$) approximations of the exact
solution for the corresponding one-dimensional viscous shock. The
one-dimensional LS-MDG-ICE($\mathcal{P}_{4}$) approximation was computed
using 16 isoparametric line cells. The location of the shock along
the line $x=0$ was computed as $y=1.4847275$ for a stand-off distance
of $0.4847275$.}
\end{figure}
\begin{figure}
\centering{}%
\begin{minipage}[c][1\totalheight][t]{0.45\columnwidth}%
\subfloat[\label{fig:Viscous-Bow-Shock-Temperature-Re-00001000-zoom}The final
grid and temperature fields of the LS-MDG-ICE$\left(\mathcal{P}_{4}\right)$
solution computed using 502 $\mathcal{P}_{4}$ isoparametric triangle
elements for the viscous Mach 5 bow shock at $\mathrm{Re}=10^{3}$
shown in Figure~\ref{fig:Viscous-Bow-Shock-Mesh} and Figure~\ref{fig:Viscous-Bow-Shock-Temperature}
respectively.]{\begin{centering}
\begin{tabular}{c}
\includegraphics[width=0.8\linewidth]{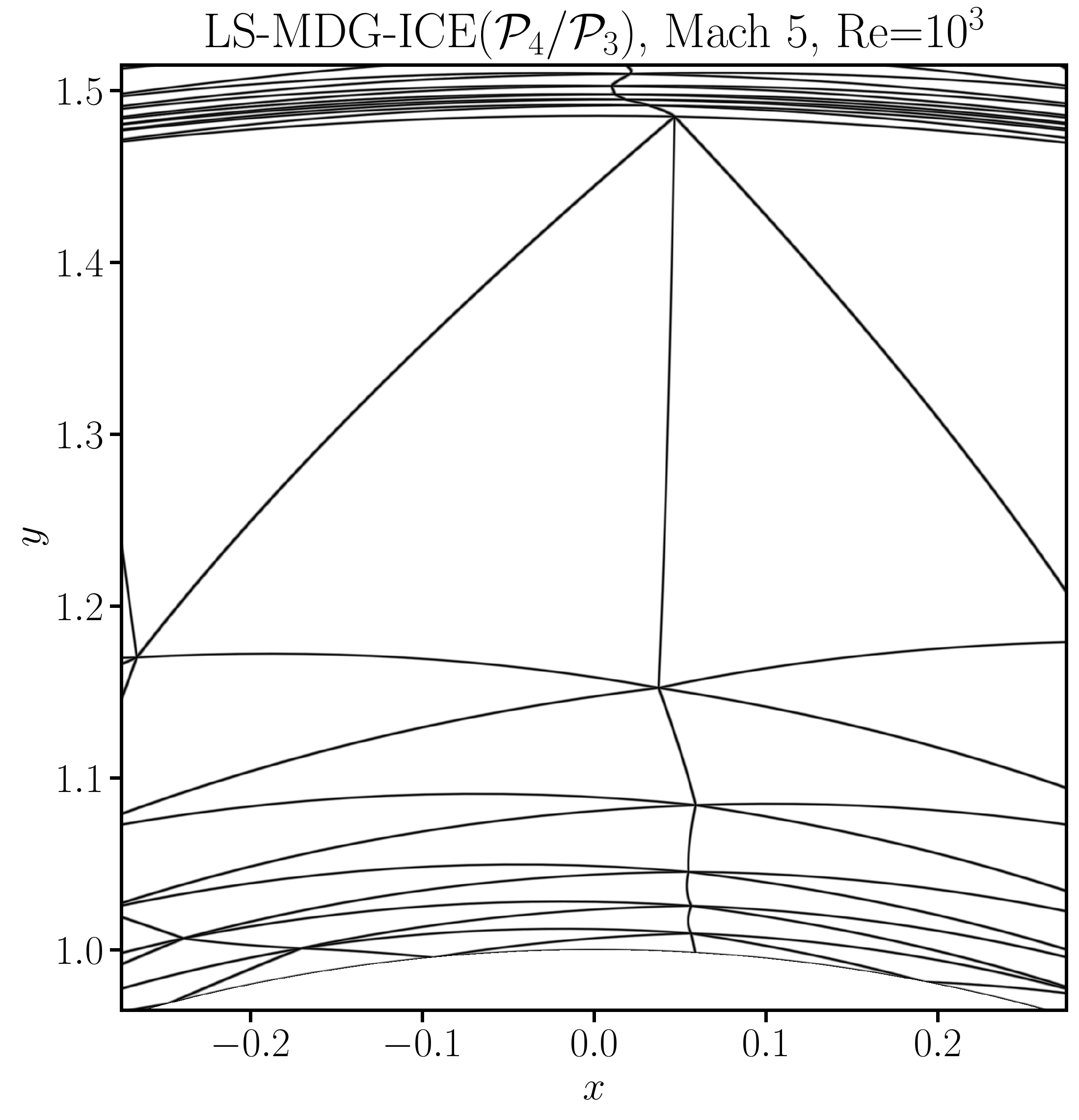}\tabularnewline
\includegraphics[width=0.8\linewidth]{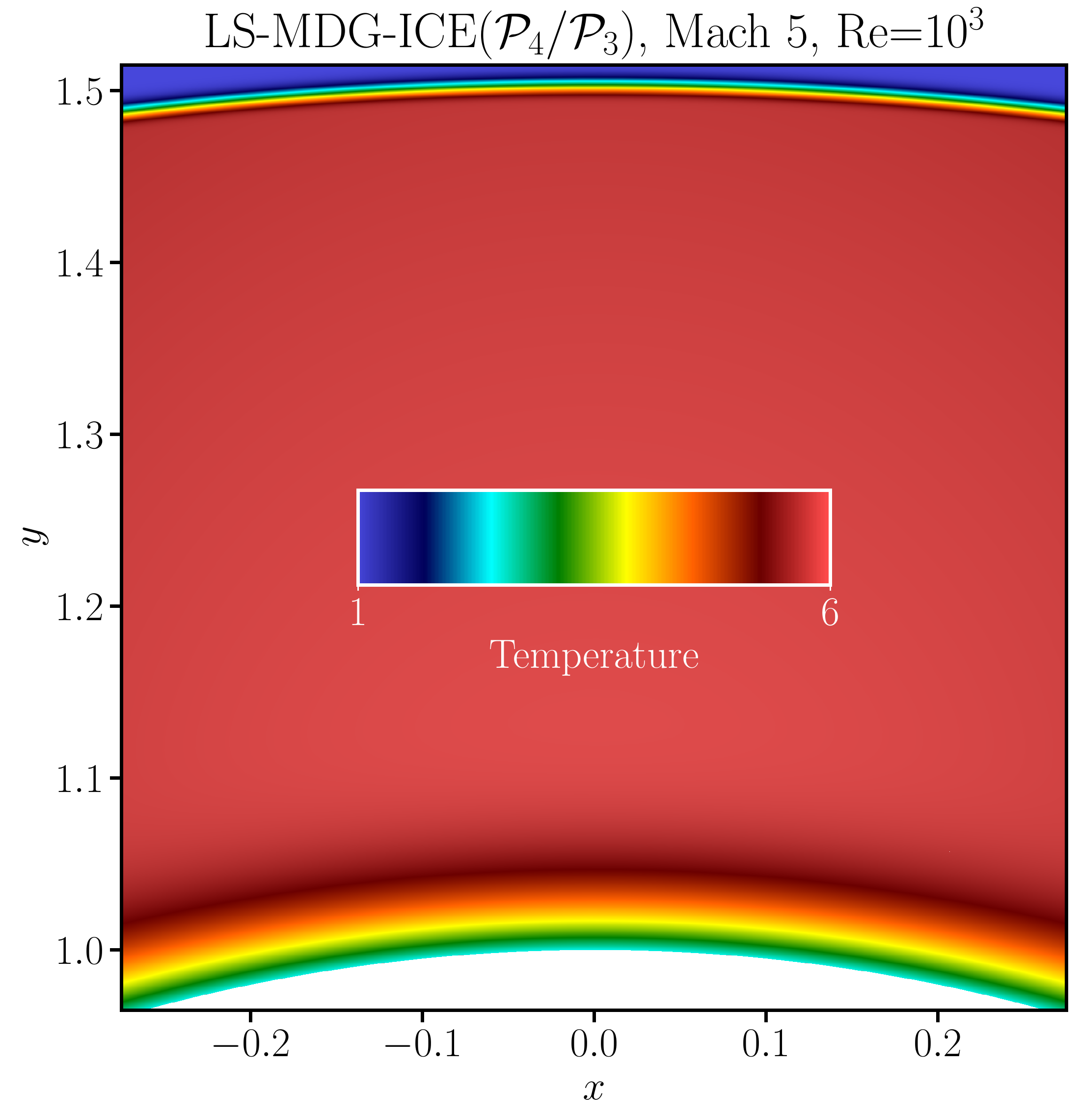}\tabularnewline
\end{tabular}
\par\end{centering}
}%
\end{minipage}%
\begin{minipage}[c][1\totalheight][t]{0.45\columnwidth}%
\subfloat[\label{fig:Viscous-Bow-Shock-Temperature-Re-00010000-zoom}The final
grid and temperature fields of the LS-MDG-ICE$\left(\mathcal{P}_{4}\right)$
solution computed using 508 $\mathcal{P}_{4}$ isoparametric triangle
elements for the viscous Mach 5 bow shock at $\mathrm{Re}=10^{4}$
shown in Figure~\ref{fig:Viscous-Bow-Shock-Re-10000-Mesh} and Figure~\ref{fig:Viscous-Bow-Shock-Re-10000-Temperature}.]{\begin{centering}
\begin{tabular}{c}
\includegraphics[width=0.8\linewidth]{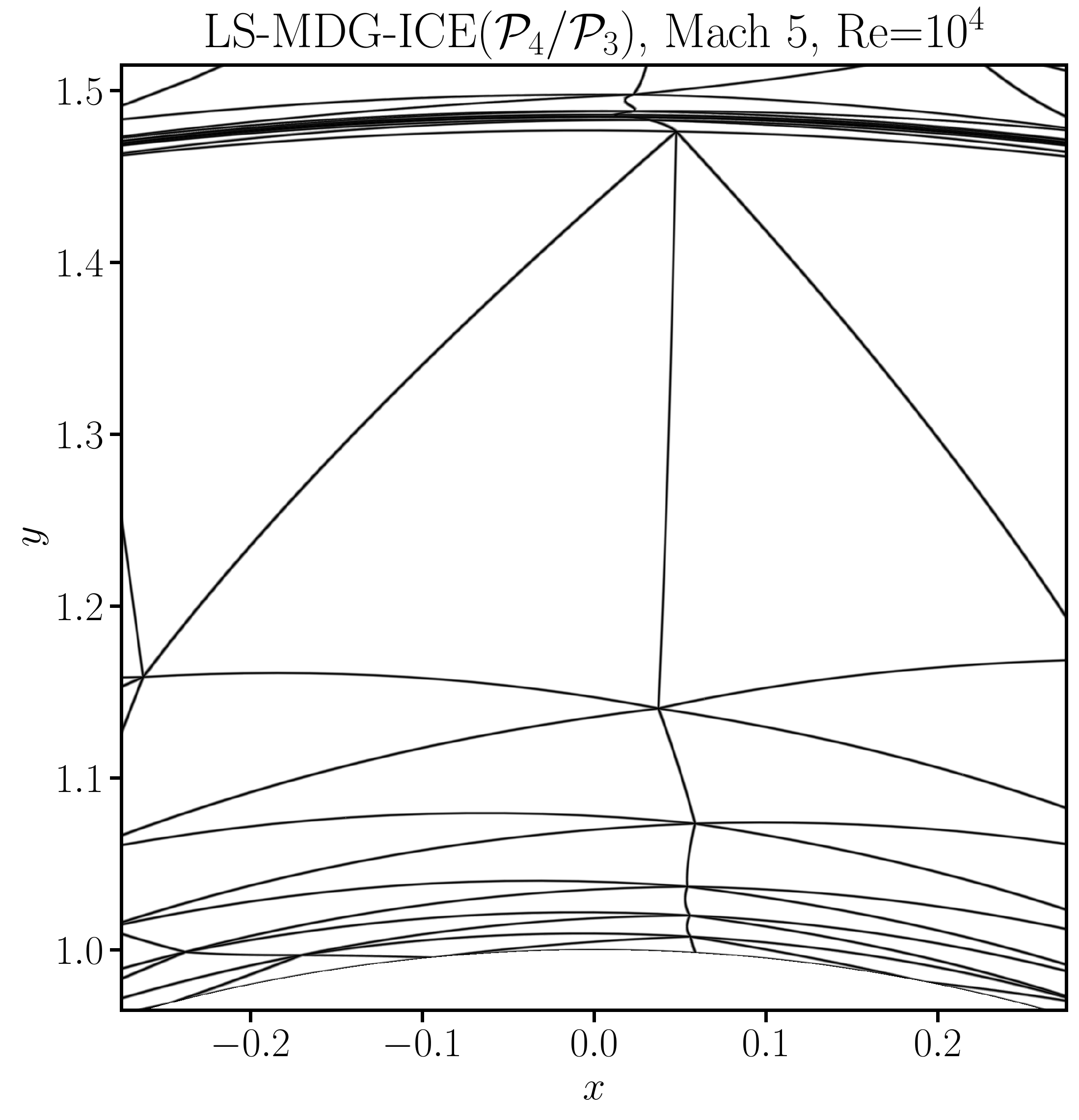}\tabularnewline
\includegraphics[width=0.8\linewidth]{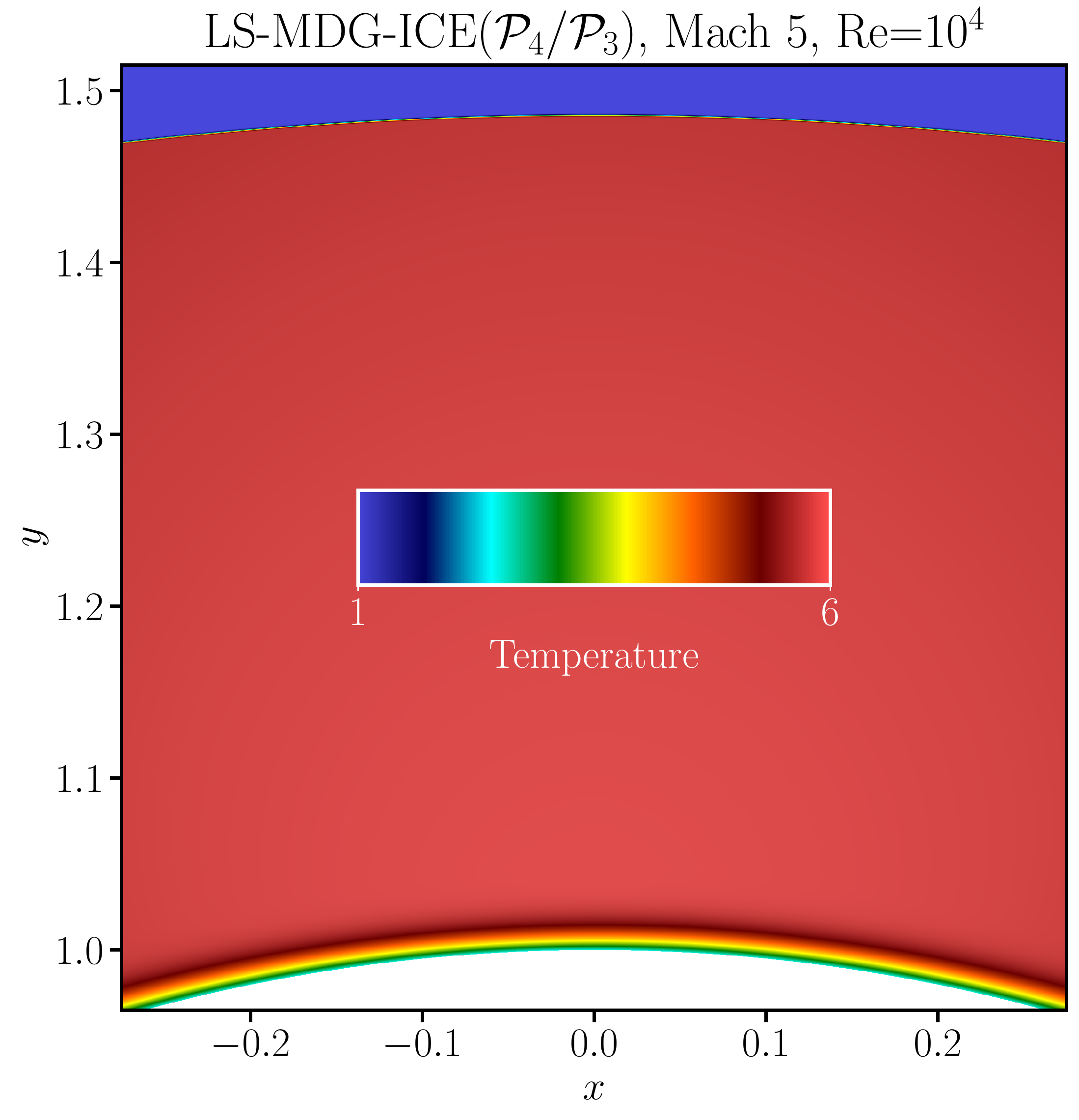}\tabularnewline
\end{tabular}
\par\end{centering}
}%
\end{minipage}\caption{\label{fig:Viscous-Bow-Shock-MDG-ICE-zoom}The final grid and temperature
fields corresponding to the LS-MDG-ICE solution computed using $\mathcal{P}_{4}$
isoparametric triangle elements for the viscous Mach 5 bow shock at
$\mathrm{Re}=10^{3}$ and $\mathrm{Re}=10^{4}$. Local edge refinement
was used to adaptively split highly anisotropic elements within the
viscous structures as they were resolved by LS-MDG-ICE.}
\end{figure}

The one-dimensional viscous shock solution, as described in Section~\ref{sec:CNS-Viscous-Shock},
is used to verify that the viscous LS-MDG-ICE formulation predicts
the correct viscous shock profile when diffusive effects are prominent.
Figure~\ref{fig:Bow-Shock-Re-1000-1D} and Figure~\ref{fig:Bow-Shock-Re-10000-1D}
compare an approximation of the exact solution for a one-dimensional
viscous shock to the centerline profiles of the Mach 5 bow shock at
$\mathrm{Re}=10^{3}$ and $\mathrm{Re}=10^{4}$, respectively, for
the following variables: temperature, $T$, normal velocity, $v_{n}$,
pressure, $p$, density, $\rho$, normal component of the normal viscous
stress tensor, $\tau_{nn}$, and normal heat flux, $q_{n}$, where
the normal is taken to be in the in the streamwise direction. For
reference and to aid comparisons to the exact solution, we also plotted
the LS-MDG-ICE($\mathcal{P}_{4}$) solution to the one-dimensional
viscous shock problem computed using 16 isoparametric line cells.
For reference, the exact and approximate values at the stagnation
point are marked on the centerline plots with the symbol $\times$.
When the exact value for the viscous flow was unavailable we approximate
it with the value corresponding to the inviscid solution, which is
not expected to differ significantly for the cases considered here,
see~\citep[Table 1]{Wil16}.

As noted in~\citep{Ker20}, the viscous shock profile corresponding
to the one-dimensional solution is expected to deviate from the bow
shock centerline profiles downstream of the viscous shock where the
viscous effects are produced by the boundary layer associated with
the blunt body centered at the origin. The jump across the viscous
shock is directly comparable to the one-dimensional solution in the
case of density, for which the viscous flux is zero. Figure~\ref{fig:Bow-Shock-1D-Density}
shows that LS-MDG-ICE accurately reproduces the density profile of
the exact shock structure with only a few anisotropic high-order cells.

The ability of LS-MDG-ICE($\mathcal{P}_{4}/\mathcal{P}_{3}$) to resolve
viscous layer using anisotropic curvilinear $r$-adaptivity is further
demonstrated in Figure~\ref{fig:Viscous-Bow-Shock-MDG-ICE-zoom}
where the solutions between the viscous shock and the blunt body for
$\mathrm{Re}=10^{3}$ and $\mathrm{Re}=10^{4}$ are compared. This
feature was showcased in~\citep[Section 3.3]{Ker20} and was further
explored in Section~\ref{sec:linear-advection-diffusion-1d} for
the case of one-dimensional advection-diffusion where it enabled LS-MDG-ICE
to achieve super-optimal convergence rates on the order of $2p$.
Here we again demonstrate it naturally extends to multiple dimensions
where LS-MDG-ICE automatically repositions the points in order to
resolve the flow field over the relevant physical length scales.

\section*{Conclusions and future work}

The Discontinuous Petrov-Galerkin (DPG) methodology of Demkowicz and
Gopalakrishnan was used to derive a least-squares formulation of MDG-ICE
with optimal test functions, which were systematically generated from
the trial spaces of both the discrete flow field and discrete geometry.
The LS-MDG-ICE formulation inherits unique and desirable features
of the (DLS) MDG-ICE formulation used in~\citep{Cor18,Ker20}, such
as the ability to fit a priori unknown interfaces and resolve sharp,
but smooth, viscous flow features using anisotropic curvilinear $r$-adaptivity.
Additionally, for both linear and nonlinear problems in one dimension,
the LS-MDG-ICE discretization was verified to converge with optimal
($p+1$) and super-optimal ($2p$) order under grid refinement for
polynomial spaces $\mathcal{P}_{2}$ through $\mathcal{P}_{5}$ using
a static, or fixed, discrete geometry and a variable discrete geometry.
We also observed optimal-order convergence for a steady, two-dimensional
linear advection problem. 

In two dimensions , the LS-MDG-ICE formulation is capable of producing
highly accurate solutions for problems involving curved interface
topologies and sharp viscous layers. For the Mach 3 inviscid bow,
LS-MDG-ICE accurately fit the curved shock and the computed stand-off
distance compared well to the value reported previously~\citep{Cor18}.
The LS-MDG-ICE($\mathcal{P}_{2}$) solution was accurate and oscillation-free
on a grid consisting of only 392 isoparametric triangle cells. Local
grid operations, i.e., edge refinement and collapse, were not required
to fit the curved interface. 

The first-order system LS-MDG-ICE formulation for viscous flows was
used to approximate a Mach 5 viscous Bow shock at $\mathrm{Re}=10^{3}$
and $\mathrm{Re}=10^{4}$. The one-dimensional analysis from our previous
work~\citep{Ker20} was reproduced for the case of $\mathrm{Re}=10^{3}$
and it was extended to the case of $\mathrm{Re}=10^{4}$. The centerline
profiles through the viscous shock were in agreement with the ODE
and one-dimensional LS-MDG-ICE($\mathcal{P}_{4}$) approximations
of the exact one-dimensional viscous shock profile. For the case of
$\mathrm{Re}=10^{3}$, the estimated stand off distance of $0.5000575$
compares well to the previously reported MDG-ICE standoff distance
of $0.49995$. For these challenging problems, we scaled the linear
elastic grid regularization by a factor similar to the one introduced
by Zahr et al.~\citep{Zah19} that is proportional to the inverse
of the element volume, which helped prevent cells from degenerating
as the solver simultaneously resolved the viscous shock and boundary
layer via anisotropic curvilinear $r$-adaptivity.

Future work will analyze the proposed method in a rigorous setting
and further study its ability to achieve optimal-order convergence
in the case of higher-dimensional, nonlinear problems. 

\section*{Acknowledgments}

This work is sponsored by the Office of Naval Research through the
Naval Research Laboratory 6.1 Computational Physics Task Area.

\section*{}

\bibliographystyle{elsarticle-num}
\bibliography{citations}

\end{document}